\documentclass[wcp]{jmlr}

\usepackage{longtable}

\usepackage{booktabs}

\usepackage{microtype}
\usepackage{graphicx}

\usepackage{cancel} 

\newcommand{\red}[1]{\textcolor{red}{#1}}
\newcommand{\blue}[1]{\textcolor{blue}{#1}}

\usepackage{amsmath}
\usepackage{amssymb}
\usepackage{mathtools}

\usepackage[capitalize,noabbrev]{cleveref}

\usepackage{comment}

\def\st{{\em s.t.~}}
\def\ie{{\em i.e.,~}}
\def\eg{{\em e.g.,~}}
\def\cf{{\em cf.,~}}

\newcommand{\wrt}{w.r.t.}

\newcommand{\tr}{tr}
\newcommand{\nn}{\nonumber} 
\newcommand{\R}{\mathbb{R}} 

\newcommand{\Indi}{\mathbb{I}} 
\newcommand{\Ns}{\mathbb{N}^\star} 
\newcommand{\pen}{\text{pen}} 
\newcommand{\bP}{\mathbb{P}} 

\newcommand{\h}{\hspace*{0.3cm}}

\DeclareMathOperator*{\argmin}{arg~min}

\DeclareMathOperator*{\id}{I} 
\newcommand{\zero}{\ensuremath{\mathbf{0}}}

\DeclareMathOperator*{\vect}{vec} 
\DeclareMathOperator*{\kl}{KL}
\DeclareMathOperator*{\tkl}{KL^{\otimes n}}
\DeclareMathOperator*{\jtkl}{JKL_{\rho}^{\otimes n}}

\DeclareMathOperator*{\card}{card}
\DeclareMathOperator*{\lasso}{Lasso}

\let\inf\relax 
\DeclareMathOperator*\inf{\vphantom{p}inf}

\DeclareMathOperator*{\adj}{Adj}

\newcommand{\norm}[1]{\left\lVert#1\right\rVert}

\DeclarePairedDelimiter\ceil{\lceil}{\rceil}

\newcommand{\E}[2]{\mathbb{E}_{#1} \left[#2\right]}

\newcommand{\vertiii}[1]{{\left\vert\kern-0.25ex\left\vert\kern-0.25ex\left\vert #1 
		\right\vert\kern-0.25ex\right\vert\kern-0.25ex\right\vert}}
\newcommand{\vertii}[1]{{\left\vert\kern-0.25ex\left\vert#1\right\vert\kern-0.25ex\right\vert}}	

\newcommand{\cA}{\mathcal{A}}

\newcommand{\cB}{\mathcal{B}}

\newcommand{\cF}{\mathcal{F}}

\newcommand{\cG}{\mathcal{G}}

\newcommand{\cM}{\mathcal{M}}

\newcommand{\cN}{\mathcal{N}}

\newcommand{\cT}{\mathcal{T}}

\newcommand{\cX}{\mathcal{X}}

\newcommand{\cY}{\mathcal{Y}}

\newcommand{\cZ}{\mathcal{Z}}

\makeatletter
\let\Ginclude@graphics\@org@Ginclude@graphics 
\makeatother

\jmlrvolume{Under Review}
\jmlryear{2024}

\title[Non-asymptotic oracle inequalities for the Lasso in high-dimensional MoE]{Non-asymptotic oracle inequalities for the Lasso in high-dimensional mixture of experts}

\author{\Name{TrungTin Nguyen} \Email{trungtin.nguyen@uq.edu.au}\\
	\addr School of Mathematics and Physics, The University of Queensland, St Lucia, QLD 4072, Australia; Univ. Grenoble Alpes, Inria, CNRS, Grenoble INP, LJK, Inria Grenoble Rhone-Alpes, 655 av. de l’Europe, 38335 Montbonnot, France.
	\AND
	\Name{Hien D Nguyen} \Email{h.nguyen5@latrobe.edu.au}\\
	\addr School of Computing, Engineering,
	and Mathematical Sciences, La Trobe University, Bundoora, VIC 3086, Australia;
	Institute of Mathematics for Industry, Kyushu University, Nishi Ward, Fukuoka 819-0395, Japan.
	\AND
	\Name{Faicel Chamroukhi}
	\Email{Faicel.chamroukhi@irt-systemx.fr}\\
	\addr IRT SystemX, Palaiseau, France.
	\AND 
	\Name{Geoffrey J McLachlan}
	\Email{g.mclachlan@uq.edu.au}\\
	\addr School of Mathematics and Physics, The University of Queensland, St Lucia, QLD 4072, Australia.
}

\begin{document}
	
	\maketitle
	
	\begin{abstract}
		We investigate the estimation properties of the mixture of experts (MoE) model in a high-dimensional setting, where the number of predictors is much larger than the sample size, and for which the literature is particularly lacking in theoretical results. We consider the class of softmax-gated Gaussian MoE (SGMoE) models, defined as MoE models with softmax gating functions and Gaussian experts, and focus on the theoretical properties of their $l_1$-regularized estimation via the Lasso. To the best of our knowledge, we are the first to investigate the $l_1$-regularization properties of SGMoE models from a non-asymptotic perspective, under the mildest assumptions, namely the boundedness of the parameter space. We provide a lower bound on the regularization parameter of the Lasso penalty that ensures non-asymptotic theoretical control of the Kullback--Leibler loss of the Lasso estimator for SGMoE models. Finally, we carry out a simulation study to empirically validate our theoretical findings.
		%
	\end{abstract}

	\section{Introduction}
	\subsection{Mixture of experts}

	MoE models, introduced in \citet{jacobs1991adaptive}, are a flexible mixture model construction for conditional density estimation and prediction. Because of their flexibility and the wealth of statistical estimation and model selection tools available, they have become widely used in statistics and machine learning.
	%
	%
	The MoE model construction allows the mixture weights (or gating functions) to depend on the explanatory variables (or predictors) together with the experts (or mixture component densities). This permits  the modeling of data arising from more complex data generating processes than those that can be analyzed using mixture models and mixture of regressions models, whose mixing parameters are independent of the covariates. 
	Finite mixture-type models have also become popular due to their universal approximation and good convergence rates for parameter and density estimation, which have been extensively studied in \citet{genovese2000rates, nguyen_convergence_2013, ho_convergence_2016, nguyen_approximation_2020,nguyen_approximation_2022}. In the same vein, recent results for parameter and conditional density estimation of MoE models have recently been published in \citet{jiang1999hierarchical, norets_approximation_2010, nguyen_universal_2016,nguyen2019approximation,nguyen2020approximationMoE,ho_convergence_2022,nguyen_demystifying_2023,nguyen_towards_2024,nguyen_general_2024}. 

	In the context of regression, \emph{softmax-gated Gaussian MoE models}, which will be referred to as \emph{SGMoE}, defined as MoE models with Gaussian experts and softmax gating functions, are a standard choice and a powerful tool for modeling more complex nonlinear relationships between response and predictor, arising from different subpopulations.
	Since each mixture weight is modeled by a softmax function of the covariates, the dependence on each feature appears both in the experts and in the gating functions, which allows one to capture more complex nonlinear relationships between the response and predictors arising from different subpopulations, compared to mixture of regressions models. This is demonstrated via numerical experiments in several works such as \citet{chamroukhi2019regularizedIJCNN,chamroukhi2019regularized,montuelle_mixture_2014}.
	The reader is referred to \citet{yuksel2012twenty,nguyen2018practical} for reviews on this topic. 
	Statistical estimation and variable selection for MoE models in the high-dimensional regression setting remain challenging. In particular, from a theoretical point of view, there is a lack of results for MoE models, where the number of explanatory variables can be much larger than the sample size.
	In such situations, we need to reduce the dimension of the problem by looking for the most relevant relationships to avoid numerical problems while ensuring identifiability.

	\subsection{Related literature}
	We focus on the use of the Lasso, originally introduced by \citet{tibshirani_regression_1996}, also known as an $l_1$-penalised maximum likelihood estimator ($l_1$-PMLE). Using $l_1$-PMLE tends to produce sparse solutions and can be viewed as a convex surrogate for the non-convex $l_0$-penalization problem. Relaxation methods have attractive computational and theoretical properties (\cf~\citealp{fan2001variable}). First introduced for the linear regression model, the Lasso estimator has since been studied and extended to many statistical problems. 
	To deal with heterogeneous high-dimensional data, several researchers have studied the Lasso for variable selection in the context of mixture of regression models, see, \eg~ \citet{khalili2007variable,stadler2010l1,Meynet:2013aa,devijver2015l1}; and \citet{lloyd2018globally}. 
	In particular, \citet{stadler2010l1} provided an $l_0$-oracle inequality, satisfied by the Lasso estimators, conditional on the restricted eigenvalue condition, namely that the Fisher information matrix is positive definite.  Furthermore, they have to introduce some margin conditions to link the Kullback--Leibler (KL) loss function to the $l_2$-norm of the parameters. 
	Another direction of studying this problem is to look at its $l_1$-regularisation properties; see \eg~\citet{massart_lasso_2011,Meynet:2013aa,devijver2015l1}. As indicated by \citet{devijver2015l1}, in contrast to results for the $l_0$-penalty, some results for the $l_1$-penalty are valid without any assumptions, either on the Gram matrix or on the bound. 

	\subsection{Main contributions}
	Our overall contributions in the paper can be summarized as follows:
	\begin{enumerate}
		\item To the best of our knowledge, \emph{we are the first to study the $l_1$-regularization properties of the SGMoE models from a non-asymptotic point of view with the mildest assumptions}.  Theorem~\ref{thm_l1_Oracle_Inequality} provides a lower bound on the regularization parameter of the Lasso penalty that ensures non-asymptotic theoretical control of the KL loss of the 
		$l_1$-PMLE estimator for SGMoE models.
		Because this result is non-asymptotic, it is valid when $n$ is fixed, while the number of predictors $p$ can grow, with respect to $n$, and can be much larger than $n$. 
		%
		
		\item Our non-asymptotic result \emph{complements the standard asymptotic results} for \emph{high-dimensional SGMoE models for feature selection} using Lasso-PMLE or more general PMLE via the Scad penalty function of \citet{khalili2010new}.  
		Specifically, 
		\citet{khalili2010new} proved both consistency in feature selection and $\sqrt{n}$ consistency of the PMLE in SGMoE models, but under several strict conditions on the regularity of the true joint density function and on the choice of tuning parameters. 
		%
		On the contrary, the \emph{only mild assumption} we use here to  obtain the order of rate convergence of \emph{the error upper bounds} in \eqref{eq_l1_Oracle_Inequality} from Theorem~\ref{thm_l1_Oracle_Inequality} is \emph{boundedness on the parameter space}, which also appeared in \citet{khalili2010new,stadler2010l1,Meynet:2013aa,devijver2015l1}.

		\item We extend non-asymptotic results for mixture of regressions models \citep{massart_lasso_2011,Meynet:2013aa,devijver2015l1} to the more general SGMoE models as defined in \eqref{eq_defMoE}, for which the \emph{theoretical analysis of the non-asymptotic result is challenging}, because, in SGMoE models,  the dependence on each feature appears both in the means of the experts and in the gating functions. 
		This  requires in particular \emph{non-trivial technical proof} we establish in this paper.
		%
		%
		
		\item Our focus in this paper is on a simplified but standard setting in which the  expert component means are linear functions with respect to the explanatory variables. Despite this linear simplification, the overall {\it SGMoE model captures  the non-linearity of the true regression function} thanks to its mixture construction. 
		We believe that \emph{the general techniques we develop here can be extended to more general experts}, such as Gaussian experts with polynomial means (see, \eg~\citealp{mendes2012convergence}), hierarchical MoE for exponential family regression models \citep{jiang1999hierarchical}, and when the covariance matrix is also parameterized as an expert function potentially depending on the covariates as in \citet{ho_convergence_2022}. 
	\end{enumerate}

	\noindent{\bf Notations.} Throughout this paper, 
	$\{1, \dots, n\}$ is abbreviated as $[n]$ for $n \in \Ns=\{1,2,\dots\}$. 
	Here, $\vect(\cdot)$ is the vectorization operator that stacks the columns of a matrix into a vector.  We denote the induced $p$-norm of a matrix $\beta$ by $\norm{\beta}_p$, $p \in \left\{1,2,\infty\right\}$,  which differs from the vector norm $\norm{\vect(\beta)}_p$. 
	For a matrix $\Sigma$, $m(\Sigma)$ and $M(\Sigma)$ denote the smallest and largest eigenvalues of $\Sigma$, respectively.
	We write $\cN\left(\cdot;v,\Sigma\right)$ for the multivariate Gaussian density with mean $v$ and with covariance matrix $\Sigma$.
	%
	%
	Given an arbitrary event $\cT$ in some probability space, we define an indicator function by: ${\Indi}_{\cT}(\omega) = 1$ if $\omega\in \cT$ and ${\Indi}_{\cT}(\omega) = 0$ if $\omega \notin \cT$.
	
	\noindent {\bf Paper organization.} In Section \ref{notationAndFramework}, we discuss the construction and framework of high-dimensional SGMoE models. In Section \ref{mainResult}, we present our main result.
	Then, we conduct a simulation study to empirically verify our theoretical results in \ref{sec_numerical_experiment}.
	Some conclusions are given in Section \ref{sec_conclusions_perspectives}. Supplementary material is devoted to proving the technical results.

	\section{Problem setup}\label{notationAndFramework}
	\subsection{High-dimensional SGMoE models}
	\label{randomCovModelSelection}
	%
	In the high-dimensional regression setting, we observe $n$ couples $\left(x_{[n]},y_{[n]}\right) \equiv \left(x_i,y_i\right)_{i \in [n]}\in \left(\cX\times\cY\right)^n \subset \left(\R^p\times\R^q\right)^n$, where typically $p\gg n$, $x_i$ is fixed and $y_i$ is a realization of the random variable $Y_i$, $i\in[n]$.
	We assume that, conditional on $x_{[n]}$, $Y_{[n]}$ are independent and identically distributed (IID) with conditional PDF $s_0\left(\cdot| x_i\right)$. Our goal is to estimate $s_0$ from the observations using the following $K$-component SGMoE models:
	\begin{align} \label{eq_defMoE}
		s_{\psi}(y| x) = \sum_{k=1}^K \frac{\exp\left(\gamma_{k0} + \gamma_k^\top x\right)}{\sum_{l=1}^K\exp\left(\gamma_{l0} + \gamma_l^\top x\right)} \cN\left(y;\beta_{k0}+ \beta_k x,\Sigma_k\right), 
	\end{align}
	with $K \in \Ns$ and unknown parameters $\psi =(\gamma,\beta,\Sigma) \equiv \left(\gamma_{k0},\gamma_k,\beta_{k0},\beta_k,\Sigma_k \right)_{k \in [K]}$ 
	in a parameter space $\Psi$.
	%
	%
	%
	For technical reasons, we require that the covariates are fixed and the boundedness assumptions on the parameter space $\Psi$.
	%
	
	The explanatory variables $x_{[n]}$ and the number of components $K$ are both fixed.
	We assume that $\cX$ is a compact subset of $\R^p$ and the observations $x_{[n]}$ are finite.  Without loss of generality, we choose to rescale $x$, so that $\norm{x}_\infty \le 1$. Therefore, we can assume that $\cX = [0,1]^p$.
	However, the arguments in our proofs are valid for covariates of on scale. 
	We assume that there exists positive constants $A_{\gamma},A_{\beta},a_{\Sigma},A_{\Sigma}$, such that $\psi \in \widetilde{\Psi}$, where
	\begin{align}\label{eq_defBoundedParameters}
		\widetilde{\Psi} &= \Big\{\psi \in \Psi \mid \max_{k \in [K]}\sup_{x\in \cX} \left(\left|\gamma_{k0}\right|+\left|\gamma_k^\top  x\right|\right) \le A_{\gamma}, \h
		\max_{z\in [q]} \max_{k \in [K]}\sup_{x \in \cX} \left(\left|\left[\beta_{k0}\right]_z\right|+\left|\left[\beta_kx\right]_z\right|\right)  \le A_{\beta}, \nn\\
		& \quad \quad a_{\Sigma} \le m\left({\Sigma}^{-1}_{k}\right) \le  M\left({\Sigma}^{-1}_{k}\right) \le A_{\Sigma}\Big\}.
	\end{align}
	
	
	\noindent {\bf Collection of SGMoE models.} In summary, given \eqref{eq_defMoE} and \eqref{eq_defBoundedParameters}, we wish to estimate $s_0$ via  the following collection of SGMoE models:
	\begin{align} \label{eq_Sm_bounded}
		S &= \left\{(x,y) \mapsto
		s_{\psi}(y| x) \mid \psi \in \widetilde{\Psi} \right\}.
	\end{align}
	In particular, to simplify the proofs, we shall assume that the true conditional PDF $s_0$ belongs to $S$. That is to say, there exists $\psi_{0} = \left(\gamma_0,\beta_0,\Sigma_0\right) \in \widetilde{\Psi}$, such that $s_0 = s_{\psi_0}$. 
	From hereon in, where there is no confusion, we will use $s_0$ and $s_{\psi_0}$, interchangeably.
	%
	

\subsection{Minimum contrast estimation}
Several loss functions have been introduced into the MLE for MoE models. For example, by using the identifiability conditions, \citet{nguyen_demystifying_2023,nguyen_towards_2024,nguyen_general_2024} established inverse bounds between Hellinger distance and some Wasserstein distances or Voronoi loss functions to accurately capture heterogeneous parameter estimation convergence rates for different classes of MoE models.
However, in this paper, our main idea is to consider the Lasso as an $l_1$-ball model selection procedure, see \eg~\citet{massart_lasso_2011}. Therefore, we follow the framework of \emph{minimum contrast estimation}, see \eg~\citet[Chapter 1]{massart_concentration_2007}, \citet{arlot_survey_2010}, and \citet{barron_risk_1999}.
In this situation, negative log-likelihood (NLL) and KL divergence are the natural choices for the density estimation problem.

\noindent
{\bf Average KL divergence.} 
To take into account the structure of conditional PDFs, for fixed explanatory variables $\left(x_i\right)_{1 \le i \le n}$, we consider the following average KL loss function:
\begin{align}\label{Extend.CarolineMeynet.eq.2.2}
	{\kl}_n(s,t) &=\frac{1}{n} \sum_{i=1}^n \kl\left(s\left(\cdot| x_i\right),t\left(\cdot| x_i\right)\right), \text{for any densities $s$ and $t$, where}\\
	\kl(s,t) &= \begin{cases}
		\int_{\R^q} \ln\left(\frac{s(y)}{t(y)}\right)s(y)dy,& \text{ if $sdy$ is absolutely continuous w.r.t
			$tdy$},\\
		+\infty, & \text{ otherwise}.\\
	\end{cases}
\end{align}

\noindent {\bf Lasso estimator.} Conditioned on $\left(x_i\right)_{1 \le i \le n}$, the MLE approach suggests estimating $s_0$ by the conditional PDF $s_\psi$ 
that minimizes the NLL:
$-\frac{1}{n}\sum_{i=1}^n \ln\left(s_\psi\left(y_i | x_i\right)\right)$. 
%
However, in high-dimensional data, we need to regularize the MLE in order to obtain reasonable estimates. Here, we first consider the $l_1$-PMLE (the Lasso estimator):
%
\begin{align}\label{eq_lasso_estimator}
	\widehat{s}^{\lasso}_{\lambda} = \argmin_{s_\psi\in S} \Big\{& -\frac{1}{n}\sum_{i=1}^n \ln\left(s_\psi\left(y_i | x_i\right)\right) + \lambda \left(\norm{\gamma}_1+\norm{\vect(\beta)}_1\right)\Big\},
\end{align}
where $\lambda \geq 0$ is a regularization parameter to be tuned, $\norm{\gamma}_1 = \sum_{k=1}^K \sum_{j=1}^p\left|\gamma_{kj}\right|$, and $\norm{\vect(\beta)}_1 = \sum_{k=1}^K  \sum_{j=1}^p \sum_{z=1}^q\left|\left[\beta_k\right]_{z,j}\right|$.
It is worth noting that these two entry-wise $l_1$ norms do not contain scalar $\gamma_{k0}$ and vector $\beta_{k0}$ bias.
These Lasso regularisation terms encourage sparsity for both gating and expert parameters. 
\section{Main result} \label{mainResult}
%
%
To simplify the statement of  Theorem~\ref{thm_l1_Oracle_Inequality}, given some constants $\kappa \geq 148$, we first define the following condition for $\lambda$ and the constant $C_{1n}$ that appears on the upper risk bounds:
\begin{align}
	&\lambda \ge \kappa \frac{K}{\sqrt{n}} C_{0n}, \text{ }C_{0n}  =  B_{0n}\left(q\ln n \sqrt{\ln(2p+1)}+1\right), 	\label{eq_lowerbound_kappa}\\
	&B_{0n} = \max\left(A_\Sigma,1+KA_G\right)\left[1+2q\sqrt{q} A_\Sigma \left( 5A^2_\beta + 4 A_\Sigma \ln n\right) \right],\nn\\
	& C_{1n} = \sqrt{2 q A_\gamma}\left(\frac{ e^{q/2-1}\pi^{q/2} }{A_\Sigma ^{q/2} } + H_{s_0}\right) + B_{0n} C_{2n}, \label{eq_l1_oracle_error_upperbound} \\
	&H_{s_0} =\max\left\{0,\ln  \left[\left(4\pi\right)^{-q/2} A_\Sigma^{q/2}\right]  \right\}, C_{2n} =  302q  K\left[1+\left(A_\gamma+  q A_\beta + \frac{ q\sqrt{q}}{a_{\Sigma}}\right)^2\right]. \label{eq_l1_oracle_smallY_error_upperbound}
\end{align}
{\bf Remark.} 
In \eqref{eq_lowerbound_kappa}, we have taken care to make dependencies explicit, not only on the tuning constant $\kappa$, but also on $n$, $p$, $q$, and $K$ as well as on $A_\beta, A_\Sigma, A_G$---all of the quantities that constrain the parameters of the model. 
See Section \ref{sec_discussion_comparison} for a more detailed description.
Note that both $C_{0n}$ and $C_{1n}$ depend on the sample size $n$ only via the $\ln n$ term. Furthermore, $H_{s_0}$ is related to the negative of the differential entropy of the true unknown conditional density $s_0 \in S$; see the supplementary material for more details.
We state our main contribution: an $l_1$-oracle inequality for the Lasso estimator for SGMoE models via Theorem~\ref{thm_l1_Oracle_Inequality}.
%
\begin{theorem}[$l_1$-oracle inequality] \label{thm_l1_Oracle_Inequality} 
	Assume that we observe $\left(x_{[n]},y_{[n]}\right)\in \left([0,1]^p\times\R^q\right)^n$, coming from an unknown conditional PDF $s_0 \equiv s_{\psi_0} \in S$, defined in \eqref{eq_Sm_bounded}. 
	Given $C_{1n}$ in \eqref{eq_l1_oracle_error_upperbound}, if $\lambda$ satisfies \eqref{eq_lowerbound_kappa},
	the Lasso estimator $\widehat{s}^{\lasso}_{\lambda}$, defined in \eqref{eq_lasso_estimator},  satisfies the $l_1$-oracle inequality:
	\begin{align} \label{eq_l1_Oracle_Inequality}
		\E{}{{\kl}_n\left(s_0,\widehat{s}^{\lasso}_{\lambda}\right)}&\le{} \frac{\kappa+1}{\kappa} \inf_{s_\psi \in S} \Big[{\kl}_n\left(s_0,s_\psi\right) +\lambda (\norm{\gamma}_1+\norm{\vect(\beta)}_1)\Big]  +\lambda
		+\sqrt{\frac{K}{n}} C_{1n}.
	\end{align}  
	%
\end{theorem}

%
\subsection{Discussion and perspectives regarding our oracle inequality}\label{sec_discussion_comparison}
%
%
{\bf The oracle model and convergence rate for the Lasso estimator.}
Theorem~\ref{thm_l1_Oracle_Inequality} characterizes the performance of Lasso estimators as $l_1$-PMLEs for SGMoE models. If the regularization parameter $\lambda$ is properly chosen, the solution to the \emph{$l_1$-penalized empirical risk minimization problem}, behaves in a manner comparable to the deterministic Lasso (the so-called \emph{oracle}). This oracle is the solution to the \emph{$l_1$-penalized true risk minimization problem}, up to an error term of order $\lambda$. Note that \emph{the best model}, denoted by $s_{\psi^*}$, is defined as the one with the smallest $l_1$-penalized risk:
\begin{align}\label{eq_define_oracle}
	\inf_{s_\psi \in S} \left[{\kl}_n\left(s_0,s_\psi\right)+\lambda (\norm{\gamma}_1+\norm{\vect(\beta)}_1)\right].
\end{align}
However, since we do not know true density $s_0$, we cannot select this best model, which we call the \emph{oracle model}. In particular, by definition, the oracle is the model in the collection that minimizes the
$l_1$-penalized risk in \eqref{eq_define_oracle}, which is generally assumed to be unknown.
From the oracle inequality of  Theorem~\ref{thm_l1_Oracle_Inequality}, we conjecture that by constructing a suitable approximation theory on a good space, we can control this $l_1$-penalized true risk to obtain the parametric convergence rate of $n^{-1/2}$ for the Lasso estimator.
A related work in this direction is \citet{massart_rates_2012}, who established convergence rates for the selected Lasso estimator of \citet{massart_lasso_2011}, for a wide range of function classes described by the interpolation spaces of \citet{barron_approximation_2008}.
Furthermore, to the best of our knowledge, Theorem 2.8 of \citet{maugis2013adaptive} is the only result in the literature that investigates the minimax estimator for Gaussian mixture models but for model selection instead of Lasso. 
%
%
We will leave the non-trivial task of obtaining Lasso extensions of this result to future work.

\noindent {\bf Implementable Lasso estimator and data-driven regularization parameter $\lambda$.} 
Note that  Theorem~\ref{thm_l1_Oracle_Inequality} ensures that there exists a sufficiently large $\lambda$ for which the estimate has good properties, but does not give an explicit value for $\lambda$. However, we at least give the lower bound on the value of $\lambda$ via the bound $\lambda \ge \kappa C(p,q,n,K)$, where $\kappa \ge 148$, although this value is obviously conservative.
Moreover, it is important to note that the Lasso estimator that appears in  Theorem~\ref{thm_l1_Oracle_Inequality} has already been implemented in practice.
Indeed, when the $l_2$-penalties are given zero weight in \citet{khalili2010new} and \citet{chamroukhi2019regularizedIJCNN,chamroukhi2019regularized}, their penalty functions and a recent result from \citet{huynh2019estimation} for generalized linear expert models belong to our framework and the $l_1$-oracle inequality from  Theorem~\ref{thm_l1_Oracle_Inequality} provides further theoretical insight for these Lasso estimators.
In particular, possible solutions for calibrating the tuning parameter $\lambda$ of the penalty from the data are the BIC~\citep{schwarz_estimating_1978} used in \citet{chamroukhi2019regularizedIJCNN,chamroukhi2019regularized} and the generalized cross-validation~\citep{stone_cross_validatory_1974} utilized in \citet{khalili2010new,khalili2007variable}.
Given the available literature, such computational strategies and numerical simulations for real data sets will not be considered and discussed further here.

\noindent {\bf Dependency on $p$, $q$, $n$, and $K$ in the lower bound of $\lambda$.} Note that we recover the same dependence of the form $\sqrt{\ln \left(2p +1\right)}$ as for the homogeneous linear regression in \citet{stadler2010l1} and of the form $\sqrt{\ln \left(2p +1\right)} \left(\ln n\right)^2/ \sqrt{n}$ for the mixture of regressions models in \citet{Meynet:2013aa}. On the contrary, the dependence on $q$ for the mixture of multivariate Gaussian regression models in \citet{devijver2015l1} has form $q^2 + q$, while here we get the form $q^2 \sqrt{q}$. The main reason is that  the class $S$ of the SGMoE model is larger, and we use a different technique to evaluate the upper bound of the uniform norm of the gradient for each element in $S$. 
%
%
In the lower bound of $\lambda$ in \eqref{eq_lowerbound_kappa}, we can potentially get the factor dependence $K^2$ instead of $K$ as in \citet{Meynet:2013aa} and \citet{devijver2015l1}. This can be explained by the fact that we used a different technique to handle the more complex model when dealing with the upper bound on the uniform norm of the gradient of $\ln s_\psi$, for $s_\psi \in S$.
We refer to \citet[Remark 5.8]{Meynet:2013aa} for some data sets where the dependence on $K$ can be reduced to the order of $\sqrt{K}$ for the mixture Gaussian regression models. The determination of optimal rates for such problems is still open. 
Furthermore, the dependence on $n$ for the homogeneous linear regression in \citet{stadler2010l1} is of the order of $n^{-1/2}$, while here we have an additional $\left(\ln n\right)^2$ factor. In fact, the same situation can be found in the $l_1$-oracle inequalities of \citet{Meynet:2013aa} and \citet{devijver2015l1}. 
As explained in \citet{Meynet:2013aa}, the use of nonlinear KL information leads to a scenario where the linearity arguments developed in \citet{stadler2010l1} with the quadratic loss function cannot be exploited. Instead, we need to use the entropy arguments to for our model, which leads to an additional $\left(\ln n\right)^2$ factor.

\noindent {\bf Multiplicative upper bound constant.} It is worth noting that the constant $1+\kappa^{-1}$ appearing in the upper bounds of  Theorem~\ref{thm_l1_Oracle_Inequality}  cannot be reduced to $1$, which is in fact the same situation as the constant from $C_1$ from~\citet[Theorem 1]{montuelle_mixture_2014}.
Note that this problem also occurred in the $l_1$-oracle inequalities of \citet{Meynet:2013aa}, and \citet{devijver2015l1}.  
Deriving an oracle inequality such that $1+ \kappa^{-1}$ can be replaced by $1$ for the KL loss is still an open problem.

\noindent {\bf Model misspecification.} 
In  Theorem~\ref{thm_l1_Oracle_Inequality}, when $s_0 \notin S$, by letting $n \rightarrow \infty$, due to the large bias from the first upper bound term, the total error upper bound of \eqref{eq_l1_Oracle_Inequality} converges to 
\begin{align} 
	&\left(1 + \frac{1}{\kappa}\right) \inf_{s_\psi \in S} \left[\lim_{n \rightarrow \infty}{\kl}_n\left(s_0,s_\psi\right)+\lambda (\norm{\gamma}_1+\norm{\vect(\beta)}_1)\right] + \lambda,\nn
\end{align}  
which may be large.
The same conclusion holds for  Theorem~\ref{thm_l1_Oracle_Inequality_Ball} when $s_0 \notin \cup_{m\in \Ns} S_m$.
Nevertheless, as we consider SGMoE models, some recent universal approximation results, see, \eg~\citet{nguyen_universal_2016,nguyen2019approximation,nguyen_approximation_2020,nguyen2020approximationMoE,nguyen_approximation_2022}, imply that if we take a sufficiently large number of mixture components $K$, that is sufficiently large class $S$, we can approximate a broad class of conditional PDFs, and thus the term on the right hand side is small for $K$ sufficiently large. This improves the error bound even when $s_0 \notin S$. 
%


%
\subsection{Comparison with the state-of-the-art}\label{sec_comparisons}

{\bf Standard asymptotic results with variable selection.}  Theorem~\ref{thm_l1_Oracle_Inequality} complements the standard asymptotic results for high-dimensional SGMoE models via feature selection using Lasso, as well asthe more general PMLE via the Scad penalty function of \citet{khalili2010new}.  
%
%
%
By extending the theoretical developments for mixture of regressions models in \citet{khalili2007variable}, standard asymptotic theorems for SGMoE are established in~\citet{khalili2010new}. Then, under several strict regularity conditions on the true joint density function and the choice of tuning parameter, the PMLE from \citet{khalili2010new}, using the Scad penalty function from \citet{fan2001variable} instead of Lasso, is proved to be both consistent in feature selection and maintains $\sqrt{n}$ consistency. 
On the contrary, the only assumption used to obtain the $n^{-1/2}$ convergence rate of the error upper bound in \eqref{eq_l1_Oracle_Inequality} from  Theorem~\ref{thm_l1_Oracle_Inequality} is boundedness on the parameter space. In fact, this boundedness condition is also required by \citet{khalili2010new}.
Furthermore, we work directly on conditional PDFs with fixed covariates rather than focusing on joint PDFs as in \citet{khalili2010new,khalili2007variable}. In future work, we shall investigate whether our proof technique used in this paper can be adapted to the problem of estimating joint PDFs when the predictors are random variables.

\noindent {\bf Boundedness assumptions on the parameter space.} 
It is worth noting that our boundedness assumptions also appeared in \citet{stadler2010l1,Meynet:2013aa,devijver2015l1}. 
They are quite natural when working with MLE \citep{maugis2011non}, at least when considering the problem of the unboundedness of the likelihood at the boundaries of the parameter space \citep{redner_mixture_1984,mclachlan2000finite}, and to prevent the likelihood from diverging. 

\subsection{Proof of Theorem \ref{thm_l1_Oracle_Inequality}}
%
%
{\bf Proof sketch of Theorem \ref{thm_l1_Oracle_Inequality}.} At a high level, the main idea here is to study the Lasso estimator,  Theorem~\ref{thm_l1_Oracle_Inequality}, as a solution of the $l_1$-ball PMLE,  Theorem~\ref{thm_l1_Oracle_Inequality_Ball}, which is defined later in this section. The proof of  Theorem~\ref{thm_l1_Oracle_Inequality_Ball} can be deduced from Propositions \ref{thm_l1_Oracle_Inequality_Ball_small} and \ref{thm_l1_Oracle_Inequality_Ball_large}, which deal with the cases for small and large values of $Y$ and are proved in the supplementary materials. 

\noindent {\bf $l_1$-ball PMLE.} We need to define the $l_1$-ball PMLE for the statement of Theorem \ref{thm_l1_Oracle_Inequality_Ball}. For this, by restricting $S$, to a suitable $l_1$-ball of the parameters $\gamma,\beta$ on the definition of $S$, we define a collection of $l_1$-ball models $S_m$, where $m\in \Ns$ is a radius of the $l_1$-ball, as follows:
\begin{align}\label{eq_model_Sm}
	S_m = \left\{s_\psi \in S \mid \norm{\gamma}_1+\norm{\vect(\beta)}_1 \le m\right\}.
\end{align}
Then for some $\eta_m \geq 0$, let $\widehat{s}_m$ be a $\eta_m$-log-likelihood estimator (LLE) in $S_m$, defined as:
\begin{align}\label{eq_model_Sm_LL}
	-\frac{1}{n} \sum_{i=1}^n \ln\left(\widehat{s}_m\left(y_i| x_i\right)\right) &\le \inf_{s_m \in S_m} \left(-\frac{1}{n} \sum_{i=1}^n\ln\left(s_m\left(y_i| x_i\right)\right) \right)+ \eta_m.
\end{align}
As is always the case, it is not sufficient to use the LLE of the estimate in each model as a criterion. It underestimates the risk of the estimate and the result is a choice of model that is too complex. In the context of the PMLE, by adding an appropriate penalty $\pen(m)$, one hopes to create a trade-off between good data fit and model complexity.
Suppose that for all $m \in \Ns$, the penalty function satisfies $\pen(m) = \lambda m$, where $\lambda$ will be determined as in \eqref{eq_lowerbound_kappa}. Then, for some $\eta \geq 0$, an $l_1$-ball PMLE is defined as $\widehat{s}_{\widehat{m}}$, where $\widehat{m}$ satisfies
\begin{align}\label{eq_model_Sm_PMLE}
	&-\frac{1}{n} \sum_{i=1}^n \ln\left(\widehat{s}_{\widehat{m}}\left(y_i| x_i\right)\right) + \pen(\widehat{m}) \le \inf_{m \in \Ns} \left(-\frac{1}{n} \sum_{i=1}^n\ln\left(\widehat{s}_m\left(y_i| x_i\right)\right) + \pen(m) \right) + \eta.
\end{align}
Note that the error terms $\eta_m$ and $\eta$ are necessary to avoid any existence problems, \eg~the infimum may not be reached. Roughly speaking, the Ekeland variational principle states that for any extended-valued lower semicontinuous function, which is bounded below, one can add a small perturbation to ensure the existence of the minimum, see \eg~\citet{borwein_techniques_2004}. This framework is also used in \citet{montuelle_mixture_2014}, and \citet{nguyen_non_asymptotic_2022}.
Next, we state an $l_1$-ball model selection via Theorem \ref{thm_l1_Oracle_Inequality_Ball}.
%
\begin{theorem}[$l_1$-ball model selection]\label{thm_l1_Oracle_Inequality_Ball}
	Assume that $\left(x_{[n]},y_{[n]}\right)\in \left([0,1]^p\times\R^q\right)^n$, come from an unknown conditional PDF $s_0 \equiv s_{\psi_0} \in S$, defined in \eqref{eq_Sm_bounded}. Given $C_{1n}$ in \eqref{eq_l1_oracle_error_upperbound},
	if $\lambda$ satisfies \eqref{eq_lowerbound_kappa}, the $l_1$-ball PMLE $\widehat{s}_{\widehat{m}}$, defined in \eqref{eq_model_Sm_PMLE},  satisfies the oracle inequality:
	\begin{align}
		&\E{}{{\kl}_n\left(s_0,\widehat{s}_{\widehat{m}}\right)}\le \left(1 + \frac{1}{\kappa}\right)\inf_{m\in\Ns}\left(\inf_{s_m \in S_m} {\kl}_n\left(s_0,s_m\right)+ \pen(m) + \eta_m\right)  + \eta  +\sqrt{\frac{K}{n}} C_{1n}.
		\label{eq_l1_Oracle_Inequality_Ball}
	\end{align}
\end{theorem}

%

%

\noindent
{\bf Proof of Theorem \ref{thm_l1_Oracle_Inequality}.} Let $\lambda > 0$ and define $\widehat{m}$ to be the smallest integer such that $\widehat{s}^{\lasso}_{\lambda}$ belongs to $S_{\widehat{m}}$, \ie~$\widehat{m} := \ceil*{\left\|\psi^{[1,2]}\right\|_{1}} \le \left\|\psi^{[1,2]}\right\|_{1} +1 $. Then using the definition of $\widehat{m}$, \eqref{eq_lasso_estimator}, \eqref{eq_model_Sm}, and $S =\bigcup_{m \in \Ns} S_m$, we get
\begin{align*}
	-\frac{1}{n}\sum_{i=1}^n \ln\left(\widehat{s}^{\lasso}_{\lambda}\left(y_i | x_i\right)\right) + \lambda \widehat{m}
	&\le -\frac{1}{n}\sum_{i=1}^n \ln\left(\widehat{s}^{\lasso}_{\lambda}\left(y_i | x_i\right)\right) + \lambda \left(\left\|\psi^{[1,2]}\right\|_{1}+1\right)\nn \\
	&=\inf_{s_\psi\in S} \left(-\frac{1}{n}\sum_{i=1}^n \ln\left(s_\psi\left(y_i | x_i\right)\right) + \lambda \left\|\psi^{[1,2]}\right\|_{1}\right) + \lambda\nn\\
	&=\inf_{m \in \Ns}\left(\inf_{s_\psi\in S_m} \left(-\frac{1}{n}\sum_{i=1}^n \ln\left(s_\psi\left(y_i | x_i\right)\right) + \lambda \left\|\psi^{[1,2]}\right\|_{1}\right)\right) + \lambda\nn\\
	&\le\inf_{m \in \Ns}\left(\inf_{s_m\in S_m} \left(-\frac{1}{n}\sum_{i=1}^n \ln\left(s_m\left(y_i | x_i\right)\right) + \lambda m\right)\right) + \lambda,
\end{align*}
which implies
\begin{align*}
	-\frac{1}{n}\sum_{i=1}^n \ln\left(\widehat{s}^{\lasso}_{\lambda}\left(y_i | x_i\right)\right) + \pen(\widehat{m}) \le  \inf_{m \in \Ns} \left(-\frac{1}{n}\sum_{i=1}^n \ln\left(\widehat{s}_m\left(y_i | x_i\right)\right) + \pen(m)\right) + \eta
\end{align*}
with $\pen(m)=\lambda m, \eta = \lambda$, and $\widehat{s}_m$ is a $\eta_m$-log-likelihood minimizer in $S_m$, with $\eta_m \ge 0$ defined by \eqref{eq_model_Sm_LL}.
Thus, $\widehat{s}^{\lasso}_{\lambda}$ satisfies \eqref{eq_model_Sm_PMLE} with $\widehat{s}^{\lasso}_{\lambda}\equiv \widehat{s}_{\widehat{m}}$, \ie
\begin{align}
	-\frac{1}{n} \sum_{i=1}^n \ln\left(\widehat{s}_{\widehat{m}}\left(y_i| x_i\right)\right) + \pen(\widehat{m})\le \inf_{m \in \Ns} \left(-\frac{1}{n} \sum_{i=1}^n\ln\left(\widehat{s}_m\left(y_i| x_i\right)\right) + \pen(m) \right) + \eta.\nn
\end{align}
Then,  Theorem~\ref{thm_l1_Oracle_Inequality_Ball} implies that
if
\begin{align*} 
	\lambda \geq \kappa\frac{KB_{0n}}{\sqrt{n}} \left(q\ln n \sqrt{\ln(2p+1)}+1\right), \\
	B_{0n} = \max\left(A_\Sigma,1+KA_G\right)\left[1+2q\sqrt{q} A_\Sigma \left( 5A^2_\beta + 4 A_\Sigma \ln n\right) \right],
\end{align*} 
for some absolute constants $\kappa \geq 148$, Theorem \ref{thm_l1_Oracle_Inequality} holds as required.

\noindent
	{\bf Proof of Theorem \ref{thm_l1_Oracle_Inequality_Ball}.} Given any $M_n >0$, this can be done by introducing an event $\cT$ and a space $F_m$ as follows:
	\begin{align}
		&\cT= \left\{\max_{i \in [n]}  \left\|Y_i\right\|_\infty=\max_{i \in [n]} \max_{z \in [q]} \left|\left[Y_i\right]_z\right| \le M_n \right\},
		%
		\cT^C= \left\{\max_{i \in [n]}  \left\|Y_i\right\|_\infty=\max_{i \in [n]} \max_{z \in [q]} \left|\left[Y_i\right]_z\right| > M_n \right\},\nn\\
		%
		&F_m = \left\{f_m = -\ln\left(\frac{s_m}{s_0}\right)= \ln(s_0)-\ln(s_m), s_m \in S_m \right\}. \nn
	\end{align}
	Conditional on $\left\{x_i\right\}_{i\in[n]}$, let $Y'_{[n]}| x_{[n]} \equiv \left(Y'_i| x_i\right)_{i \in [n]}$ be IID random samples from $Y$ arising from the conditional PDF $s_0\left(\cdot| x_i\right), i \in [n]$. They are independence copies of the sample $Y_{[n]}| x_{[n]}$.
	By taking into account the definition of average KL loss function from \eqref{Extend.CarolineMeynet.eq.2.2}, the conditional expectation property, we obtain	
	\begin{align}
		{\kl}_n(s_0,\widehat{s}_{\widehat{m}}) 
		&=  \E{Y'_{[n]}| x_{[n]}}{ 	\frac{1}{n} \sum_{i=1}^n \widehat{f}_{\widehat{m}}\left(Y_i| x_i\right) | \cT } {\Indi}_{\cT} + \E{Y'_{[n]}| x_{[n]}}{ 	\frac{1}{n} \sum_{i=1}^n \widehat{f}_{\widehat{m}}\left(Y_i| x_i\right) | \cT^C } {\Indi}_{\cT^C}\nn\\
		&\equiv	\left({\kl}_n(s_0,\widehat{s}_{\widehat{m}}) | \cT\right){\Indi}_{\cT} + 	\left({\kl}_n(s_0,\widehat{s}_{\widehat{m}})  |\cT^C\right){\Indi}_{\cT^c}.
		\label{eq_AKL_conditional_explanation}
	\end{align}
	%

	From now on, when there is no confusion,
	the expectation of \eqref{eq_AKL_conditional_explanation} is written as follows:
	\begin{align}
		\E{}{{\kl}_n\left(s_0,\widehat{s}_{\widehat{m}}\right)}
		%
		& = \E{Y_{[n]}}{\left({\kl}_n(s_0,\widehat{s}_{\widehat{m}}) | \cT\right){\Indi}_{\cT}} +\E{Y_{[n]}}{\left({\kl}_n(s_0,\widehat{s}_{\widehat{m}})  |\cT^C\right){\Indi}_{\cT^c}}\nn\\
		&  \equiv \E{}{{\kl}_n\left(s_0,\widehat{s}_{\widehat{m}}\right) {\Indi}_{\cT}} + \E{}{{\kl}_n\left(s_0,\widehat{s}_{\widehat{m}}\right) {\Indi}_{\cT^C}}.   \label{eq_AKL_conditional_explanation_simpler}
	\end{align}
	Therefore, on the basis of the above remark \eqref{eq_AKL_conditional_explanation_simpler}, Theorem~\ref{thm_l1_Oracle_Inequality_Ball} is proved by obtaining the upper bound for each of the following terms using Propositions \ref{thm_l1_Oracle_Inequality_Ball_small} and \ref{thm_l1_Oracle_Inequality_Ball_large}:
	\begin{align}
		&\E{}{{\kl}_n\left(s_0,\widehat{s}_{\widehat{m}}\right)}  =  \E{}{{\kl}_n\left(s_0,\widehat{s}_{\widehat{m}}\right) {\Indi}_{\cT}} + \E{}{{\kl}_n\left(s_0,\widehat{s}_{\widehat{m}}\right) {\Indi}_{\cT^C}}. \nn
	\end{align}

	Given some constants $\kappa \geq 148$, we need to define the following condition for $\lambda$, given some constant $M_n > 0$:
	\begin{align}
		%
		%
		&\lambda \ge \kappa \frac{K}{\sqrt{n}} C_{3n}, \text{ }C_{3n}  =  B_{n}\left(q\ln n \sqrt{\ln(2p+1)}+1\right), 	\label{eq_lowerbound_kappa_smallY}\\
		&B_{n} = \max\left(A_\Sigma,1+KA_G\right)\left[1+ q\sqrt{q}\left(M_n+ A_\beta\right)^2 A_\Sigma\right]
		\label{eq_Bn}.
	\end{align}

	\begin{proposition}[Small values of $Y$]\label{thm_l1_Oracle_Inequality_Ball_small}
		Assume that  $\left(x_{[n]},y_{[n]}\right)\in \left([0,1]^p\times\R^q\right)^n$ comes from an unknown conditional PDF $s_0 \equiv s_{\psi_0} \in S$ defined in \eqref{eq_Sm_bounded}. Given $C_{2n}$ and $B_{n}$ defined as in \eqref{eq_l1_oracle_smallY_error_upperbound} and \eqref{eq_Bn}, respectively,
		if $\lambda$ satisfies \eqref{eq_lowerbound_kappa_smallY}, the $l_1$-ball PMLE $\widehat{s}_{\widehat{m}}$, defined in \eqref{eq_model_Sm_PMLE},  satisfies: 
		%
		\begin{align}
			&\E{}{{\kl}_n\left(s_0,\widehat{s}_{\widehat{m}}\right){\Indi}_{\cT}} \le{}\left(1+\kappa^{-1}\right)\inf_{m\in\Ns}\left(\inf_{s_m \in S_m} {\kl}_n\left(s_0,s_m\right)+ \pen(m) + \eta_m\nn\right)  +\eta + \sqrt{\frac{K}{n}} B_{n} C_{2n}.
		\end{align}
	\end{proposition}
	\begin{proposition}[Large values of $Y$]\label{thm_l1_Oracle_Inequality_Ball_large}
		Consider $s_0,\cT$, and $\widehat{s}_{\widehat{m}}$ as defined in Proposition \ref{thm_l1_Oracle_Inequality_Ball_small} and $H_{s_0}$ as defined in \eqref{eq_l1_oracle_smallY_error_upperbound}. 
		Then,
		\begin{align*}		
			&\E{}{{\kl}_n\left(s_0,\widehat{s}_{\widehat{m}}\right){\Indi}_{\cT^C}} \le  \left(\frac{ e^{q/2-1}\pi^{q/2} }{A_\Sigma ^{q/2} } + H_{s_0}\right)  \sqrt{2 K n q A_\gamma} e^{-\frac{M_n^2-2M_n A_\beta}{4 A_\Sigma}}.
		\end{align*}	
	\end{proposition}
	Proposition \ref{thm_l1_Oracle_Inequality_Ball_small} constitutes our important technical contribution.
	Via Lemma \ref{lem_control_deviation}, the main idea to prove Proposition \ref{thm_l1_Oracle_Inequality_Ball_small} is to control the following deviation on the event $\cT$:
	\begin{align}
		&\sup_{f_{m} \in F_{m}} \left|\nu_n\left(-f_{m}\right)\right|\equiv \sup_{f_{m} \in F_{m}} \left|\frac{1}{n} \sum_{i=1}^n \left\{f_m\left(Y_i| x_i\right)-\E{}{f_m\left(Y_i| x_i\right)}\right\}\right|, 	\label{eq_define_deviation_sketch}\\
		&F_m = \left\{f_m = -\ln\left(\frac{s_m}{s_0}\right)= \ln(s_0)-\ln(s_m), s_m \in S_m \right\} \label{eq_Fm}.
	\end{align}
	%

	%
	\begin{lemma}[Control of deviation]\label{lem_control_deviation}
		%
		For each $m' \in \Ns$, let
		\begin{align}
			\Delta_{m'} &=  m' \sqrt{\ln(2p+1)}\ln n+2\sqrt{K} \left(A_\gamma+  q A_\beta + \frac{ q\sqrt{q}}{a_{\Sigma}}\right)
			\label{eq_Delta_m}.
		\end{align}		
		Then, 
		on the event $\cT$, 
		for all $m'\in \Ns$, and for all $t>0$, with probability greater than $1-e^{-t}$,
		\begin{align}\label{eq_model_Sm_LL1.Gaussian}
			&\sup_{f_{m'} \in F_{m'}} \left|\nu_n\left(-f_{m'}\right)\right|\le\frac{4KB_n}{\sqrt{n}} \left[37 q \Delta_{m'} + \sqrt{2}\left(A_\gamma+  q A_\beta + \frac{ q\sqrt{q}}{a_{\Sigma}}\right)\sqrt{t}\right].
		\end{align}
	\end{lemma}
	The proof of Lemma \ref{lem_control_deviation} appears in the supplementary material
	and follows the arguments developed in the proof of \citet[Theorem 7.11]{massart_concentration_2007}. The proof of Proposition \ref{thm_l1_Oracle_Inequality_Ball_small} is in the spirit of Vapnik's method of structural risk minimization, first established in \citet{vapnik1982estimation} and briefly summarized in Section 8.2 of \citet{massart_concentration_2007}. In particular, we use concentration inequalities combined with symmetrization arguments to obtain an upper bound on the empirical process in expectation from \eqref{eq_define_deviation_sketch}.  Our technique combines Vapnik's structural risk minimization paradigm (\eg\citealp{vapnik1982estimation}) and model selection theory for conditional density estimation (\eg\citealp{cohen2011conditional}), which extend the density estimation results of~\citet{massart_concentration_2007}.


\section{Numerical experiments}\label{sec_numerical_experiment}
In this section, we empirically validate the convergence rate of the error upper bound in \eqref{eq_l1_Oracle_Inequality} from  Theorem~\ref{thm_l1_Oracle_Inequality} in our SGMoE models. For simplicity, we only perform a simulation study to illustrate
the convergence rates when $\cX \subset \R^p, \quad p = 6$, and $\cY \subset \R^q, \quad q = 1$.
All the following simulations were performed in
R 4.3.2 on a standard Unix machine.
We construct simulated data sets sampling from the true conditional density, $s_0$, belongs to the class of SGMoE models $S$:
\begin{align} 
	s_{0}(y| x) = \sum_{k=1}^2 \frac{\exp\left(\gamma_{0k0} + \gamma_{0k}^\top x\right)}{\sum_{l=1}^2\exp\left(\gamma_{0l0} + \gamma_{0l}^\top x\right)} \cN\left(y;\beta_{0k0}+ \beta_{0k} x,\Sigma_{0k}\right).\nn
\end{align}
Here, the true parameters for the true SGMoE model are given by:
\begin{align*}
	(\gamma_{010}, \gamma_{01})^\top &= (1,2,0,0,-1,0,0)^\top; \quad \Sigma_{01}=\Sigma_{02} = 1
	;\\
	(\beta_{010}, \beta_{01})^\top &= (0,0,1.5,0,0,0,1)^\top; \quad (\beta_{020}, \beta_{02})^\top = (0,1,-1.5,0,0,2,0)^\top.
\end{align*}
Here we implement the $l_1$-PMLE using the EM algorithm for the SGMoE model with coordinate ascent algorithm for updating the gating network, following the strategy of \cite{chamroukhi2019regularizedIJCNN,chamroukhi2019regularized}.

\begin{figure}[ht!]
	\centering
	\includegraphics[width=.75\linewidth]{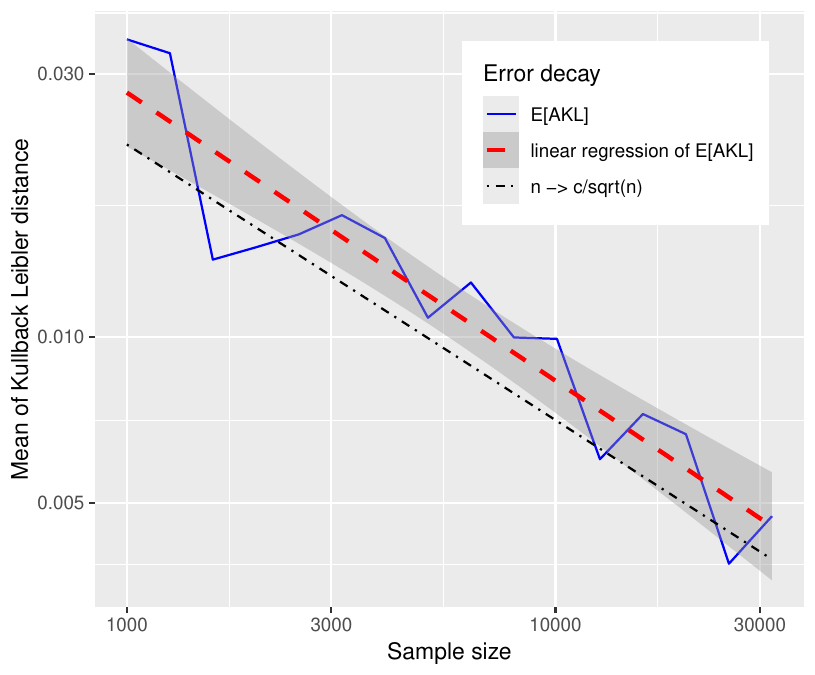}
	\caption{Average KL divergence between the true and selected densities based on the Lasso estimator, represented in a log-log scale, using $100$ different choices of sample size $n$ between $1000$ and $32000$ over $n_t = 20$ trials and $n_y = 30$. A free least-square regression with confidence interval and a regression with slope $-1/2$ were added to stress the two different behavior for each graph.}
	\label{fig_error_decay}
\end{figure}

We want to empirically validate the convergence rate of the error upper bound in terms of the KL divergence, which cannot be exactly calculated in the case of Gaussian mixtures.
Thus, we assess the divergence through a Monte Carlo simulation, given that we have knowledge of the true density. It is important to mention that the variability in this randomized approximation has been shown to be minimal in practice, a fact that is corroborated by the numerical experiments conducted by \cite{nguyen_non_asymptotic_2022,montuelle_mixture_2014}. Specifically, we calculate the Monte Carlo approximation for the average KL divergence ${\kl}_n(s_0,\widehat{s}^{\lasso}_{\lambda})$ as described below:
\begin{align*}
	\frac{1}{n} \sum_{i=1}^n \kl\left(s_0\left(\cdot| x_i\right),\widehat{s}^{\lasso}_{\lambda}\left(\cdot| x_i\right)\right)
	\approx \frac{1}{nn_y}\sum_{i=1}^n \sum_{j=1}^{n_y} \ln\left(\frac{s_0\left(y_{ij}| x_i\right)}{\widehat{s}^{\lasso}_{\lambda}\left(y_{ij}| x_i\right)}\right).
\end{align*}
Here $(y_{ij})_{j\in[n_y]}$ are drawn from $s_0\left(\cdot| x_i\right)$. 
Then $\E{}{{\kl}_n\left(s_0,\widehat{s}^{\lasso}_{\lambda}\right)}$ is approximated again by averaging over $n_t$ Monte Carlo trials. Therefore, the simulated data used for approximation can be written as $(x_{i},y_{ij})_{t}$ with $i \in [n], j\in[n_y], t \in[n_t]$. 
Figure \ref{fig_error_decay} illustrates that the error decreases with order \( C_{1n}\sqrt{K}/\sqrt{n} \), as predicted theoretically in Theorem~\ref{thm_l1_Oracle_Inequality}, as the sample size \( n \) increases when applying the penalty based on our criterion.

\section{Conclusions}
\label{sec_conclusions_perspectives}
To the best of our knowledge, we are the first to establish an $l_1$-oracle inequality for  SGMoE models from a non-asymptotic perspective, under the mildest assumptions, namely the boundedness of the parameter space, which provides a lower bound for the regularization parameter of the Lasso, while ensuring non-asymptotic theoretical control of the KL loss of the estimator.
We further conduct a simulation study to empirically confirm our theoretical results.
We believe that our contribution assists in  further popularizing MoE models by providing a theoretical basis for their application to heterogeneous high-dimensional data.

\bibliography{acml2024}

\begin{thebibliography}{56}
\providecommand{\natexlab}[1]{#1}
\providecommand{\url}[1]{\texttt{#1}}
\expandafter\ifx\csname urlstyle\endcsname\relax
  \providecommand{\doi}[1]{doi: #1}\else
  \providecommand{\doi}{doi: \begingroup \urlstyle{rm}\Url}\fi

\bibitem[Arlot and Celisse(2010)]{arlot_survey_2010}
Sylvain Arlot and Alain Celisse.
\newblock A survey of cross-validation procedures for model selection.
\newblock \emph{Statistics Surveys}, 4:\penalty0 40--79, January 2010.

\bibitem[Barron et~al.(1999)Barron, Birg{\'e}, and Massart]{barron_risk_1999}
Andrew Barron, Lucien Birg{\'e}, and Pascal Massart.
\newblock Risk bounds for model selection via penalization.
\newblock \emph{Probability theory and related fields}, 113:\penalty0 301--413,
  1999.

\bibitem[Barron et~al.(2008)Barron, Cohen, Dahmen, and
  DeVore]{barron_approximation_2008}
Andrew~R. Barron, Albert Cohen, Wolfgang Dahmen, and Ronald~A. DeVore.
\newblock Approximation and learning by greedy algorithms.
\newblock \emph{The Annals of Statistics}, 36\penalty0 (1):\penalty0 64 -- 94,
  2008.

\bibitem[Borwein and Zhu(2004)]{borwein_techniques_2004}
Jonathan~M Borwein and Qiji~J Zhu.
\newblock \emph{Techniques of {Variational} {Analysis}}.
\newblock Springer New York, 2004.

\bibitem[Chamroukhi and Huynh(2018)]{chamroukhi2019regularizedIJCNN}
Faicel Chamroukhi and Bao~Tuyen Huynh.
\newblock Regularized maximum-likelihood estimation of mixture-of-experts for
  regression and clustering.
\newblock In \emph{2018 International Joint Conference on Neural Networks
  (IJCNN)}, pages 1--8, 2018.

\bibitem[Chamroukhi and Huynh(2019)]{chamroukhi2019regularized}
Faicel Chamroukhi and Bao~Tuyen Huynh.
\newblock Regularized maximum likelihood estimation and feature selection in
  mixtures-of-experts models.
\newblock \emph{Journal de la Soci{\'e}t{\'e} Fran{\c{c}}aise de Statistique},
  160\penalty0 (1):\penalty0 57--85, 2019.

\bibitem[Cohen and {Le Pennec}(2011)]{cohen2011conditional}
S~X Cohen and Erwan {Le Pennec}.
\newblock Conditional density estimation by penalized likelihood model
  selection and applications.
\newblock \emph{Technical report, INRIA}, 2011.

\bibitem[Cover(1999)]{cover1999elements}
Thomas~M Cover.
\newblock \emph{Elements of information theory}.
\newblock John Wiley \& Sons, 1999.

\bibitem[Devijver(2015)]{devijver2015l1}
Emilie Devijver.
\newblock {An $l_1$-oracle inequality for the Lasso in multivariate finite
  mixture of multivariate Gaussian regression models}.
\newblock \emph{ESAIM: PS}, 19:\penalty0 649--670, 2015.

\bibitem[Duistermaat and Kolk(2004)]{duistermaat2004multidimensional}
Johannes~Jisse Duistermaat and Johan~AC Kolk.
\newblock \emph{Multidimensional real analysis I: differentiation}, volume~86.
\newblock Cambridge University Press, 2004.

\bibitem[Fan and Li(2001)]{fan2001variable}
Jianqing Fan and Runze Li.
\newblock Variable selection via nonconcave penalized likelihood and its oracle
  properties.
\newblock \emph{Journal of the American statistical Association}, 96\penalty0
  (456):\penalty0 1348--1360, 2001.

\bibitem[Genovese and Wasserman(2000)]{genovese2000rates}
Christopher~R Genovese and Larry Wasserman.
\newblock {Rates of convergence for the Gaussian mixture sieve}.
\newblock \emph{The Annals of Statistics}, 28\penalty0 (4):\penalty0
  1105--1127, aug 2000.

\bibitem[Golub and Van~Loan(2012)]{golub2012matrix}
Gene~H Golub and Charles~F Van~Loan.
\newblock \emph{Matrix computations}, volume~3.
\newblock JHU press, 2012.

\bibitem[Ho and Nguyen(2016)]{ho_convergence_2016}
Nhat Ho and XuanLong Nguyen.
\newblock Convergence rates of parameter estimation for some weakly
  identifiable finite mixtures.
\newblock \emph{The Annals of Statistics}, 44\penalty0 (6):\penalty0 2726 --
  2755, 2016.

\bibitem[Ho et~al.(2022)Ho, Yang, and Jordan]{ho_convergence_2022}
Nhat Ho, Chiao-Yu Yang, and Michael~I. Jordan.
\newblock Convergence {Rates} for {Gaussian} {Mixtures} of {Experts}.
\newblock \emph{Journal of Machine Learning Research}, 23\penalty0
  (323):\penalty0 1--81, 2022.

\bibitem[Horn and Johnson(2012)]{horn2012matrix}
Roger~A Horn and Charles~R Johnson.
\newblock \emph{{Matrix analysis}}.
\newblock Cambridge University Press, 2012.

\bibitem[Huynh and Chamroukhi(2019)]{huynh2019estimation}
Bao~Tuyen Huynh and Faicel Chamroukhi.
\newblock {Estimation and feature selection in mixtures of generalized linear
  experts models}.
\newblock \emph{arXiv preprint arXiv:1907.06994}, 2019.

\bibitem[Jacobs et~al.(1991)Jacobs, Jordan, Nowlan, and
  Hinton]{jacobs1991adaptive}
Robert~A Jacobs, Michael~I Jordan, Steven~J Nowlan, and Geoffrey~E Hinton.
\newblock Adaptive mixtures of local experts.
\newblock \emph{Neural computation}, 3\penalty0 (1):\penalty0 79--87, 1991.

\bibitem[Jensen(1906)]{jensen1906fonctions}
J~L W~V Jensen.
\newblock {Sur les fonctions convexes et les in{\'{e}}galit{\'{e}}s entre les
  valeurs moyennes}.
\newblock \emph{Acta Mathematica}, 30\penalty0 (1):\penalty0 175--193, 1906.

\bibitem[Jiang and Tanner(1999)]{jiang1999hierarchical}
Wenxin Jiang and Martin~A Tanner.
\newblock Hierarchical mixtures-of-experts for exponential family regression
  models: approximation and maximum likelihood estimation.
\newblock \emph{Annals of Statistics}, pages 987--1011, 1999.

\bibitem[Khalili(2010)]{khalili2010new}
Abbas Khalili.
\newblock New estimation and feature selection methods in mixture-of-experts
  models.
\newblock \emph{Canadian Journal of Statistics}, 38\penalty0 (4):\penalty0
  519--539, 2010.

\bibitem[Khalili and Chen(2007)]{khalili2007variable}
Abbas Khalili and Jiahua Chen.
\newblock Variable selection in finite mixture of regression models.
\newblock \emph{Journal of the american Statistical association}, 102\penalty0
  (479):\penalty0 1025--1038, 2007.

\bibitem[Lloyd-Jones et~al.(2018)Lloyd-Jones, Nguyen, and
  McLachlan]{lloyd2018globally}
Luke~R Lloyd-Jones, Hien~D Nguyen, and Geoffrey~J McLachlan.
\newblock A globally convergent algorithm for lasso-penalized mixture of linear
  regression models.
\newblock \emph{Computational Statistics \& Data Analysis}, 119:\penalty0
  19--38, 2018.

\bibitem[Magnus and Neudecker(2019)]{magnus2019matrix}
Jan~R Magnus and Heinz Neudecker.
\newblock \emph{Matrix differential calculus with applications in statistics
  and econometrics}.
\newblock John Wiley \& Sons, 2019.

\bibitem[Mansuripur(1987)]{mansuripur1987introduction}
Masud Mansuripur.
\newblock \emph{{Introduction to information theory}}.
\newblock Prentice-Hall, Inc., 1987.

\bibitem[Massart(2007)]{massart_concentration_2007}
Pascal Massart.
\newblock \emph{Concentration {Inequalities} and {Model} {Selection}: {Ecole}
  d'{Eté} de {Probabilités} de {Saint}-{Flour} {XXXIII}-2003}.
\newblock Springer, 2007.

\bibitem[Massart and Meynet(2011)]{massart_lasso_2011}
Pascal Massart and Caroline Meynet.
\newblock {The Lasso as an $l_1$-ball model selection procedure}.
\newblock \emph{Electronic Journal of Statistics}, 5:\penalty0 669 -- 687,
  2011.

\bibitem[Massart and Meynet(2012)]{massart_rates_2012}
Pascal Massart and Caroline Meynet.
\newblock Some {Rates} of {Convergence} for the {Selected} {Lasso} {Estimator}.
\newblock In \emph{Algorithmic {Learning} {Theory}}, pages 17--33, Berlin,
  Heidelberg, 2012.
\newblock ISBN 978-3-642-34106-9.

\bibitem[Maugis and Michel(2011)]{maugis2011non}
Cathy Maugis and Bertrand Michel.
\newblock A non asymptotic penalized criterion for gaussian mixture model
  selection.
\newblock \emph{ESAIM: Probability and Statistics}, 15:\penalty0 41--68, 2011.

\bibitem[Maugis-Rabusseau and Michel(2013)]{maugis2013adaptive}
Cathy Maugis-Rabusseau and Bertrand Michel.
\newblock {Adaptive density estimation for clustering with Gaussian mixtures}.
\newblock \emph{ESAIM: Probability and Statistics}, 17:\penalty0 698--724,
  2013.

\bibitem[McLachlan and Peel(2000)]{mclachlan2000finite}
G~J McLachlan and D~Peel.
\newblock \emph{{Finite Mixture Models}}.
\newblock John Wiley \& Sons, 2000.

\bibitem[Mendes and Jiang(2012)]{mendes2012convergence}
Eduardo~F Mendes and Wenxin Jiang.
\newblock On convergence rates of mixtures of polynomial experts.
\newblock \emph{Neural computation}, 24\penalty0 (11):\penalty0 3025--3051,
  2012.

\bibitem[Meynet(2013)]{Meynet:2013aa}
C~Meynet.
\newblock {An $l_1$-oracle inequality for the Lasso in finite mixture Gaussian
  regression models}.
\newblock \emph{ESAIM: Probability and Statistics}, 17:\penalty0 650--671,
  2013.

\bibitem[Montuelle and Le~Pennec(2014)]{montuelle_mixture_2014}
Lucie Montuelle and Erwan Le~Pennec.
\newblock Mixture of {Gaussian} regressions model with logistic weights, a
  penalized maximum likelihood approach.
\newblock \emph{Electronic Journal of Statistics}, 8\penalty0 (1):\penalty0
  1661--1695, 2014.

\bibitem[Nguyen and Chamroukhi(2018)]{nguyen2018practical}
Hien~D Nguyen and Faicel Chamroukhi.
\newblock Practical and theoretical aspects of mixture-of-experts modeling: An
  overview.
\newblock \emph{Wiley Interdisciplinary Reviews: Data Mining and Knowledge
  Discovery}, 8\penalty0 (4):\penalty0 e1246, 2018.

\bibitem[Nguyen et~al.(2016)Nguyen, Lloyd-Jones, and
  McLachlan]{nguyen_universal_2016}
Hien~D Nguyen, Luke~R Lloyd-Jones, and Geoffrey~J McLachlan.
\newblock A universal approximation theorem for mixture-of-experts models.
\newblock \emph{Neural computation}, 28\penalty0 (12):\penalty0 2585--2593,
  2016.

\bibitem[Nguyen et~al.(2019)Nguyen, Chamroukhi, and
  Forbes]{nguyen2019approximation}
Hien~D Nguyen, Faicel Chamroukhi, and Florence Forbes.
\newblock {Approximation results regarding the multiple-output Gaussian gated
  mixture of linear experts model}.
\newblock \emph{Neurocomputing}, 366:\penalty0 208--214, 2019.
\newblock ISSN 0925-2312.

\bibitem[Nguyen et~al.(2021)Nguyen, Nguyen, Chamroukhi, and
  McLachlan]{nguyen2020approximationMoE}
Hien~D Nguyen, TrungTin Nguyen, Faicel Chamroukhi, and Geoffrey~John McLachlan.
\newblock {Approximations of conditional probability density functions in
  Lebesgue spaces via mixture of experts models}.
\newblock \emph{Journal of Statistical Distributions and Applications},
  8\penalty0 (1):\penalty0 13, 2021.

\bibitem[Nguyen et~al.(2023)Nguyen, Nguyen, and Ho]{nguyen_demystifying_2023}
Huy Nguyen, TrungTin Nguyen, and Nhat Ho.
\newblock Demystifying {Softmax} {Gating} {Function} in {Gaussian} {Mixture} of
  {Experts}.
\newblock In \emph{Thirty-seventh {Conference} on {Neural} {Information}
  {Processing} {Systems}}, 2023.

\bibitem[Nguyen et~al.(2024{\natexlab{a}})Nguyen, Akbarian, Nguyen, and
  Ho]{nguyen_general_2024}
Huy Nguyen, Pedram Akbarian, TrungTin Nguyen, and Nhat Ho.
\newblock A {General} {Theory} for {Softmax} {Gating} {Multinomial} {Logistic}
  {Mixture} of {Experts}.
\newblock In \emph{Forty-first {International} {Conference} on {Machine}
  {Learning}}, 2024{\natexlab{a}}.

\bibitem[Nguyen et~al.(2024{\natexlab{b}})Nguyen, Nguyen, Nguyen, and
  Ho]{nguyen_towards_2024}
Huy Nguyen, TrungTin Nguyen, Khai Nguyen, and Nhat Ho.
\newblock Towards {Convergence} {Rates} for {Parameter} {Estimation} in
  {Gaussian}-gated {Mixture} of {Experts}.
\newblock In \emph{Proceedings of {The} 27th {International} {Conference} on
  {Artificial} {Intelligence} and {Statistics}}, volume 238, pages 2683--2691,
  May 2024{\natexlab{b}}.

\bibitem[Nguyen et~al.(2020)Nguyen, Nguyen, Chamroukhi, and
  McLachlan]{nguyen_approximation_2020}
TrungTin Nguyen, Hien~D. Nguyen, Faicel Chamroukhi, and Geoffrey~J. McLachlan.
\newblock Approximation by finite mixtures of continuous density functions that
  vanish at infinity.
\newblock \emph{Cogent Mathematics \& Statistics}, 7\penalty0 (1):\penalty0
  1750861, January 2020.

\bibitem[Nguyen et~al.(2022{\natexlab{a}})Nguyen, Chamroukhi, Nguyen, and
  McLachlan]{nguyen_approximation_2022}
TrungTin Nguyen, Faicel Chamroukhi, Hien~D. Nguyen, and Geoffrey~J. McLachlan.
\newblock Approximation of probability density functions via location-scale
  finite mixtures in {Lebesgue} spaces.
\newblock \emph{Communications in Statistics - Theory and Methods}, pages
  1--12, May 2022{\natexlab{a}}.

\bibitem[Nguyen et~al.(2022{\natexlab{b}})Nguyen, Nguyen, Chamroukhi, and
  Forbes]{nguyen_non_asymptotic_2022}
TrungTin Nguyen, Hien~Duy Nguyen, Faicel Chamroukhi, and Florence Forbes.
\newblock A non-asymptotic approach for model selection via penalization in
  high-dimensional mixture of experts models.
\newblock \emph{Electronic Journal of Statistics}, 16\penalty0 (2):\penalty0
  4742 -- 4822, 2022{\natexlab{b}}.

\bibitem[Nguyen(2013)]{nguyen_convergence_2013}
XuanLong Nguyen.
\newblock Convergence of latent mixing measures in finite and infinite mixture
  models.
\newblock \emph{The Annals of Statistics}, 41\penalty0 (1):\penalty0 370--400,
  2013.

\bibitem[Norets(2010)]{norets_approximation_2010}
Andriy Norets.
\newblock Approximation of conditional densities by smooth mixtures of
  regressions.
\newblock \emph{The Annals of Statistics}, 38\penalty0 (3):\penalty0 1733 --
  1766, 2010.

\bibitem[Redner and Walker(1984)]{redner_mixture_1984}
Richard~A Redner and Homer~F Walker.
\newblock Mixture densities, maximum likelihood and the {EM} algorithm.
\newblock \emph{SIAM review}, 26\penalty0 (2):\penalty0 195--239, 1984.

\bibitem[Schwarz(1978)]{schwarz_estimating_1978}
Gideon Schwarz.
\newblock Estimating the dimension of a model.
\newblock \emph{The Annals of Statistics}, 6\penalty0 (2):\penalty0 461--464,
  1978.

\bibitem[Stadler et~al.(2010)Stadler, Buhlmann, and van~de Geer]{stadler2010l1}
N~Stadler, P~Buhlmann, and S~van~de Geer.
\newblock $l_1$-penalization for mixture regression models.
\newblock \emph{TEST}, 19:\penalty0 209--256, 2010.

\bibitem[Stone(1974)]{stone_cross_validatory_1974}
M.~Stone.
\newblock Cross-{Validatory} {Choice} and {Assessment} of {Statistical}
  {Predictions}.
\newblock \emph{Journal of the Royal Statistical Society: Series B
  (Methodological)}, 36\penalty0 (2):\penalty0 111--133, 1974.

\bibitem[Tibshirani(1996)]{tibshirani_regression_1996}
Robert Tibshirani.
\newblock Regression shrinkage and selection via the {Lasso}.
\newblock \emph{Journal of the Royal Statistical Society: Series B
  (Methodological)}, 58\penalty0 (1):\penalty0 267--288, 1996.

\bibitem[Van Der~Vaart and Wellner(1996)]{van1996weak}
AW~Van Der~Vaart and JA~Wellner.
\newblock \emph{Weak Convergence and Empirical Processes: With Applications to
  Statistics Springer Series in Statistics}, volume~58.
\newblock Springer, 1996.

\bibitem[Vapnik(1982)]{vapnik1982estimation}
Vladimir Vapnik.
\newblock \emph{Estimation of Dependences Based on Empirical Data (Springer
  Series in Statistics)}.
\newblock Springer-Verlag, 1982.

\bibitem[Wainwright(2019)]{wainwright2019high}
Martin~J Wainwright.
\newblock \emph{High-dimensional statistics: A non-asymptotic viewpoint},
  volume~48.
\newblock Cambridge University Press, 2019.

\bibitem[Williams and Rasmussen(2006)]{williams2006gaussian}
Christopher~K Williams and Carl~Edward Rasmussen.
\newblock \emph{{Gaussian processes for machine learning}}, volume~2.
\newblock MIT press Cambridge, MA, 2006.

\bibitem[Yuksel et~al.(2012)Yuksel, Wilson, and Gader]{yuksel2012twenty}
S~E Yuksel, J~N Wilson, and P~D Gader.
\newblock {Twenty Years of Mixture of Experts}.
\newblock \emph{IEEE Transactions on Neural Networks and Learning Systems},
  23\penalty0 (8):\penalty0 1177--1193, 2012.

\end{thebibliography}


\newpage
\appendix

\textbf{\Large{Supplement to
	``Non-asymptotic oracle inequalities for the Lasso in high-dimensional mixture of experts''}}
In this supplement, we provide proofs for  Theorem~\ref{thm_l1_Oracle_Inequality_Ball},and Propositions \ref{thm_l1_Oracle_Inequality_Ball_small} and \ref{thm_l1_Oracle_Inequality_Ball_large} in Appendices~\ref{l1-OracleInequality.Ball.Proof}, \ref{thm_l1_Oracle_Inequality_Ball_small.Proof} and \ref{thm_l1_Oracle_Inequality_Ball_large.Proof}, respectively. 
We then provide proofs for the remaining lemmas and provide further related technical results in Appendices~\ref{sec_proof_Lemma} and \ref{technicalResult}, respectively.

\section{Proof of main results} \label{sec_proof_Prop} 

First, it is important to note that the negative of the differential entropy (see, \eg~\citet[Chapter 9]{mansuripur1987introduction}) of the true unknown conditional density $s_0 \in S$, defined in \eqref{eq_Sm_bounded}, is finite, see more in Lemma \ref{lem_differentialEntropy_SGaME}, which is proved in Appendix \ref{sec_proof_lem_differentialEntropy_SGaME}.
\begin{lemma}
	\label{lem_differentialEntropy_SGaME}
	There exists a constant $H_{s_0} =\max\left\{0,\ln C_{s_0} \right\}$, where $C_{s_0} = \left(4\pi\right)^{-q/2} A_\Sigma^{q/2},$ \st 
	\begin{align}
		\max\left\{0,\sup_{x \in \cX} \int_{\R^q} \ln\left(s_0\left(y|x\right)\right) s_0\left(y|x\right) dy  \right\} \le H_{s_0} < \infty.\label{eq_differentialEntropy_finite}
	\end{align}
\end{lemma}

We then introduce some definitions and notations that we shall use in the proofs.
\subsection*{Additional notation}
For any measurable function $f:\R^q \rightarrow \R$
, consider its empirical norm
\begin{align*}
	\norm{f}_n:= \sqrt{\frac{1}{n}\sum_{i=1}^n f^2\left(Y_i| x_i\right)},
\end{align*}
and its conditional expectation	
\begin{align*}
	\E{Y|X = x}{f} := 
	\E{}{f\left(Y| X\right)| X = x} = \int_{\R^q}f(y| x)s_0(y| x)dy.
\end{align*}
%
Furthermore, we also define its empirical process
\begin{align}\label{eq_model_Sm1}
	P_n(f):=\frac{1}{n} \sum_{i=1}^n f\left(Y_i| x_i\right),
\end{align}
%
with expectation	
%
\begin{align}\label{eq_Pf_expectation_empirical}
	P(f) = \frac{1}{n} \sum_{i=1}^n  \E{Y_i|X_i = x_i}{f\left(Y_i| X_i\right)| X_i = x_i}=\frac{1}{n} \sum_{i=1}^n \int_{\R^q} f\left(y| x_i\right)s_0\left(y | x_i\right)dy, 
\end{align}
and the recentered process
\begin{align}\label{eq_model_Sm2}
	\nu_n(f):=P_n(f) - P(f) = \frac{1}{n} \sum_{i=1}^n \left[f\left(Y_i | x_i\right) - \int_{\R^q} f\left(y| x_i\right)s_0\left(y | x_i\right)dy\right].
\end{align}
For all $m\in \Ns$, recall that we consider the model
\begin{align*}
	S_m &= \left\{s_\psi \in S| \norm{\gamma}_1+\norm{\vect(\beta)}_1 \le m\right\} \equiv \left\{s_\psi \in S| \left\|\psi^{[1,2]}\right\|_{1} \le m\right\},\\
	F_m &= \left\{f_m = -\ln\left(\frac{s_m}{s_0}\right)= \ln(s_0)-\ln(s_m), s_m \in S_m \right\}.
\end{align*}

By using the basic properties of the infimum: for every $\epsilon > 0$, there exists $x_{\epsilon}\in A$, such that $x_\epsilon < \inf A+ \epsilon$. Then let $\delta_{\kl}>0$ for all $m\in \Ns$, and let $\eta_m \geq 0$. It holds that there exist two functions $\widehat{s}_m$ and $\overline{s}_m$ in $S_m$, such that
\begin{align}
	P_n\left(-\ln \widehat{s}_m\right) &\le \inf_{s_m \in S_m} P_n\left(-\ln s_m\right) + \eta_m, \text{and}	 \label{eq_model_Sm4}\\
	{\kl}_n\left(s_0,\overline{s}_m\right) &\le \inf_{s_m \in S_m} {\kl}_n\left(s_0,s_m\right) + \delta_{\kl}.\label{eq_model_Sm5}
\end{align}
Define 
\begin{align}\label{eq_model_Sm6}
	\widehat{f}_m := -\ln \left(\frac{\widehat{s}_m}{s_0}\right),\text{and }  \overline{f}_m:= -\ln \left(\frac{\overline{s}_m}{s_0}\right).
\end{align}
Let $\eta \geq 0$ and fix $m \in \Ns$. Further, define
\begin{align}\label{eq_model_Sm7}
	\widehat{\cM}(m) = \left\{m'\in\Ns| P_n\left(-\ln\widehat{s}_{m'}\right) + \pen(m') \le P_n\left(-\ln\widehat{s}_{m}\right) + \pen(m)+\eta \right\}.
\end{align}

\subsection{Proof of Theorem~\ref{thm_l1_Oracle_Inequality_Ball}} \label{l1-OracleInequality.Ball.Proof}
Let $M_n > 0$ and $\kappa \geq 148$. Assume that, for all $m \in \Ns$, the penalty function satisfies $\pen(m) = \lambda m$, with
\begin{align}\label{Extend.CarolineMeynet.eq.4.39.Gaussian}
	\lambda \geq \kappa\frac{KB_n}{\sqrt{n}} \left(q\ln n \sqrt{\ln(2p+1)}+1\right).
\end{align} 
We derive, from Propositions \ref{thm_l1_Oracle_Inequality_Ball_small} and \ref{thm_l1_Oracle_Inequality_Ball_large}, that any penalized likelihood estimator $\widehat{s}_{\widehat{m}}$ with $\widehat{m}$, satisfying
\begin{align*}
	-\frac{1}{n} \sum_{i=1}^n \ln\left(\widehat{s}_{\widehat{m}}\left(y_i| x_i\right)\right) + \pen(\widehat{m})\le \inf_{m \in \Ns} \left(-\frac{1}{n} \sum_{i=1}^n\ln\left(\widehat{s}_m\left(y_i| x_i\right)\right) + \pen(m) \right) + \eta,
\end{align*}
for some $\eta \geq 0$, yields
\begin{align}\label{eq_model_Sm_PMLE0.Gaussian}
	\E{}{{\kl}_n\left(s_0,\widehat{s}_{\widehat{m}}\right)}
	=&\E{}{{\kl}_n\left(s_0,\widehat{s}_{\widehat{m}}\right){\Indi}_{\cT}} +\E{}{{\kl}_n\left(s_0,\widehat{s}_{\widehat{m}}\right){\Indi}_{\cT^c}}\nn\\
	\le{}& \left(1+\kappa^{-1}\right)\inf_{m\in\Ns}\left(\inf_{s_m \in S_m} {\kl}_n\left(s_0,s_m\right)+ \pen(m) + \eta_m\nn\right)\\
	&  + \frac{302K^{3/2}qB_n}{\sqrt{n}}\left(1+\left(A_\gamma+  q A_\beta + \frac{ q\sqrt{q}}{a_{\Sigma}}\right)^2\right)+ \eta\nn\\
	&+
	\left(\frac{ e^{q/2-1}\pi^{q/2} }{A_\Sigma ^{q/2} } + H_{s_0}\right)  \sqrt{2 K n q A_\gamma} e^{-\frac{M_n^2-2M_n A_\beta}{4 A_\Sigma}}.
\end{align}
To obtain inequality \eqref{eq_l1_Oracle_Inequality_Ball}, it only remains to optimize the inequality \eqref{eq_model_Sm_PMLE0.Gaussian}, with respect $M_n$. Since the two terms depending on $M_n$, in \eqref{eq_model_Sm_PMLE0.Gaussian}, have opposite monotonicity with respect to $M_n$, we are looking for a value of $M_n$ such that these two terms are the same order with respect to $n$. Consider the positive solution $M_n = A_\beta + \sqrt{A_\beta^2+4 A_\Sigma \ln n}$ of the equation $\frac{X\left(X-2 A_\beta\right)}{4A_\Sigma} - \ln n = 0$. Then, on the one hand, 
\begin{align*}
	e^{-\frac{M_n^2-2M_n A_\beta}{4 A_\Sigma}}\sqrt{n} = e^{-\ln n} \sqrt{n} = \frac{1}{\sqrt{n}}.
\end{align*}
On the other hand, using the inequality $(a+b)^2 \le 2 (a^2+b^2)$, we have 
\begin{align*}
	B_n &= \max\left(A_\Sigma,1+KA_G\right)\left(1+ q\sqrt{q}\left(M_n+ A_\beta\right)^2 A_\Sigma\right)
	\nn\\
	&= \max\left(A_\Sigma,1+KA_G\right)\left(1+ q\sqrt{q} A_\Sigma \left( 2A_\beta + \sqrt{A_\beta^2+4 A_\Sigma \ln n}\right)^2 \right)\nn\\
	&\le \max\left(A_\Sigma,1+KA_G\right)\left(1+ 2q\sqrt{q} A_\Sigma \left( 5A^2_\beta + 4 A_\Sigma \ln n\right) \right),
\end{align*}
hence \eqref{eq_model_Sm_PMLE0.Gaussian} implies \eqref{eq_l1_Oracle_Inequality_Ball}.
\subsection{Proof of Proposition \ref{thm_l1_Oracle_Inequality_Ball_small}}\label{thm_l1_Oracle_Inequality_Ball_small.Proof}
For every $m' \in \widehat{\cM}(m)$, from \eqref{eq_model_Sm7}, \eqref{eq_model_Sm6}, and \eqref{eq_model_Sm4}, we obtain
\begin{align*}
	P_n\left(\widehat{f}_{m'}\right)+ \pen(m') &= P_n\left(\ln(s_0)-\ln \left(\widehat{s}_{m'}\right)\right)+ \pen(m')\nn \text{ (using \eqref{eq_model_Sm6})}\\
	&\le P_n\left(\ln(s_0)-\ln \left(\widehat{s}_{m}\right)\right)+ \pen(m) + \eta\nn \text{ (using \eqref{eq_model_Sm7})}\\
	&\le P_n\left(\ln(s_0)-\ln \left(\overline{s}_{m}\right)\right)+ \eta_m+ \pen(m) + \eta\nn\\
	& \quad  \text{ (using \eqref{eq_model_Sm4} with $\overline{s}_m$ in $S_m$ and linearity of $P_n$)}\nn\\
	&=P_n\left(\overline{f}_{m}\right)+ \pen(m) + \eta_m+ \eta  \text{ (using \eqref{eq_model_Sm6})}.
\end{align*}
By the definition of recentered process, $\nu_n\left(\cdot\right)$, in \eqref{eq_model_Sm2}, it holds that
\begin{align*}
	P\left(\widehat{f}_{m'}\right)+ \pen(m')\le P\left(\overline{f}_{m}\right)+ \pen(m)  + \nu_n\left(\overline{f}_m\right) - \nu_n\left(\widehat{f}_{m'}\right)+ \eta + \eta_m.
\end{align*}
Taking into account \eqref{Extend.CarolineMeynet.eq.2.2} and \eqref{eq_model_Sm1}, we obtain	
\begin{align*}
	{\kl}_n(s_0,\widehat{s}_{m'}) &= \frac{1}{n} \sum_{i=1}^n \int_{\R^q} \ln\left(\frac{s_0\left(y| x_i\right)}{\widehat{s}_{m'}\left(y| x_i\right)}\right)s_0\left(y| x_i\right) dy = \frac{1}{n} \sum_{i=1}^n \int_{\R^q} \widehat{f}_{m'}\left(y| x_i\right)s_0\left(y| x_i\right) dy\nn \text{ (using \eqref{eq_model_Sm6})}\\
	&=
	P\left(\widehat{f}_{m'}\right) \text{ (using \eqref{eq_Pf_expectation_empirical})}.
\end{align*}
Similarly, we also obtain ${\kl}_n(s_0,\overline{s}_{m})
=   P\left(\overline{f}_{m}\right)$.
Hence, \eqref{eq_model_Sm5} implies that
\begin{align}\label{eq_model_Sm8}
	{\kl}_n(s_0,\widehat{s}_{m'})+ \pen(m')
	&\le {\kl}_n(s_0,\overline{s}_{m})+ \pen(m)  + \nu_n\left(\overline{f}_m\right) - \nu_n\left(\widehat{f}_{m'}\right)+ \eta + \eta_m\nn\\
	&\le\inf_{s_m \in S_m} {\kl}_n\left(s_0,s_m\right)+ \pen(m)  + \nu_n\left(\overline{f}_m\right) \nn\\
	& \quad - \nu_n\left(\widehat{f}_{m'}\right) + \eta_m  + \delta_{\kl}+ \eta.
\end{align}
All that remains is to control the deviation of $ - \nu_n\left(\widehat{f}_{m'}\right) = \nu_n\left(-\widehat{f}_{m'}\right)$. To handle the randomness of $\widehat{f}_{m'}$, we shall control the deviation of $\sup_{f_{m'} \in F_{m'}} \nu_n\left(-f_{m'}\right)$, since $\widehat{f}_{m'}\in F_{m'}$. Such control is provided by Lemma \ref{lem_control_deviation}.
From \eqref{eq_model_Sm8} and \eqref{eq_model_Sm_LL1.Gaussian}, we derive that on the event $\cT$,  for all $m \in \Ns$, and $t>0$, with probability larger than $1-e^{-t}$,
\begin{align}
	{\kl}_n(s_0,\widehat{s}_{m'})+ \pen(m')
	\le{}&\inf_{s_m \in S_m} {\kl}_n\left(s_0,s_m\right)+ \pen(m)  + \nu_n\left(\overline{f}_m\right) - \nu_n\left(\widehat{f}_{m'}\right) + \eta_m  + \delta_{\kl}+ \eta \nn\\
	\le{}&\inf_{s_m \in S_m} {\kl}_n\left(s_0,s_m\right)+ \pen(m)  + \nu_n\left(\overline{f}_m\right) + \eta_m  + \delta_{\kl}+ \eta\nn\\
	&+\frac{4KB_n}{\sqrt{n}} \left[37 q \Delta_{m'} + \sqrt{2}\left(A_\gamma+  q A_\beta + \frac{ q\sqrt{q}}{a_{\Sigma}}\right)\sqrt{t}\right]\nn\\
	\le&\inf_{s_m \in S_m} {\kl}_n\left(s_0,s_m\right)+ \pen(m)  + \nu_n\left(\overline{f}_m\right) + \eta_m  + \delta_{\kl}+ \eta\nn\\
	&+\frac{4KB_n}{\sqrt{n}} \left[37 q \Delta_{m'} + \frac{1}{2}\left(A_\gamma+  q A_\beta + \frac{ q\sqrt{q}}{a_{\Sigma}}\right)^2+t\right], \text{  if }m'\in \widehat{\cM}(m). \label{eq_model_Sm_LL2}
\end{align}
Here, we use the fact that $\widehat{f}_{m'} \in F_{m'}$ and get the last inequality using the fact that $2ab \le a^2 + b^2 $ for $b = \sqrt{t}$, and $a =  \left(A_\gamma+  q A_\beta+ \frac{ q\sqrt{q}}{a_{\Sigma}}\right)/\sqrt{2}$.

It remains to sum up the tail bounds \eqref{eq_model_Sm_LL2} over all possible values of $m \in \Ns$ and $m'\in \widehat{\cM}(m)$. To get an inequality valid on a set of high probability, we need to adequately choose the value of the parameter $t$, depending on $m\in \Ns$ and $m'\in \widehat{\cM}(m)$. Let $z > 0$, for all $m \in \Ns$ and $m' \in \widehat{\cM}(m)$, and apply \eqref{eq_model_Sm_LL2} to obtain $t = z+m+m'$. Then, on the event $\cT$, for all $m \in \Ns$, 
with probability larger than $1-e^{-\left(z+m+m'\right)}$,
\begin{align}
	{\kl}_n(s_0,\widehat{s}_{m'})+ \pen(m')
	\le{}&\inf_{s_m \in S_m} {\kl}_n\left(s_0,s_m\right)+ \pen(m)  + \nu_n\left(\overline{f}_m\right) + \eta_m  + \delta_{\kl}+ \eta\nn\\
	&+\frac{4KB_n}{\sqrt{n}} \left[37 q \Delta_{m'} + \frac{1}{2}\left(A_\gamma+  q A_\beta + \frac{ q\sqrt{q}}{a_{\Sigma}}\right)^2+\left(z+m+m'\right)\right],\nn\\
	& \hspace{-4.5cm}\text{if $m'\in \widehat{\cM}(m)$}. \label{eq_model_Sm_LL3.Gaussian}
\end{align}
Here, \eqref{eq_model_Sm_LL3.Gaussian} is equivalent to
\begin{align}
	{\kl}_n(s_0,\widehat{s}_{m'}) - \nu_n\left(\overline{f}_m\right)
	\le{}&\inf_{s_m \in S_m} {\kl}_n\left(s_0,s_m\right)+ \left[\pen(m)+\frac{4KB_n}{\sqrt{n}}m\right] + \eta_m  + \delta_{\kl}+ \eta\nn\\
	&+ \left[\frac{4KB_n}{\sqrt{n}} \left(37 q\Delta_{m'}+m'\right) - \pen(m') \right]\nn\\
	& + \frac{4KB_n}{\sqrt{n}} \left[ \frac{1}{2}\left(A_\gamma+  q A_\beta + \frac{ q\sqrt{q}}{a_{\Sigma}}\right)^2+z\right].\label{eq_model_Sm_LL3.2.Gaussian}
\end{align}

Note that with probability larger than $1-e^{-z}$, \eqref{eq_model_Sm_LL3.Gaussian} holds simultaneously for all $m \in \Ns$ and $m' \in \widehat{\cM}(m)$. Indeed, by defining the event $${\cap}_{(m,m')\in \Ns \times \widehat{\cM}(m)}\Omega_{m,m'}=\left\{w: w\in \Omega \text{ such that the event in \eqref{eq_model_Sm_LL3.Gaussian} holds}\right\},$$
%
it holds that, on the event $\cT$, 
\begin{eqnarray*}
	\bP\left(\cap_{(m,m') \in \Ns \times \widehat{\cM}(m)}\Omega_{m,m'}\right) &=& 1- \bP\left(\cup_{(m,m') \in \Ns \times \widehat{\cM}(m)}\Omega^C_{m,m'}\right)\nn\\
	&\geq& 1- \sum_{(m,m') \in \Ns \times \widehat{\cM}(m)}\bP\left(\Omega^C_{m,m'}\right)\nn\\
	&\geq&  1-\sum_{(m,m') \in \Ns \times \widehat{\cM}(m)}e^{-\left(z+m+m'\right)}\nn\\
	&\geq&  1-\sum_{(m,m') \in \Ns \times \Ns}e^{-\left(z+m+m'\right)}\nn \\
	&=&  1-e^{-z}\left(\sum_{m \in \Ns }e^{-m}\right)^2
	\ge 1- e^{-z},
\end{eqnarray*}
where we get the last inequality by using the the geometric series
\begin{align*}
	\sum_{m=1}^{\infty}\left(e^{-1}\right)^m =\sum_{m=0}^{\infty}\left(e^{-1}\right)^m -1 = \frac{1}{1-e^{-1}}-1=\frac{e}{e-1} -1=\frac{1}{e-1} < 1.
\end{align*} 

Taking into account \eqref{eq_Delta_m}, we get
\begin{align}
	{\kl}_n(s_0,\widehat{s}_{m'}) - \nu_n\left(\overline{f}_m\right)
	\le{}&\inf_{s_m \in S_m} {\kl}_n\left(s_0,s_m\right)+ \left[\pen(m)+\frac{4KB_n}{\sqrt{n}}m\right] + \eta_m  + \delta_{\kl}+ \eta\nn\\
	&+ \left[\frac{4KB_n}{\sqrt{n}} \left(37q\ln n\sqrt{\ln(2p+1)}+1\right)m' - \pen(m') \right]\nn\\
	& + \frac{4KB_n}{\sqrt{n}} \Bigg[\frac{1}{2}\left(A_\gamma+  q A_\beta + \frac{ q\sqrt{q}}{a_{\Sigma}}\right)^2 + 74q\sqrt{K} \left(A_\gamma+  q A_\beta + \frac{ q\sqrt{q}}{a_{\Sigma}}\right)+z \Bigg]. \label{eq_model_Sm_LL4.Gaussian}
\end{align}
Now, let $\kappa \geq 1$ and assume that $\pen(m)=\lambda m$, for all $ m \in \Ns$ with 
\begin{align}\label{eq.initialConditionLambda}
	\lambda \geq \kappa\frac{4KB_n}{\sqrt{n}} \left(37q\ln n \sqrt{\ln(2p+1)}+1\right).
\end{align}
Then, \eqref{eq_model_Sm_LL4.Gaussian} implies
\begin{align*}
	{\kl}_n(s_0,\widehat{s}_{m'}) - \nu_n\left(\overline{f}_m\right)
	\le{}&\inf_{s_m \in S_m} {\kl}_n\left(s_0,s_m\right)+ \left[\lambda m+\frac{4KB_n}{\sqrt{n}}m\right] + \eta_m  + \delta_{\kl}+ \eta\nn\\
	&+ \left[\underbrace{\frac{4KB_n}{\sqrt{n}} \left(37q\ln n \sqrt{\ln(2p+1)}+1\right)}_{\le \lambda \kappa^{-1}}m' - \lambda m' \right]\nn\\
	& + \frac{4KB_n}{\sqrt{n}} \Bigg[\frac{1}{2}\left(A_\gamma+  q A_\beta + \frac{ q\sqrt{q}}{a_{\Sigma}}\right)^2 + 74q\sqrt{K} \left(A_\gamma+  q A_\beta + \frac{ q\sqrt{q}}{a_{\Sigma}}\right)+z \Bigg]\nn\\
	\le{}&\inf_{s_m \in S_m} {\kl}_n\left(s_0,s_m\right)+ \left[\pen(m)+\underbrace{\frac{4KB_n}{\sqrt{n}}m}_{\le \kappa^{-1}\pen(m)}\right] + \eta_m  + \delta_{\kl}+ \eta\nn\\
	&+ \underbrace{\left[\lambda \kappa^{-1}m' - \lambda m' \right]}_{\le 0}\nn\\
	& + \frac{4KB_n}{\sqrt{n}} \Bigg[\frac{1}{2}\left(A_\gamma+  q A_\beta + \frac{ q\sqrt{q}}{a_{\Sigma}}\right)^2 + 74q\sqrt{K} \left(A_\gamma+  q A_\beta + \frac{ q\sqrt{q}}{a_{\Sigma}}\right)+z \Bigg]\nn\\
	\le{}&\inf_{s_m \in S_m} {\kl}_n\left(s_0,s_m\right)+ \left(1+\kappa^{-1}\right)\pen(m) + \eta_m  + \delta_{\kl}+ \eta\nn\\
	& + \frac{4KB_n}{\sqrt{n}} \Bigg[\frac{1}{2}\left(A_\gamma+  q A_\beta + \frac{ q\sqrt{q}}{a_{\Sigma}}\right)^2 + 74q\sqrt{K} \left(A_\gamma+  q A_\beta + \frac{ q\sqrt{q}}{a_{\Sigma}}\right)+z \Bigg].\nn\\
\end{align*}
Next, using the inequality $2ab \le \beta^{-1}a^2 + \beta^{-1}b^2$ for $a = \sqrt{K}$, $b = K\left(A_\gamma+  q A_\beta + \frac{ q\sqrt{q}}{a_{\Sigma}}\right)$, and $\beta = \sqrt{K}$, and the fact that $K \le K^{3/2}$, for all $K \in \Ns$, it follows that
\begin{align}
	&{\kl}_n(s_0,\widehat{s}_{m'}) - \nu_n\left(\overline{f}_m\right)
	\le \inf_{s_m \in S_m} {\kl}_n\left(s_0,s_m\right)+ \left(1+\kappa^{-1}\right)\pen(m) + \eta_m  + \delta_{\kl}+ \eta\nn\\
	& + \frac{4B_n}{\sqrt{n}} \Bigg[ \frac{qK^{3/2}}{2}\left( A_\gamma+  q A_\beta + \frac{ q\sqrt{q}}{a_{\Sigma}}\right)^2+\underbrace{74q\sqrt{K}K\left(A_\gamma+  q A_\beta + \frac{ q\sqrt{q}}{a_{\Sigma}}\right)}_{37q\times 2ab}+Kz\Bigg]\nn\\
	\le{}&\inf_{s_m \in S_m} {\kl}_n\left(s_0,s_m\right)+ \left(1+\kappa^{-1}\right)\pen(m) + \eta_m  + \delta_{\kl}+ \eta\nn\\
	& + \frac{4B_n}{\sqrt{n}} \left[ 37qK^{1/2} + \frac{75qK^{3/2}}{2}\left(A_\gamma+  q A_\beta + \frac{ q\sqrt{q}}{a_{\Sigma}}\right)^2+ Kz\right].
	\label{eq_model_Sm_LL5.Gaussian}
\end{align}
By \eqref{eq_model_Sm_PMLE} and \eqref{eq_model_Sm7}, $\widehat{m}$ belongs to $\widehat{\cM}(m)$, for all $m\in\Ns$,
so we deduce from \eqref{eq_model_Sm_LL5.Gaussian} that on the event $\cT$, for all $z >0$, with probability greater than $1-e^{-z}$,
\begin{align}\label{eq_model_Sm_LL6.Gaussian}
	{\kl}_n(s_0,\widehat{s}_{\widehat{m}}) - \nu_n\left(\overline{f}_m\right)
	\le{}& \inf_{m\in\Ns}\left(\inf_{s_m \in S_m} {\kl}_n\left(s_0,s_m\right)+ \left(1+\kappa^{-1}\right)\pen(m) + \eta_m\nn\right) + \eta +  \delta_{\kl}\\
	& + \frac{4B_n}{\sqrt{n}} \left[ 37qK^{1/2} + \frac{75qK^{3/2}}{2}\left(A_\gamma+  q A_\beta + \frac{ q\sqrt{q}}{a_{\Sigma}}\right)^2+ Kz\right].
\end{align}

Note that for any non-negative random variable $Z$ and any $a > 0,\E{}{Z} = a \int_{z \ge 0}\bP(Z > az)dz$. 
Indeed, if we let $t = az$, then $dz= a dt$ and
\begin{align}
	a \int_{z \ge 0}\bP(Z > az)dz &=  a\int_{0}^\infty \int_{az}^\infty f_Z(u)du dz = \int_{0}^\infty \int_{t}^\infty f_Z(u)du dt = \int_{0}^\infty \int_{0}^u f_Z(u)dt du \nn\\
	& = \int_{0}^\infty f_Z(u) \int_{0}^u dt du = \int_{0}^\infty f_Z(u) u du = \E{}{Z}. \label{eq_expectation_nonnegativeRV}
\end{align}
Then, we define the following random variable \wrt~the random response $Y_{[n]} := \left(Y_i\right)_{i \in [n]}$:
\begin{align*}
	Z &:= {\kl}_n(s_0,\widehat{s}_{\widehat{m}}) - \nu_n\left(\overline{f}_m\right) - \left[\inf_{m\in\Ns}\left(\inf_{s_m \in S_m} {\kl}_n\left(s_0,s_m\right)+ \left(1+\kappa^{-1}\right)\pen(m) + \eta_m\nn\right) + \eta +  \delta_{\kl}\right]\nn\\
	& \quad - \frac{4B_n}{\sqrt{n}} \left[ 37qK^{1/2} + \frac{75qK^{3/2}}{2}\left(A_\gamma+  q A_\beta + \frac{ q\sqrt{q}}{a_{\Sigma}}\right)^2\right].
\end{align*}
Then by \eqref{eq_model_Sm_LL6.Gaussian}, on the even $\cT$, it holds that $P(Z \le az) \ge 1-e^{-z}$ and if $Z \le 0$ then $P(Z < a z) = 1 \ge 1-e^{-z}$, for all $z > 0$, where $a =  \frac{4B_nK}{\sqrt{n}} > 0$. Therefore, it is sufficient to consider $Z \ge 0$ and it holds that $P(Z > az|\cT) \le e^{-z}$. In this case, by \eqref{eq_expectation_nonnegativeRV} and the fact that $\bP(\cT) \le 1$, it holds that
\begin{align}
	\E{Y_{[n]}}{Z{\Indi}_{\cT}} &=\bP(\cT)\E{Y_{[n]}}{Z | \cT} \le \E{Y_{[n]}}{Z | \cT} 
	\le a \int_{z \ge 0}e^{-z}dz = a.\label{eq_upperBound_Exp_Z}
\end{align}
%

Then, by integrating  \eqref{eq_model_Sm_LL6.Gaussian} over $z>0$ using \eqref{eq_upperBound_Exp_Z}, the fact that 
\begin{align}
	\E{Y_{[n]}}{\nu_n\left(\overline{f}_m\right)} &=\E{Y_{[n]}}{P_n\left(\overline{f}_m\right)} - \E{Y_{[n]}}{P\left(\overline{f}_m\right) } = 0,
\end{align}
$\delta_{\kl} >0$ can be chosen arbitrary small, and $\E{Y_{[n]}}{{\Indi}_{\cT}} = \bP\left(\cT\right) \le 1$,  we obtain that
\begin{align}\label{eq_model_Sm_LL6.1.Gaussian}
	\E{}{{\kl}_n(s_0,\widehat{s}_{\widehat{m}}){\Indi}_{\cT}}
	\le{}& \left[\inf_{m\in\Ns}\left(\inf_{s_m \in S_m} {\kl}_n\left(s_0,s_m\right)+ \left(1+\kappa^{-1}\right)\pen(m) + \eta_m\nn\right) + \eta\right]\E{Y_{[n]}}{{\Indi}_{\cT}}\nn\\
	&+ \frac{4B_n}{\sqrt{n}} \left[ 37qK^{1/2} + \frac{75qK^{3/2}}{2}\left(A_\gamma+  q A_\beta + \frac{ q\sqrt{q}}{a_{\Sigma}}\right)^2+ K\right]\E{Y_{[n]}}{{\Indi}_{\cT}}\nn\\
	\le{}& \inf_{m\in\Ns}\left(\inf_{s_m \in S_m} {\kl}_n\left(s_0,s_m\right)+ \left(1+\kappa^{-1}\right)\pen(m) + \eta_m\nn\right) + \eta\\
	&+ \frac{4B_n}{\sqrt{n}} \left[ 37qK^{3/2} + \frac{75qK^{3/2}}{2}\left(A_\gamma+  q A_\beta + \frac{ q\sqrt{q}}{a_{\Sigma}}\right)^2+ qK^{3/2}\right]\nn\\
	\le{}& \inf_{m\in\Ns}\left(\inf_{s_m \in S_m} {\kl}_n\left(s_0,s_m\right)+ \left(1+\kappa^{-1}\right)\pen(m) + \eta_m\nn\right) + \eta\\
	&  + \frac{302K^{3/2}qB_n}{\sqrt{n}}\left(1+\left(A_\gamma+  q A_\beta + \frac{ q\sqrt{q}}{a_{\Sigma}}\right)^2\right).
\end{align}

\subsection{Proof of Proposition \ref{thm_l1_Oracle_Inequality_Ball_large}}\label{thm_l1_Oracle_Inequality_Ball_large.Proof}
By the Cauchy-Schwarz inequality,
\begin{align} \label{eq_model_Sm_LL7}
	\E{}{{\kl}_n\left(s_0,\widehat{s}_{\widehat{m}}\right){\Indi}_{\cT^C}} \le \sqrt{\E{}{{\kl}_n^2\left(s_0,\widehat{s}_{\widehat{m}}\right)}}\sqrt{\E{}{{\Indi}^2_{\cT^C}}} = \sqrt{\E{}{{\kl}_n^2\left(s_0,\widehat{s}_{\widehat{m}}\right)}}\sqrt{\bP\left(\cT^C\right)}.
\end{align}
We seek to bound the two terms of the right-hand side of \eqref{eq_model_Sm_LL7}.

For the first term, let us bound $\kl\left(s_0\left(\cdot| x\right),s_\psi\left(\cdot| x\right)\right)$, for all $s_\psi \in S$ and $x \in \cX$. Let $s_\psi \in S$ and $x \in \cX$. 
%
%
Then, we obtain
\begin{align}\label{eq_model_Sm_LL8}
	\kl\left(s_0\left(\cdot | x\right),s_\psi\left(\cdot | x\right)\right) ={}& \int_{\R^q} \ln\left(\frac{s_0(y| x)}{s_\psi(y| x)}\right)s_0(y| x)dy\nn\\
	={}& \int_{\R^q} \ln\left(s_0(y| x)\right)s_0(y| x)dy - \int_{\R^q} \ln\left(s_\psi(y| x)\right)s_0(y| x)dy\nn\\
	\le{}& -\int_{\R^q} \ln\left(s_\psi(y| x)\right)s_0(y| x)dy+H_{s_0}, \forall x \in \cX \left(\text{using \eqref{eq_differentialEntropy_finite}}\right). 
\end{align}
Since
\begin{align*}
	a_G:=\frac{\exp\left(-A_{\gamma}\right)}{\sum_{l=1}^K\exp\left(A_{\gamma}\right)} \le  \sup_{x\in \cX, \gamma \in  \widetilde{\Gamma}} \frac{\exp\left(\gamma_{k0} + \gamma_k^\top x\right)}{\sum_{l=1}^K\exp\left(\gamma_{l0} + \gamma_l^\top x\right)} \le \frac{\exp\left(A_{\gamma}\right)}{\sum_{l=1}^K\exp\left(-A_{\gamma}\right)}=:A_G,
\end{align*}
there exists deterministic positive constants $a_G,A_G$, such that
\begin{align}\label{eq.defBoundedGating}
	a_G \le \sup_{x\in \cX, \gamma \in  \widetilde{\Gamma}}{g}_{k}\left(x;\gamma\right) \le A_G.
\end{align}
Here,  the softmax gating function ${g}_{k}\left(x;\gamma\right)$ is described as 
\begin{align}\label{eq.def.SoftmaxGatingParameter}
	{g}_{k}\left(x;\gamma\right)= \frac{\exp\left(w_k(x)\right)}{\sum_{l=1}^K\exp\left(w_l(x)\right)}, w_k(x) = \gamma_{k0} + \gamma_k^\top x, \gamma = \left(\gamma_{k0},\gamma_k^\top \right)_{k \in [K]} \in \Gamma = \R^{(p+1)K}.
\end{align}
Thus, for all $y \in \R^q$, 
\begin{align}
	&\ln\left(s_\psi(y| x)\right)s_0(y| x)\nn\\
	&\geq\ln \left[\sum_{k=1}^K  \frac{a_G\det(\Sigma_{k}^{-1})^{1/2}}{\left(2 \pi \right)^{q/2}} \exp\left(-\left(y^\top  \Sigma_{k}^{-1} y + \left(\beta_{k0}+ \beta_k x\right)^\top  \Sigma_{k}^{-1}\left(\beta_{k0}+ \beta_k x\right)\right)\right)\right] \nn \\
	& \hspace{0.5cm} \times \sum_{k=1}^K  \frac{a_G\det(\Sigma_{0,k}^{-1})^{1/2}}{\left(2 \pi \right)^{q/2} } \exp\left(-\left(y^\top  \Sigma_{0,k}^{-1} y + \left(\beta_{0,k0}+\beta_{0,k}x\right)^\top  \Sigma_{0,k}^{-1} \left(\beta_{0,k0}+\beta_{0,k}x\right)\right)\right)\nn\\
	& \left(\text{using \eqref{eq.defBoundedGating} and $-(a-b)^\top  A (a-b)/2 \geq -(a^\top A a+b^\top A b)$,~\eg $a = y, b = \beta_{k0}+ \beta_k x$, $A = \Sigma_{k}$}\right)\nn\\
	&\geq\ln \left[\sum_{k=1}^K  \frac{a_G a^{q/2}_{\Sigma}}{\left(2 \pi \right)^{q/2}} \exp\left(-\left(y^\top  \Sigma_{k}^{-1} y + \left(\beta_{k0}+ \beta_k x\right)^\top  \Sigma_{k}^{-1}\left(\beta_{k0}+ \beta_k x\right)\right)\right)\right] \nn \\
	& \hspace{0.5cm} \times \sum_{k=1}^K  \frac{a_G a^{q/2}_{\Sigma}}{\left(2 \pi \right)^{q/2} }  \exp\left(-\left(y^\top  \Sigma_{0,k}^{-1} y + \left(\beta_{0,k0}+\beta_{0,k}x\right)^\top  \Sigma_{0,k}^{-1} \left(\beta_{0,k0}+\beta_{0,k}x\right)\right)\right) \left(\text{using \eqref{eq_defBoundedParameters}}\right)\nn\\
	&\geq\ln \left[K  \frac{a_G a^{q/2}_{\Sigma}}{\left(2 \pi \right)^{q/2}} \exp\left(-\left(y^\top  y + qA^2_\beta\right) A_\Sigma \right)\right] \times K  \frac{a_G a^{q/2}_{\Sigma}}{\left(2 \pi \right)^{q/2} } \exp\left(-\left(y^\top  y + qA^2_\beta\right) A_\Sigma \right) \left(\text{using \eqref{eq_defBoundedParameters}}\right) \label{eq_model_Sm_LL9},
\end{align}
where, in the last inequality, we use the fact that for all $u \in \R^q$. By using the eigenvalue decomposition of $\Sigma_1 = P^\top  D P$,
\begin{align*}
	\left|u^\top \Sigma_1 u\right| = \left|u^\top  P^\top D P u\right| \le \left\| P u\right\|_2  \le  M(D)\left\| Pu\right\|^2_2 \le A_{\Sigma} \left\| u\right\|^2_2 \le A_\Sigma q\left\| u\right\|^2_\infty,
\end{align*}
where in the last inequality, we used the fact that \eqref{eq.vector.infinity2.equiv}. Therefore, setting $u=\sqrt{2A_\Sigma}y$ and $h(t) = t \ln t$, for all $t \in \R$, and noticing that $h(t) \geq h\left(e^{-1}\right) = -e^{-1}$, for all $t \in \R$, and from \eqref{eq_model_Sm_LL8} and \eqref{eq_model_Sm_LL9}, we get that
\begin{align}\label{Extend.CarolineMeynet.eq.4.30}
	&\kl\left(s_0\left(\cdot| x\right),s_\psi\left(\cdot| x\right)\right) - H_{s_0}\nn\\
	&\le -\int_{\R^q} \left[\ln \left[K  \frac{a_\gamma a^{q/2}_{\Sigma}}{\left(2 \pi \right)^{q/2}} \exp\left(-\left(y^\top  y + qA^2_\beta\right) A_\Sigma \right)\right] K  \frac{a_\gamma a^{q/2}_{\Sigma}}{\left(2 \pi \right)^{q/2} } \exp\left(-\left(y^\top  y + qA^2_\beta\right) A_\Sigma \right) \right)dy\nn\\
	&= - \frac{K  a_\gamma a^{q/2}_{\Sigma} e^{-qA^2_\beta A_\Sigma}}{\left(2 A_\Sigma \right)^{q/2} }\int_{\R^q} \left[\ln \left(K  \frac{a_\gamma a^{q/2}_{\Sigma}}{\left(2 \pi \right)^{q/2}}\right) - qA^2_\beta A_\Sigma - \frac{u^\top u}{2} \right] \frac{e^{-\frac{u^\top u}{2}}}{\left(2 \pi \right)^{q/2}} du\nn\\
	&= - \frac{K  a_\gamma a^{q/2}_{\Sigma} e^{-qA^2_\beta A_\Sigma}}{\left(2 A_\Sigma \right)^{q/2} } \E{U}{\left[\ln \left(K  \frac{a_\gamma a^{q/2}_{\Sigma}}{\left(2 \pi \right)^{q/2}}\right) - qA^2_\beta A_\Sigma - \frac{U^\top U}{2} \right]} \left(\text{with } U \sim \cN_q(0,{\id}_q)\right)\nn\\
	&= - \frac{K  a_\gamma a^{q/2}_{\Sigma} e^{-qA^2_\beta A_\Sigma}}{\left(2 A_\Sigma \right)^{q/2} } \left[\ln \left(K  \frac{a_\gamma a^{q/2}_{\Sigma}}{\left(2 \pi \right)^{q/2}}\right) - qA^2_\beta A_\Sigma - \frac{q}{2} \right]\nn\\
	&= - \frac{K  a_\gamma a^{q/2}_{\Sigma} e^{-qA^2_\beta A_\Sigma - \frac{q}{2}}}{\left(2 \pi \right)^{q/2}\left( A_\Sigma \right)^{q/2} } e^{q/2}\pi^{q/2} \ln \left(\frac{Ka_\gamma a^{q/2}_{\Sigma}e^{- qA^2_\beta A_\Sigma - \frac{q}{2}}}{\left(2 \pi \right)^{q/2}}\right) \le \frac{ e^{q/2-1}\pi^{q/2} }{A_\Sigma ^{q/2} },
\end{align}
where we used the fact that $t\ln(t) \geq -e^{-1}$, for all $t \in \R$.

Then, for all $s_\psi \in S$,
\begin{align*}
	{\kl}_n\left(s_0,s_\psi\right) = \frac{1}{n}\sum_{i=1}^n{\kl}\left(s_0\left(\cdot| x_i\right),s_\psi\left(.| x_i\right)\right) \le\frac{ e^{q/2-1}\pi^{q/2} }{A_\Sigma ^{q/2} }+H_{s_0},
\end{align*}
and note that $\widehat{s}_{\widehat{m}} \in S$, thus
\begin{align}\label{Extend.CarolineMeynet.eq.4.31}
	\sqrt{\E{}{{\kl}_n^2\left(s_0,\widehat{s}_{\widehat{m}}\right)}} \le \frac{ e^{q/2-1}\pi^{q/2} }{A_\Sigma ^{q/2} } + H_{s_0}.
\end{align}
We now provide an upper bound for $\bP\left(\cT^C\right)$:
\begin{align}\label{Extend.CarolineMeynet.eq.4.32}
	%
	\bP\left(\cT^C\right) \le \sum_{i=1}^n \bP\left(\left\|Y_i\right\|_\infty>M_n\right).
\end{align}
For all $i \in [n]$, 
$$Y_i|x_i \sim \sum_{k=1}^K {g}_{k}\left(x_i;\gamma\right)\cN_q\left(\beta_{k0}+\beta_kx_i,\Sigma_k\right),$$ 
so we see from \eqref{Extend.CarolineMeynet.eq.4.32} that we need to provide an upper bound on $\bP\left(\left|Y_x\right|>M_n\right)$, with $$Y_x\sim \sum_{k=1}^K {g}_{k}\left(x;\gamma\right) \cN_q\left(\beta_{k0}+\beta_kx,\Sigma_{k}\right), x \in \cX.$$
First, using Chernoff's inequality for a centered Gaussian variable (see Lemma \ref{lem.Chernoff.Inequality}), and the fact that $\psi$ belongs to the bounded space $\widetilde{\Psi}$ (defined by \eqref{eq_defBoundedParameters}), and that $\sum_{k=1}^K {g}_{k}\left(x;\gamma\right) = 1$, we get
\begin{align} \label{Extend.CarolineMeynet.eq.4.36}
	&\bP\left(\left\|Y_x\right\|_\infty>M_n\right) \nn\\ 
	\nn\\
	&= \sum_{k=1}^K \frac{{g}_{k}\left(x;\gamma\right)}{\left(2 \pi \right)^{q/2} \det(\Sigma_k)^{1/2}}\int_{\left\{\left\|y\right\|_\infty>M_n\right\}}   \exp\left(-\frac{\left(y-\left(\beta_{k0}+\beta_kx\right)\right)^\top  \Sigma_k^{-1} \left(y-\left(\beta_{k0}+\beta_kx\right)\right)}{2}\right)dy\nn\\
	&= \sum_{k=1}^K{g}_{k}\left(x;\gamma\right)\bP\left(\left\|Y_{x,k}\right\|_\infty>M_n\right)
	\le \sum_{k=1}^K{g}_{k}\left(x;\gamma\right)\sum_{z=1}^q \bP\left(\left|\left[Y_{x,k}\right]_z\right|>M_n\right) \nn\\
	&= \sum_{k=1}^K{g}_{k}\left(x;\gamma\right)\sum_{z=1}^q \left(\bP\left(\left[Y_{x,k}\right]_z<-M_n\right) + \bP\left(\left[Y_{x,k}\right]_z>M_n\right)\right) \nn\\
	&= \sum_{k=1}^K{g}_{k}\left(x;\gamma\right)\sum_{z=1}^q \left(\bP\left(U > \frac{M_n-\left[\beta_{k0}+\beta_kx\right]_z}{\left[\Sigma_{k}\right]^{1/2}_{z,z}}\right) + \bP\left(U < \frac{-M_n-\left[\beta_{k0}+\beta_kx\right]_z}{\left[\Sigma_{k}\right]^{1/2}_{z,z}}\right)\right) \nn\\
	&= \sum_{k=1}^K{g}_{k}\left(x;\gamma\right)\sum_{z=1}^q \left(\bP\left(U > \frac{M_n-\left[\beta_{k0}+\beta_kx\right]_z}{\left[\Sigma_{k}\right]^{1/2}_{z,z}}\right) + \bP\left(U > \frac{M_n+\left[\beta_{k0}+\beta_kx\right]_z}{\left[\Sigma_{k}\right]^{1/2}_{z,z}}\right)\right) \nn\\
	&\le \sum_{k=1}^K {g}_{k}\left(x;\gamma\right)\sum_{z=1}^q \left[e^{-\frac{1}{2} \left(\frac{M_n-\left[\beta_{k0}+\beta_kx\right]_z}{\left[\Sigma_{k}\right]^{1/2}_{z,z}}\right)^2}+e^{-\frac{1}{2} \left(\frac{M_n+\left[\beta_{k0}+\beta_kx\right]_z}{\left[\Sigma_{k}\right]^{1/2}_{z,z}}\right)^2}\right] \left(\text{using Lemma \ref{lem.Chernoff.Inequality}, \eqref{eq.Chernoff-bound.Gaussian}}\right) \nn\\
	&\le 2 \sum_{k=1}^K {g}_{k}\left(x;\gamma\right)\sum_{z=1}^q e^{-\frac{1}{2} \left(\frac{M_n-\left|\left[\beta_{k0}+\beta_kx\right]_z\right|}{\left[\Sigma_{k}\right]^{1/2}_{z,z}}\right)^2} 
	\le 2 \sum_{k=1}^K {g}_{k}\left(x;\gamma\right)\sum_{z=1}^q e^{-\frac{1}{2} \frac{M^2_n-2 M_n\left|\left[\beta_{k0}+\beta_kx\right]_z\right|+\left|\left[\beta_{k0}+\beta_kx\right]\right|^2_z}{\left[\Sigma_{k}\right]_{z,z}}}\nn\\
	&\le 2 K A_\gamma q e^{-\frac{M_n^2-2M_n A_\beta}{2 A_\Sigma}},
\end{align}
where 
$$Y_{x,k} \sim \cN_q\left(\beta_{k0}+\beta_kx,\Sigma_{k}\right),\left[Y_{x,k}\right]_z \sim \cN\left(\left[\beta_{k0}+\beta_kx\right]_z,\left[\Sigma_{k}\right]_{z,z}\right),\text{ and } U = \frac{\left[Y_{x,k}\right]_z-\left[\beta x\right]_z}{\left[\Sigma_{k}\right]^{1/2}_{z,z}}\sim \cN(0,1),$$
and using the facts that $e^{-\frac{1}{2} \frac{\left|\left[\beta_{k0}+\beta_kx\right]\right|^2_z}{A_\Sigma}}
\le 1$ and $\max_{ 1 \le z \le q}\left|\left[\Sigma_{k}\right]_{z,z}\right| \le \left\|\Sigma_{k}\right\|_2 = M\left(\Sigma_{k}\right) = m\left(\Sigma^{-1}_{k}\right) \le A_\Sigma$.
We derive from \eqref{Extend.CarolineMeynet.eq.4.32} and \eqref{Extend.CarolineMeynet.eq.4.36} that
\begin{align}\label{Extend.CarolineMeynet.eq.4.37}
	\bP\left(\cT^c\right) \le 2 K n q A_\gamma e^{-\frac{M_n^2-2M_n A_\beta}{2 A_\Sigma}},
\end{align}
and finally from \eqref{eq_model_Sm_LL7},
\eqref{Extend.CarolineMeynet.eq.4.31},
and \eqref{Extend.CarolineMeynet.eq.4.37}, we obtain

\begin{align} \label{Extend.CarolineMeynet.eq.4.38}
	\E{}{{\kl}_n\left(s_0,\widehat{s}_{\widehat{m}}\right){\Indi}_{\cT^C}} \le  \left(\frac{ e^{q/2-1}\pi^{q/2} }{A_\Sigma ^{q/2} } + H_{s_0}\right)  \sqrt{2 K n q A_\gamma} e^{-\frac{M_n^2-2M_n A_\beta}{4 A_\Sigma}}.
\end{align}

\section{Proofs of technical lemmas\label{sec_proof_Lemma}}
\subsection{Proof of Lemma \ref{lem_control_deviation}} \label{sec_proof_lem_control_deviation}
Let $m \in \Ns$, on the event $\cT$, to control the deviation
\begin{align}\label{eq_define_deviation}
	\sup_{f_{m} \in F_{m}} \left|\nu_n\left(-f_{m}\right)\right|
	&= \sup_{f_{m} \in F_{m}} \left|\frac{1}{n} \sum_{i=1}^n \left\{f_m\left(Y_i| x_i\right)-\E{}{f_m\left(Y_i| x_i\right)}\right\}\right|,
\end{align}	
we shall use concentration and symmetrization arguments. We shall first use the following concentration inequality, which is an adaption of \citet[Theorem 4.10]{wainwright2019high}.

\begin{lemma}[Theorem 4.10 from \citet{wainwright2019high}]\label{Extend.CarolineMeynet.Lemma.5.1}
	Let $Z_1,\ldots,Z_n$ be independent random variables with values in some space $\cZ$ and let $\cF$ be a class of integrable real-valued functions with domain on $\cZ$. 
	Assume that 
	\begin{align}\label{eq_empirical_uniform_upperbound}
		\sup_{f \in \cF} \norm{f}_\infty  \le R_n \quad  \text{for some non-random constants $R_n < \infty$}.
	\end{align}	
	Then, for all $t>0$,
	\begin{align}\label{Extend.CarolineMeynet.eq.5.2}
		&\bP\left(\sup_{f \in \cF} \left|\frac{1}{n} \sum_{i=1}^n \left[f(Z_i) - \E{}{f\left(Z_i\right)}\right]\right|> \E{}{\sup_{f \in \cF} \left|\frac{1}{n} \sum_{i=1}^n \left[f(Z_i) - \E{}{f\left(Z_i\right)}\right]\right|}+2 \sqrt{2} R_n \sqrt{\frac{t}{n}}\right)
		\le e^{-t}.
	\end{align}
	That is, with probability greater than $1-e^t$,
	\begin{align}\label{Extend.CarolineMeynet.eq.5.2.equiv}
		\sup_{f \in \cF} \left|\frac{1}{n} \sum_{i=1}^n \left[f(Z_i) - \E{}{f\left(Z_i\right)}\right]\right|\le \E{}{\sup_{f \in \cF} \left|\frac{1}{n} \sum_{i=1}^n \left[f(Z_i) - \E{}{f\left(Z_i\right)}\right]\right|}+2 \sqrt{2} R_n \sqrt{\frac{t}{n}}
	\end{align}
\end{lemma}
Then, we propose to bound $\E{}{\sup_{f \in \cF} \left|\frac{1}{n} \sum_{i=1}^n \left[f(Z_i) - \E{}{f\left(Z_i\right)}\right]\right|}$ due to the following symmetrization argument. The proof of this result can be found in \citet{van1996weak}.
\begin{lemma}[See Lemma 2.3.6 in \citet{van1996weak}]\label{Extend.CarolineMeynet.Lemma.5.2} Let $Z_1,\ldots,Z_n$ be independent random variables with values in some space $\cZ$ and let $\cF$ be a class of real-valued functions on $\cZ$. Let $\left(\epsilon_1, \ldots, \epsilon_n\right)$ be a Rademacher sequence independent of $\left(Z_1,\ldots,Z_n\right)$. Then,
	\begin{align} \label{Extend.CarolineMeynet.eq.5.3}
		\E{}{\sup_{f \in \cF} \left|\frac{1}{n} \sum_{i=1}^n \left[f(Z_i) - \E{}{f\left(Z_i\right)}\right]\right|} \le 2 \E{}{\sup_{f \in \cF} \left|\frac{1}{n} \sum_{i=1}^n \epsilon_i f(Z_i)\right|}.
	\end{align}
\end{lemma}
From \eqref{Extend.CarolineMeynet.eq.5.3}, the problem is to provide an upper bound on $$\E{}{\sup_{f \in \cF} \left|\frac{1}{n} \sum_{i=1}^n \epsilon_i f(Z_i)\right|}.$$ To do so, we shall apply the following lemma, which is adapted from Lemma 6.1 in \citet{massart_concentration_2007}.
\begin{lemma}[See Lemma 6.1 in \citet{massart_concentration_2007}] \label{Extend.CarolineMeynet.Lemma.5.3}
	Let $Z_1,\ldots,Z_n$ be independent random variables with values in some space $\cZ$ and let $\cF$ be a class of real-valued functions on $\cZ$. Let $(\epsilon_1, \ldots, \epsilon_n)$ be a Rademacher sequence, independent of $(Z_1,\ldots,Z_n)$. Define $R_n$, a non-random constant, such that 
	\begin{align}\label{Extend.CarolineMeynet.eq.5.4}
		\sup_{f \in \cF} \norm{f}_n \le R_n.
	\end{align}
	Then, for all $S\in \Ns$, 
	\begin{align} \label{Extend.CarolineMeynet.eq.5.5}
		\E{}{\sup_{f \in \cF} \left|\frac{1}{n} \sum_{i=1}^n \epsilon_i f(Z_i)\right|} \le R_n\left(\frac{6}{\sqrt{n}} \sum_{s=1}^S 2^{-s} \sqrt{\ln \left[1+M\left(2^{-s}R_n,\cF,\norm{.}_n\right)\right]}+2^{-S}\right),
	\end{align}
	where $M\left(\delta,\cF,\norm{.}_n\right)$ stands for the $\delta$-packing number (see Definition \ref{def.delta-packing.Number}) of the set of functions $\cF$, equipped with the metric induced by the norm $\norm{\cdot}_n$.
\end{lemma}
We are now able to prove Lemma \ref{lem_control_deviation}. Indeed, given any fixed values $x_1,\ldots,x_n \in \cX$, in order to control $\sup_{f_m \in F_m} \left|\nu_n\left(-f_m\right)\right|{\Indi}_{\cT}$ from the \eqref{eq_define_deviation}, we would like to apply Lemmas \ref{Extend.CarolineMeynet.Lemma.5.1}--\ref{Extend.CarolineMeynet.Lemma.5.3}.
On the one hand, we see from \eqref{Extend.CarolineMeynet.eq.5.4} that we need an upper bound of $ \sup_{f_m \in F_m} \left\|f_m\right\|_\infty {\Indi}_{\cT}$. On the other hand, we see from \eqref{Extend.CarolineMeynet.eq.5.5} that on the event $\cT$, we need to bound the entropy of the set of functions $F_m$, equipped with the metric induced by the norm $\norm{\cdot}_n$. Such bounds are provided by the two following lemmas.

Recall that given $M_n>0$, we considered the event
\begin{align}
	\cT= \left\{\max_{i \in [n]}  \left\|Y_i\right\|_\infty=\max_{i \in [n]} \max_{z \in [q]} \left|\left[Y_i\right]_z\right| \le M_n \right\}, \label{eq_define_event_M}
\end{align}
let $B_n = \max\left(A_\Sigma,1+KA_G\right)\left(1+ q\sqrt{q}\left(M_n+ A_\beta\right)^2 A_\Sigma\right).$
\begin{lemma}\label{lem_control_fm_T}
	On the event $\cT$, for all $m\in \Ns$,
	\begin{align}
		\label{Extend.CarolineMeynet.eq.5.6.Gaussian}
		\sup_{f_m \in F_m} \left\|f_m\right\|_\infty  \le 2KB_n \left(A_\gamma+  q A_\beta + \frac{ q\sqrt{q}}{a_{\Sigma}}\right)=:R_n.	
	\end{align}
\end{lemma}
{\bf Proof of  Lemma \ref{lem_control_fm_T}.}
See Appendix \ref{Extend.CarolineMeynet.ProofLemma.5.4}.

%
%

\begin{lemma}\label{Extend.CarolineMeynet.Lemma.5.5}
	Let $\delta >0$ and $m \in \Ns$. 
	On the event $\cT$, 
	we have the following upper bound of the $\delta$-packing number of the set of functions $F_m$, equipped with the metric induced by the norm $\norm{\cdot}_n$:
	\begin{align*}
		M\left(\delta,F_m,\norm{\cdot}_n\right)
		&\le  \left(2p+1\right)^{\frac{72 B^2_n q^2K^2 m^2}{\delta^2}}\left(1+ \frac{18B_nKqA_\beta}{\delta}\right)^K\left(1+ \frac{18B_nKA_\gamma}{\delta}\right)^K\left(1+ \frac{18B_nKq \sqrt{q}}{a_{\Sigma}\delta}\right)^K.
	\end{align*}
\end{lemma}
{\bf Proof of  Lemma \ref{Extend.CarolineMeynet.Lemma.5.5}.}
See Appendix \ref{Extend.CarolineMeynet.ProofLemma.5.5}.

\begin{lemma}[Lemma 5.9 from \citet{Meynet:2013aa}]\label{Extend.CarolineMeynet.Lemma.5.9}
	Let $\delta > 0$ and $\left(x_{ij}\right)_{i \in [n];j=1,\ldots,p} \in \R^{np}$. There exists a family $\cB$ of $\left(2p+1\right)^{\norm{x}^2_{\max,n}/\delta^2}$ vectors in $\R^p$, such that for all $\beta \in \R^p$, with $\norm{\beta}_1 \le 1$, where $\norm{x}^2_{\max,n} =  \frac{1}{n} \sum_{i=1}^n \max_{j \in \left\{1,\ldots,p\right\}}x^2_{ij}$, there exists $\beta' \in \cB$, such that
	\begin{align*}
		\frac{1}{n}\sum_{i=1}^n \left(\sum_{j=1}^p \left(\beta_j - \beta'_j\right)x_{ij}\right)^2 \le \delta^2.
	\end{align*}
\end{lemma}
{\bf Proof of Lemma \ref{Extend.CarolineMeynet.Lemma.5.9}.} See in the proof of \citet[Lemma 5.9]{Meynet:2013aa}.

Via the upper bounds provided in Lemmas \ref{lem_control_fm_T} and \ref{Extend.CarolineMeynet.Lemma.5.5}, we can apply Lemma \ref{Extend.CarolineMeynet.Lemma.5.3} to get an upper bound of
\begin{align}
	\E{}{\sup_{f_m \in \cF_m} \left|\frac{1}{n} \sum_{i=1}^n \epsilon_i f_m(Y_i| x_i)\right|} \text{ on the event $\cT$.}
\end{align}

In order to provide such an upper bound, Lemmas \ref{Extend.CarolineMeynet.Lemma.5.1} and \ref{Extend.CarolineMeynet.Lemma.5.3} can be utilized via defining a suitable class of integrable real-valued functions as follows:
\begin{align}
	\cF := \left\{f:=f_m{\Indi}_{\left|f_m\right| \le R_n}: f_m \in F_m, \right\}, \quad Z_i := Y_i | x_i, \forall i \in [n].\label{eq_define_F}
\end{align}
Indeed, by definition, it holds that
\begin{align}
	\sup_{f \in \cF} \|f\|_n \le \sup_{f \in \cF} \|f\|_\infty = \sup_{f \in \cF} \sup_{z \in \cZ}\left|f(z)\right|
	=  \sup_{f_m \in \cF_m} \sup_{z \in \cZ}\left|f_m\left(z\right){\Indi}_{\left|f_m\left(z\right)\right| \le R_n}\right| \le R_n. 
\end{align}
Note that the last inequality is valid since if $\left|f_m\left(z\right)\right| \le R_n$ then $\left|f_m\left(z\right){\Indi}_{\left|f_m\left(z\right)\right| \le R_n}\right|  = \left|f_m\left(z\right)\right| \le R_n $. Otherwise, if $\left|f_m\left(z\right)\right| > R_n$, then $\left|f_m\left(z\right){\Indi}_{\left|f_m\left(z\right)\right| \le R_n}\right|  = \left|f_m\left(z\right) \times0\right| = 0 \le R_n $.
We thus obtain the following results.
\begin{lemma}\label{Extend.CarolineMeynet.Lemma.5.6}
	Let $m\in\Ns$, consider $\left(\epsilon_1,\ldots,\epsilon_n\right)$, a Rademacher sequence independent of $\left(Y_1,\ldots,Y_n\right)$. Then, on the event $\cT$, it holds that
	\begin{align}
		&\E{}{\sup_{f_m \in \cF_m} \left|\frac{1}{n} \sum_{i=1}^n \epsilon_i f_m(Y_i| x_i)\right|} \le  \frac{74KB_nq}{\sqrt{n}} \Delta_m, \text{ where}\label{Extend.CarolineMeynet.eq.5.7.Gaussian}\nn\\
		&	\Delta_m:=   m \sqrt{\ln(2p+1)}\ln n+2\sqrt{K} \left(A_\gamma+  q A_\beta + \frac{ q\sqrt{q}}{a_{\Sigma}}\right).
	\end{align}
	%
	\end{lemma}
	{\bf Proof of Lemma \ref{Extend.CarolineMeynet.Lemma.5.6}.} See Appendix \ref{Extend.CarolineMeynet.ProofLemma.5.6}.
	%
	
	%
	%

	We now return to the proof of the Lemma \ref{lem_control_deviation}.

	Finally, on the event $\cT$, using \eqref{eq_define_F} and Lemma \ref{lem_control_fm_T}, for all $m \in \Ns$ and $t>0$, with probability greater than $1-e^{-t}$, we obtain
	\begin{align}
		\sup_{f_{m} \in F_{m}} \left|\nu_n\left(-f_{m}\right)\right|
		&= \sup_{f \in \cF} \left|\frac{1}{n} \sum_{i=1}^n \left\{f\left(Z_i\right)-\E{}{f\left(Z_i\right)}\right\}\right| \\
		& \le \E{}{\sup_{f \in \cF} \left|\frac{1}{n} \sum_{i=1}^n \left[f(Z_i) - \E{}{f\left(Z_i\right)}\right]\right|}+2 \sqrt{2} R_n \sqrt{\frac{t}{n}} \left(\text{using Lemma \ref{Extend.CarolineMeynet.Lemma.5.1}} \right) \\
		& \le \E{}{\sup_{f \in \cF} \left|\frac{1}{n} \sum_{i=1}^n \left[f(Z_i) - \E{}{f\left(Z_i\right)}\right]\right|}+2 \sqrt{2} R_n \sqrt{\frac{t}{n}} \left( \text{since $ \le 1$}\right) \\
		& \le 2 \E{}{\sup_{f \in \cF} \left|\frac{1}{n} \sum_{i=1}^n \epsilon_i f(Z_i)\right|}+2 \sqrt{2} R_n \sqrt{\frac{t}{n}}\left(\text{Lemma \ref{Extend.CarolineMeynet.Lemma.5.2}} \right) \\
		& = 2\E{}{\sup_{f_m \in \cF_m} \left|\frac{1}{n} \sum_{i=1}^n \epsilon_i f_m(Y_i| x_i)\right|}+2 \sqrt{2} R_n \sqrt{\frac{t}{n}} \\
		&\le  \frac{148KB_n q}{\sqrt{n}} \Delta_m+ 4 \sqrt{2} KB_n \left(A_\gamma+  q A_\beta + \frac{ q\sqrt{q}}{a_{\Sigma}}\right) \sqrt{\frac{t}{n}} \\
		&\left(\text{using Lemma \ref{Extend.CarolineMeynet.Lemma.5.6} and $R_n = 2KB_n \left(A_\gamma+  q A_\beta + \frac{ q\sqrt{q}}{a_{\Sigma}}\right)$}\right) \\
		&\le \frac{4KB_n}{\sqrt{n}} \left[37 q \Delta_m + \sqrt{2}\left(A_\gamma+  q A_\beta + \frac{ q\sqrt{q}}{a_{\Sigma}}\right)\sqrt{t}\right].
	\end{align}

\subsection{Proofs of Lemmas \ref{lem_control_fm_T}--\ref{Extend.CarolineMeynet.Lemma.5.6}}
The proofs of Lemmas \ref{lem_control_fm_T}--\ref{Extend.CarolineMeynet.Lemma.5.5} require an upper bound on the uniform norm of the gradient of $\ln s_\psi$, for $s_\psi \in S$. We begin by providing such an upper bound.
\begin{lemma}\label{Upperbound.GradientLog.Lemma}
	Given $s_\psi$, as described in \eqref{eq_Sm_bounded}, it holds that
	\begin{align}
		\sup_{x\in \cX} \sup_{\psi \in \widetilde{\Psi}} \left\|\frac{\partial \ln\left(s_\psi(\cdot| x)\right)}{\partial \psi}\right\|_{\infty} &\le G(\cdot),\nn\\
		G: \R^q \ni y \mapsto G(y)&= \max\left(A_\Sigma,1+KA_G\right)\left(1+ q\sqrt{q}\left(\norm{y}_\infty+ A_\beta\right)^2 A_\Sigma\right).\label{Upperbound.GradientLog}
	\end{align}
\end{lemma}

{\bf Proof of Lemma \ref{Upperbound.GradientLog.Lemma}.}
Let $s_\psi \in S$, with 
$\psi = \left(\gamma,\beta,\Sigma\right)$. From now on, we consider any $x \in \cX$,  any $y\in \R^q$, and any $k \in [K]$. We can write
\begin{align*}
	\ln \left(s_\psi(y | x)\right) &= \ln \left(\sum_{k=1}^K {g}_{k}\left(x;\gamma\right)\cN\left(y;\beta_{k0}+\beta_kx,\Sigma_k\right)\right) = \ln \left(\sum_{k=1}^K f_k(x,y)\right),\\
	{g}_{k}\left(x;\gamma\right)&=  \frac{\exp\left(w_k(x)\right)}{\sum_{l=1}^K\exp\left(w_l(x)\right)}, w_k(x) = \gamma_{k0} + \gamma_k^\top x,\\
	\cN\left(y;\beta_{k0}+\beta_kx,\Sigma_k\right) &= \frac{1}{\left(2 \pi \right)^{q/2} \det(\Sigma_k)^{1/2}} \exp\left(-\frac{\left(y-\left(\beta_{k0}+\beta_kx\right)\right)^\top  \Sigma_k^{-1} \left(y-\left(\beta_{k0}+\beta_kx\right)\right)}{2}\right),\\
	f_k (x,y) &= {g}_{k}\left(x;\gamma\right)\cN\left(y;\beta_{k0}+\beta_kx,\Sigma_k\right)\nn\\
	&=\frac{{g}_{k}\left(x;\gamma\right)}{\left(2 \pi \right)^{q/2} \det(\Sigma_k)^{1/2}} \exp\left[-\frac{1}{2}\left(y-\left(\beta_{k0}+\beta_kx\right)\right)^\top  \Sigma_k^{-1} \left(y-\left(\beta_{k0}+\beta_kx\right)\right)\right].
\end{align*}
By using the chain rule, for all $l \in [K]$,
\begin{align*}
	\frac{\partial \ln \left(s_\psi(y | x)\right)}{\partial \gamma_{l0}}
	&= \sum_{k=1}^K\frac{f_k (x,y)}{{g}_{k}\left(x;\gamma\right)\sum_{k=1}^K f_k (x,y)}\frac{\partial {g}_{k}\left(x;\gamma\right) }{  \partial w_l(x)}\underbrace{\frac{\partial  w_l(x) }{  \partial \gamma_{l0}}}_{=1},\text{and}\\
	\frac{\partial \ln \left(s_\psi(y | x)\right)}{\partial \left(\gamma^\top _{l}x\right)}
	&= \sum_{k=1}^K\frac{f_k (x,y)}{{g}_{k}\left(x;\gamma\right)\sum_{k=1}^K f_k (x,y)}\frac{\partial {g}_{k}\left(x;\gamma\right) }{  \partial w_l(x)}\underbrace{\frac{\partial  w_l(x) }{  \partial \left(\gamma^\top _{l}x\right)}}_{=1}.
\end{align*}
Furthermore,
\begin{align*} 
	\frac{\partial {g}_{k}\left(x;\gamma\right) }{  \partial w_l(x)} &= \frac{\partial }{  \partial w_l(x)}\left(\frac{\exp\left(w_k(x)\right)}{\sum_{l=1}^K\exp\left(w_l(x)\right)}\right)\nn\\
	&= \frac{\delta_{lk} \exp\left(w_k(x)\right)}{\sum_{l=1}^K\exp\left(w_l(x)\right)}-\frac{\exp\left(w_k(x)\right)}{\sum_{l=1}^K\exp\left(w_l(x)\right)}\frac{\exp\left(w_l(x)\right)}{\sum_{l=1}^K\exp\left(w_l(x)\right)}= {g}_{k}\left(x;\gamma\right)\left(\delta_{lk} - {g}_{l}\left(x;\gamma\right)\right), \nn\\
	& \quad 
	\text{ where }\delta_{lk}  = 
	\begin{cases}
		1& \text{ if } l=k,\\
		0& \text{ if } l \neq k.
	\end{cases}
\end{align*}
Therefore, we obtain 
\begin{align*}
	\left|\frac{\partial \ln \left(s_\psi(y | x)\right)}{\partial \left(\gamma_l^\top  x\right)}\right|
	&=\left|\frac{\partial \ln \left(s_\psi(y | x)\right)}{\partial \gamma_{l0}}\right|=\left| \sum_{k=1}^K\frac{f_k (x,y)}{{g}_{k}\left(x;\gamma\right)\sum_{k=1}^K f_k (x,y)}{g}_{k}\left(x;\gamma\right)\left(\delta_{lk} - {g}_{l}\left(x;\gamma\right)\right)\right|\nn\\
	&= \left| \sum_{k=1}^K\frac{f_k (x,y)}{\sum_{k=1}^K f_k (x,y)}\left(\delta_{lk} - {g}_{l}\left(x;\gamma\right)\right)\right|\le\left| \sum_{k=1}^K\left(\delta_{lk} - {g}_{l}\left(x;\gamma\right)\right)\right| 
	\nn\\
	&= \left|1- \sum_{k=1}^K{g}_{l}\left(x;\gamma\right)\right|
	= \left|1- K{g}_{l}\left(x;\gamma\right)\right| \le 1+K{g}_{l}\left(x;\gamma\right) \le 1 + K A_G \left( \text{using \eqref{eq.defBoundedGating}}\right).
\end{align*}
Similarly, by using the fact that  $\psi $ belongs to the bounded space $\widetilde{\Psi}$, $f_l(x,y)/\sum_{k=1}^K f_k(x,y) \le 1$,
\begin{align*}
	&\left\|\frac{\partial \ln \left(s_\psi(y | x)\right)}{\partial \beta_{l0}}\right\|_\infty
	=
	\left\|\frac{\partial \ln \left(s_\psi(y | x)\right)}{\partial \left(\beta_l x\right)}\right\|_\infty\nn\\
	&= \left\| \frac{f_l(x,y)}{\sum_{k=1}^K f_k(x,y)}\frac{\partial}{\partial\left(\beta_{l0}+\beta_lx\right)} \left[-\frac{1}{2}\left(y-\left(\beta_{l0}+\beta_lx\right)\right)^\top  \Sigma_l^{-1} \left(y-\left(\beta_{l0}+\beta_lx\right)\right)\right]\right\|_\infty\nn\\
	&\le \left\| \frac{\partial}{\partial\left(\beta_{l0}+\beta_lx\right)} \left[-\frac{1}{2}\left(y-\left(\beta_{l0}+\beta_lx\right)\right)^\top  \Sigma_l^{-1} \left(y-\left(\beta_{l0}+\beta_lx\right)\right)\right]\right\|_\infty\nn\\
	&=\left\|\Sigma_l^{-1}\left(y-\left(\beta_{l0}+\beta_lx\right)\right)\right\|_\infty \le \left\|\Sigma_l^{-1}\right\|_\infty\left\|\left(y-\left(\beta_{l0}+\beta_lx\right)\right)\right\|_\infty \left(\text{using \eqref{norm.Equality}}\right)\nn\\
	&\le \sqrt{q}\left\|\Sigma_l^{-1}\right\|_2\left(\left\|y\right\|_\infty + \left\|\beta_{l0}+\beta_lx\right\|_\infty\right) \left(\text{using \eqref{eq.Ineq.Infty-2-norm}}\right)\nn\\
	&\le \sqrt{q}M\left(\Sigma_l^{-1}\right)\left(\left\|y\right\|_\infty + \left\|\beta_{l0}+\beta_lx\right\|_\infty\right) \left(\text{using \eqref{eq.Ineq.2-norm_Vect-Infty}}\right)\nn\\
	&\le  \sqrt{q}A_\Sigma\left(\left\|y\right\|_\infty + A_\beta\right) \left(\text{using \eqref{eq_defBoundedParameters}}\right).
\end{align*}

Now, we need to calculate the gradient \wrt~to the covariance matrices of the Gaussian experts. To do this, we need the following result: given any $l \in [K]$, $v_l = \beta_{l0} + \beta_l x$, it holds that
\begin{align}
	\frac{\partial }{ \partial \Sigma_l}\cN\left(x;v_l,\Sigma_l\right)
	&=\cN\left(x;v_l,\Sigma_l\right) \underbrace{\frac{1}{2}\left[ \Sigma^{-1}_l\left(x-v_l\right)  \left(x-v_l\right)^\top \Sigma^{-1}_l -  \left(\Sigma^{-1}_l\right)^\top \right]}_{T\left(x,v_l,\Sigma_l\right)},\label{eq.derivative.CovarianceGaussian}
\end{align}
noting that
\begin{align}
	\frac{\partial }{ \partial \Sigma_l}\left(\left(x-v_l\right)^\top  \Sigma^{-1}_l \left(x-v_l\right)\right) &=- \Sigma^{-1}_l\left(x-v_l\right)  \left(x-v_l\right)^\top \Sigma^{-1}_l \left(\text{using Lemma \ref{derivative.QuadraticPrecision}}\right),\\
	\frac{\partial }{ \partial \Sigma_l} \left(\det(\Sigma_l)\right) &= \det(\Sigma_l) \left(\Sigma^{-1}_l\right)^\top  \left(\text{using Jacobi formula, Lemma \ref{Jacobi's formula}}\right). \label{eq.derivative.determinant}
\end{align}
For any $l \in [K]$,
\begin{align*}
	\left|\frac{\partial \ln \left(s_\psi(y | x)\right)}{\partial \left(\left[\Sigma_l\right]_{z_1,z_2}\right)}\right| &\le
	\left\|	\frac{\partial \ln \left(s_\psi(y | x)\right)}{\partial \Sigma_l}\right\|_2 \left(\text{using \eqref{eq.Ineq.2-norm_Vect-Infty}}\right)\nn\\
	&= \left|\frac{f_l(x,y)}{\sum_{k=1}^K f_k(x,y)}\right|\left\| \frac{\partial}{\partial\Sigma_l} \left[-\frac{1}{2}\left(y-\left(\beta_{l0}+\beta_lx\right)\right)^\top  \Sigma_l^{-1} \left(y-\left(\beta_{l0}+\beta_lx\right)\right)\right]\right\|_2\nn\\
	&\le \left\| \frac{\partial}{\partial\Sigma_l} \left[-\frac{1}{2}\left(y-\left(\beta_{l0}+\beta_lx\right)\right)^\top  \Sigma_l^{-1} \left(y-\left(\beta_{l0}+\beta_lx\right)\right)\right]\right\|_2\nn\\
	&= \frac{1}{2} \left\|\Sigma_l^{-1}\left(y-\left(\beta_{l0}+\beta_lx\right)\right)  \left(y-\left(\beta_{l0}+\beta_lx\right)\right)^\top \Sigma_l^{-1} -  \left(\Sigma^{-1}_l\right)^\top \right\|_2\ \left(\text{using \eqref{eq.derivative.CovarianceGaussian}}\right)\nn\\
	&\le \frac{1}{2} \left[A_\Sigma+ \sqrt{q}\left\|\left(y-\left(\beta_{l0}+\beta_lx\right)\right)  \left(y-\left(\beta_{l0}+\beta_lx\right)\right)^\top  \right\|_\infty A^2_{\Sigma}\right]\left( \text{using \eqref{eq.Ineq.Infty-2-norm}}\right)\nn\\
	&\le \frac{1}{2}\left[A_\Sigma+q\sqrt{q}\left(\norm{y}_\infty + A_{\beta}\right)^2 A^2_{\Sigma}\right] \left(\text{using \eqref{eq_defBoundedParameters}}\right),
\end{align*}
where, in the last inequality given $a = y - \left(\beta_{l0}+\beta_lx\right)$, we use the fact that
\begin{align*}
	\left\|aa^\top \right\|_{\infty} = \max_{1 \le i \le q} \sum_{j=1}^q \left|[aa^\top ]_{i,j}\right| = \max_{1 \le i \le q} \sum_{j=1}^q \left|a_{i}a_{j}\right| = \max_{1 \le i \le q}\left|a_{i}\right| \sum_{j=1}^q \left|a_{j}\right|\le q \left\| a\right\|^2_{\infty}.
\end{align*}
Thus,
\begin{align*}
	&\sup_{x\in \cX} \sup_{\psi \in \widetilde{\Psi}} \left\|\frac{\partial \ln\left(s_\psi(y| x)\right)}{\partial \psi}\right\|_{\infty}\nn\\
	&\le \max\Bigg[1 + K A_G,\sqrt{q}\left(\left\|y\right\|_\infty  + A_\beta\right)A_\Sigma, \frac{1}{2}\left[A_\Sigma+q\sqrt{q}\left(\norm{y}_\infty + A_{\beta}\right)^2 A^2_{\Sigma}\right] \Bigg]\nn\\
	&\le \max\Bigg[1+KA_G,\max\left(A_\Sigma,1\right)\left(1+ q\sqrt{q}\left(\norm{y}_\infty+ A_\beta\right)^2 A_\Sigma\right)\Bigg]\nn\\
	&\le  \max\left(A_\Sigma,1+KA_G\right)\left(1+ q\sqrt{q}\left(\norm{y}_\infty+ A_\beta\right)^2 A_\Sigma\right)
	=: G(y),
\end{align*}
where we use the fact that
\begin{align*}
	\sqrt{q}\left(\left\|y\right\|_\infty  + A_\beta\right)A_\Sigma=:\theta &\le 1+ \theta^2=1+ q\left(\norm{y}_\infty+ A_\beta\right)^2 A^2_\Sigma \nn\\
	&\le\max\left(A_\Sigma,1\right)\left(1+ q\sqrt{q}\left(\norm{y}_\infty+ A_\beta\right)^2 A_\Sigma\right).
\end{align*}

\subsubsection{Proof of Lemma \ref{lem_control_fm_T}} \label{Extend.CarolineMeynet.ProofLemma.5.4}
Let $m\in \Ns$ and $f_m \in F_m$. By \eqref{eq_Fm}, there exists $s_m \in S_m$, such that $f_m = -\ln \left(s_m/s_0\right)$. For all $x\in \cX$, let $\psi(x) = \left(\gamma_{k0},\gamma_k^\top x,\beta_{k0},\beta_kx,\Sigma_k\right)_{k\in [K]}$ be the parameters of $s_m\left(
\cdot| x\right)$. 
In our case, we approximate $f(\psi) = \ln\left(s_\psi\left(y_i | x_i\right)\right)$ around $\psi_0(x_i)$ by the $n=0^{th}$ degree Taylor polynomial of $f(\psi)$. That is,

\begin{align*}
	\left|\ln\left(s_m\left(y_i | x_i\right)\right)-\ln\left(s_0\left(y_i | x_i\right)\right)\right| &=:\left|f(\psi) - f(\psi_0) \right| = \left|R_0(\psi)\right| \left(\text{defined in Lemma \ref{Taylor'sInequality}}\right)\nn\\
	&\le \sup_{x\in \cX} \sup_{\psi \in \widetilde{\Psi}} \left\|\frac{\partial \ln\left(s_\psi(y_i| x)\right)}{\partial \psi}\right\|_\infty \left\|\psi\left(x_i\right)-\psi_{0}\left(x_i\right)\right\|_1.
\end{align*}

First applying Taylor's inequality and then Lemma \ref{Upperbound.GradientLog.Lemma} on the event $\cT$. For all $i \in [n]$, it holds that
\begin{align*}
	&\left|f_m\left(y_i | x_i\right)\right|{\Indi}_{\cT}
	=\left|\ln\left(s_m\left(y_i | x_i\right)\right)-\ln\left(s_0\left(y_i | x_i\right)\right)\right|{\Indi}_{\cT} \le \sup_{x\in \cX} \sup_{\psi \in \widetilde{\Psi}} \left\|\frac{\partial \ln\left(s_\psi(y_i| x)\right)}{\partial \psi}\right\|_\infty \left\|\psi\left(x_i\right)-\psi_{0}\left(x_i\right)\right\|_1{\Indi}_{\cT}\nn\\
	&\le \underbrace{\max\left(A_\Sigma,1+KA_G\right)\left(1+ q\sqrt{q}\left(M_n+ A_\beta\right)^2 A_\Sigma\right)}_{=:B_n}\left\|\psi\left(x_i\right)-\psi_{0}\left(x_i\right)\right\|_1 \left(\text{using Lemma \ref{Upperbound.GradientLog.Lemma}}\right)\nn\\
	&\le B_n \sum_{k=1}^K \left(\left|\gamma_{k0}  - \gamma_{0,k0}\right|+\left|\gamma_k^\top  x_i - \gamma_{0,k}^\top  x_i\right| +\left\|\beta_{k0} - \beta_{0,k0}\right\|_1+\left\|\beta_k x_i - \beta_{0,k} x_i\right\|_1+\left\|\vect\left(\Sigma_k - \Sigma_{0,k}\right)\right\|_1\right)\nn\\
	&\le 2 B_n \sum_{k=1}^K \left(\left|\gamma_{k0}\right|+\left|\gamma_k^\top  x_i\right|+\left\|\beta_{k0} \right\|_1+\left\|\beta_k x_i\right\|_1+q\left\|\Sigma_k \right\|_1\right)\left(\text{using \eqref{eq.Ineq.1-norm_Vect-Infty}}\right)\nn\\
	&\le 2KB_n \left( A_\gamma+q\left\|\beta_{k0} \right\|_\infty+q\left\|\beta_k x_i\right\|_\infty+q\sqrt{q}\left\|\Sigma_k \right\|_2\right) \left(\text{using \eqref{eq_defBoundedParameters}, \eqref{eq.vector.infinity1.equiv}, \eqref{eq.vector.12.equiv}, \eqref{eq.Ineq.1-2-norm}}\right) \nn\\
	&\le 2KB_n \left(A_\gamma+  q A_\beta + \frac{ q\sqrt{q}}{a_{\Sigma}}\right) \left(\text{using \eqref{eq_defBoundedParameters}}\right).
\end{align*}
Therefore, on the event $\cT$,
\begin{align*}
	\sup_{f_m \in F_m} \left\|f_m\right\|_\infty  \le 2KB_n \left(A_\gamma+  q A_\beta + \frac{ q\sqrt{q}}{a_{\Sigma}}\right)=:R_n.
\end{align*}

\subsubsection{Proof of Lemma \ref{Extend.CarolineMeynet.Lemma.5.5}}\label{Extend.CarolineMeynet.ProofLemma.5.5}
Let $m\in \Ns$, $f^{[1]}_m \in F_m$, and  $x\in [0,1]^p$. By \eqref{eq_Fm}, there exists $s^{[1]}_m \in S_m$, such that $f^{[1]}_m = -\ln \left(s^{[1]}_m/s_0\right)$. Introduce the notation $s^{[2]}_m \in S$ and  $f^{[2]}_m = -\ln \left(s^{[2]}_m/s_0\right)$. Let
\begin{align*}
	\psi^{[1]}(x) = \left(\gamma^{[1]}_{k0},\gamma^{[1]}_kx,\beta^{[1]}_{k0},\beta^{[1]}_kx,\Sigma^{[1]}_k\right)_{k\in [K]}, \text{and } \psi^{[2]}(x) = \left(\gamma^{[2]}_{k0},\gamma^{[2]}_kx,\beta^{[2]}_{k0},\beta^{[2]}_kx,\Sigma^{[2]}_k\right)_{k\in [K]},
\end{align*}
be the parameters of the PDFs $s^{[1]}_m\left(\cdot| x\right)$ and $s^{[2]}_m\left(\cdot| x\right)$, respectively. By applying Taylor's inequality and then Lemma \ref{Upperbound.GradientLog.Lemma} on the event $\cT$, for all $i \in [n]$, it holds that
\begin{align*}
	&\left|f^{[1]}_m\left(y_i | x_i\right)-f^{[2]}_{m}\left(y_i | x_i\right)\right|= \left|\ln\left(s^{[1]}_m\left(y_i | x_i\right)\right)-\ln\left(s^{[2]}_m\left(y_i | x_i\right)\right)\right|\nn\\
	&\le \sup_{x\in \cX} \sup_{\psi \in \widetilde{\Psi}} \left|\frac{\partial \ln\left(s_\psi(y_i| x)\right)}{\partial \psi}\right|\left\|\psi^{[1]}\left(x_i\right)-\psi^{[1]}\left(x_i\right)\right\|_1 \left(\text{using Taylor's inequality in Lemma \ref{Taylor'sInequality}}\right)\nn\\
	&\le \underbrace{\max\left(A_\Sigma,C(p,K)\right)\left(1+ q\sqrt{q}\left(M_n+ A_\beta\right)^2 A_\Sigma\right)}_{B_n}\left\|\psi^{[1]}\left(x_i\right)-\psi^{[2]}\left(x_i\right)\right\|_1 \left(\text{using Lemma \ref{Upperbound.GradientLog.Lemma}}\right)\nn\\
	&\le B_n \sum_{k=1}^K \Bigg(\left|\gamma^{[1]}_{k0}  - \gamma^{[2]}_{k0}\right|+\left|\gamma^{[1]^\top }_k x_i - \gamma^{[2]^\top }_k x_i\right| \nn\\
	&\hspace{4cm} +\left\|\beta_{k0}^{[1]} - \beta^{[2]}_{k0}\right\|_1+\left\|\beta^{[1]}_k x_i - \beta^{[2]}_{k} x_i\right\|_1+\left\|\vect\left(\Sigma^{[1]}_k - \Sigma^{[2]}_{k}\right)\right\|_1\Bigg)\nn.
\end{align*}
By the Cauchy-Schwarz inequality, $\left(\sum_{i=1}^m a_i\right)^2 \le m \sum_{i=1}^m a^2_i$ ($m \in\Ns$), we get
\begin{align*}
	&\left|f^{[1]}_m\left(y_i | x_i\right)-f^{[2]}_{m}\left(y_i | x_i\right)\right|^2\nn\\
	&\le 3B^2_n \left[\left(\sum_{k=1}^K\left|\gamma^{[1]^\top }_k x_i - \gamma^{[2]^\top }_k x_i\right|\right)^2+\left(\sum_{k=1}^K\sum_{z=1}^q\left|\left[\beta^{[1]}_k x_i\right]_z - \left[\beta^{[2]}_k x_i\right]_z\right|\right)^2 \right]\nn\\
	&\quad + 3B^2_n \left(\left\|\beta^{[1]}_{0}  - \beta^{[2]}_{0}\right\|_1+\left\|\gamma^{[1]}_{0}  - \gamma^{[2]}_{0}\right\|_1+\left\|\vect\left(\Sigma^{[1]} - \Sigma^{[2]}\right)\right\|_1\right)^2\nn\\
	&\le 3B^2_n \left[K\sum_{k=1}^K\left(\sum_{j=1}^p\gamma^{[1]^\top }_{kj} x_{ij} - \sum_{j=1}^p\gamma^{[2]^\top }_{kj} x_{ij}\right)^2 Kq\sum_{k=1}^K\sum_{z=1}^q\left(\sum_{j=1}^p\left[\beta^{[1]}_k\right]_{z,j}x_{ij} - \sum_{j=1}^p\left[\beta^{[2]}_k\right]_{z,j}x_{ij}\right)^2 \right] \nn\\
	&\quad + 3B^2_n \left(\left\|\beta^{[1]}_{0}  - \beta^{[2]}_{0}\right\|_1+\left\|\gamma^{[1]}_{0}  - \gamma^{[2]}_{0}\right\|_1+\left\|\vect\left(\Sigma^{[1]} - \Sigma^{[2]}\right)\right\|_1\right)^2,
\end{align*}
and
\begin{align*}
	\norm{f^{[1]}_m-f^{[2]}_m}^2_n
	&=\frac{1}{n}\sum_{i=1}^n\left|f^{[1]}_m\left(y_i | x_i\right)-f^{[2]}_{m}\left(y_i | x_i\right)\right|^2\nn\\ 
	&\le 3B^2_n\underbrace{K\sum_{k=1}^K\frac{1}{n}\sum_{i=1}^n\left(\sum_{j=1}^p\gamma^{[1]}_{kj} x_{ij} - \sum_{j=1}^p\gamma^{[2]}_{kj} x_{ij}\right)^2}_{=:a}\nn\\
	&\hspace{0.1cm}+3B^2_n\underbrace{Kq\sum_{k=1}^K\sum_{z=1}^q\frac{1}{n}\sum_{i=1}^n\left(\sum_{j=1}^p\left[\beta^{[1]}_k\right]_{z,j}x_{ij} - \sum_{j=1}^p\left[\beta^{[2]}_k\right]_{z,j}x_{ij}\right)^2}_{=:b}\nn\\
	&\hspace{0.1cm}+3B^2_n\left(\left\|\beta^{[1]}_{0}  - \beta^{[2]}_{0}\right\|_1+\left\|\gamma^{[1]}_{0}  - \gamma^{[2]}_{0}\right\|_1+\left\|\vect\left(\Sigma^{[1]} - \Sigma^{[2]}\right)\right\|_1\right)^2.
\end{align*}
So, for all $\delta > 0$, if 
\begin{align*}
	&a \le \delta^2/\left(36 B^2_n\right),
	b \le \delta^2/\left(36 B^2_n\right),
	\left\|\beta^{[1]}_{0}  - \beta^{[2]}_{0}\right\|_1 \le \delta/\left(18 B_n\right), \nn\\
	&\left\|\gamma^{[1]}_{0}  - \gamma^{[2]}_{0}\right\|_1 \le \delta/\left(18 B_n\right),\text{and}
	\left\|\vect\left(\Sigma^{[1]} - \Sigma^{[2]}\right)\right\|_1 \le \delta/\left(18 B_n\right),
\end{align*}
then $\norm{f^{[1]}_m-f^{[2]}_m}^2_n \le \delta^2/4$.
To bound $a$ and $b$, we can write
\begin{align*}
	a &= Km^2\sum_{k=1}^K\frac{1}{n}\sum_{i=1}^n\left(\sum_{j=1}^p\frac{\gamma^{[1]}_{kj}}{m} x_{ij} - \sum_{j=1}^p\frac{\gamma^{[2]}_{kj}}{m} x_{ij}\right)^2,\text{and } \nn\\
	b &= Kqm^2\sum_{k=1}^K\sum_{z=1}^q\frac{1}{n}\sum_{i=1}^n\left(\sum_{j=1}^p\frac{\left[\beta^{[1]}_k\right]_{z,j}}{m}x_{ij} - \sum_{j=1}^p\frac{\left[\beta^{[2]}_k\right]_{z,j}}{m}x_{ij}\right)^2.
\end{align*}
Then, we apply Lemma \ref{Extend.CarolineMeynet.Lemma.5.9} to obtain $\frac{\gamma^{[1]}_{k,.}}{m} =\left(\frac{\gamma^{[1]}_{kj}}{m}\right)_{j\in[q]}$ and  $\frac{\left[\beta^{[1]}_k\right]_{z,.}}{m} = \left(\frac{\left[\beta^{[1]}_k\right]_{z,j}}{m}\right)_{j\in[q]}$, for all $k\in[K],z\in [q]$. Since $s^{[1]}_m \in S_m$, and using \eqref{eq_model_Sm}, we have  $\left\|\gamma^{[1]}_k\right\| \le m$ and $\left\|\vect\left(\beta^{[1]}_k\right)\right\|_1 \le m$, which leads to $\sum_{j=1}^p\left|\frac{\gamma^{[1]}_{kj}}{m}\right|\le 1$ and $\sum_{z=1}^q\sum_{j=1}^p \left|\frac{\left[\beta^{[1]}_k\right]_{z,j}}{m}\right| \le 1$, respectively. Furthermore, given $x \in \cX	=[0,1]^p$, we have $\norm{x}^2_{\max,n} = 1$. Thus, there exist families $\cA$ of $\left(2p+1\right)^{36 B^2_n K^2 m^2/\delta^2}$ vectors and $\cB$ of $\left(2p+1\right)^{16 B^2_n q^2K^2 m^2/\delta^2}$ vectors of $\R^p$, such that for all $k\in[K]$, $z \in [q]$, $\gamma^{[1]}_{k,.}$, and$\left[\beta^{[1]}_k\right]_{z,.}$, there exist $\gamma^{[1]}_{k,.} \in \cA$ and $\left[\beta^{[2]}_k\right]_{z,.}\in \cB$, such that
\begin{align*}
	&\frac{1}{n}\sum_{i=1}^n\left(\sum_{j=1}^p\frac{\gamma^{[1]}_{kj}}{m} x_{ij} - \sum_{j=1}^p\frac{\gamma^{[2]}_{kj}}{m} x_{ij}\right)^2 \le \frac{\delta^2}{36 B^2_n K^2 m^2}, \text{and }\nn\\
	&\frac{1}{n}\sum_{i=1}^n\left(\sum_{j=1}^p\frac{\left[\beta^{[1]}_k\right]_{z,j}}{m}x_{ij} - \sum_{j=1}^p\frac{\left[\beta^{[2]}_k\right]_{z,j}}{m}x_{ij}\right)^2 \le \frac{\delta^2}{36 B^2_n q^2 K^2 m^2},
\end{align*}
which leads to $a \le \delta^2/36B^2_n$ and $b \le \delta^2/36B^2_n$.
Moreover, \eqref{eq_defBoundedParameters} leads to
\begin{align*}
	\left\|\beta^{[1]}_0\right\|_1 &= \sum_{k=1}^K \left\|\beta^{[1]}_{0k}\right\|_1 \le Kq \left\|\beta^{[1]}_{0k}\right\|_\infty \le Kq A_\beta \left(\text{using \eqref{eq.vector.infinity1.equiv}}\right), \nn\\
	\left\|\gamma^{[1]}_0\right\|_1 &= \sum_{k=1}^K \left|\gamma^{[1]}_{0k}\right| \le K A_\gamma,
	\left\|\vect\left(\Sigma^{[1]} \right)\right\|_1 
	\le \frac{K q\sqrt{q}}{a_{\Sigma}}.
\end{align*}
Therefore, on the event $\cT$,
\begin{align*}
	M\left(\delta,F_m,\norm{\cdot}_n\right)
	&\le N\left(\delta/2,F_m,\norm{\cdot}_n\right)\left( \text{using Lemma \ref{lem.delta-packingCovering.Lemma}}\right)\nn\\
	&\le\card(\cA)\card(\cB)N\left(\frac{\delta}{18B_n},B_1^K\left(KqA_\beta\right),\norm{\cdot}_1\right)\nn\\
	& \quad \times N\left(\frac{\delta}{18B_n},B_1^K\left(KA_\gamma\right),\norm{\cdot}_1\right) N\left(\frac{\delta}{18B_n},B_1^K\left(\frac{K q\sqrt{q}}{a_{\Sigma}}\right),\norm{\cdot}_1\right)\nn\\
	&\le \left(2p+1\right)^{\frac{72 B^2_n q^2K^2 m^2}{\delta^2}}\left(1+ \frac{18B_nKqA_\beta}{\delta}\right)^K\left(1+ \frac{18B_nKA_\gamma}{\delta}\right)^K\left(1+ \frac{18B_nKq \sqrt{q}}{a_{\Sigma}\delta}\right)^K.
\end{align*}
\subsubsection{Proof of Lemma \ref{Extend.CarolineMeynet.Lemma.5.6}}\label{Extend.CarolineMeynet.ProofLemma.5.6}
Let $m\in\Ns$. From Lemma \ref{lem_control_fm_T}, on the event $\cT$, it holds that
\begin{align}
	\sup_{f_m \in F_m} \left\|f_m\right\|_n  \le 2KB_n \left(A_\gamma+  q A_\beta + \frac{ q\sqrt{q}}{a_{\Sigma}}\right)=:R_n.\label{eq_define_deviation8.Gaussian}
\end{align}
From Lemma \ref{Extend.CarolineMeynet.Lemma.5.5}, on the event $\cT$ for all $S \in \Ns$, $\text{ with } \delta = 2^{-s}R_n$,
\begin{align}\label{eq_define_deviation9.Gaussian}
	&\sum_{s=1}^S 2^{-s}\sqrt{\ln \left[1+ M\left(2^{-s}R_n,F_m,\norm{\cdot}_n\right)\right]}\le{} \sum_{s=1}^S 2^{-s}\sqrt{\ln \left[ 2M\left(\delta,F_m,\norm{\cdot}_n\right)\right]}\nn\\
	\le{}& \sum_{s=1}^S 2^{-s} \Bigg[\sqrt{\ln 2}+ \frac{6\sqrt{2} B_n q K m }{\delta} \sqrt{\ln \left(2p+1\right)} \nn\\
	& \quad +\sqrt{K \ln \left[\left(1+ \frac{18B_nKqA_\beta}{\delta}\right)\left(1+ \frac{18B_nKA_\gamma}{\delta}\right)\left(1+ \frac{18B_nKq \sqrt{q}}{a_{\Sigma}\delta}\right)\right]}\Bigg]\nn\\
	\le{}& \sum_{s=1}^S 2^{-s} \Bigg[\sqrt{\ln 2}+ \frac{2^{s}6 \sqrt{2} B_n q K m }{R_n} \sqrt{\ln \left(2p+1\right)} \nn\\
	& \quad +\sqrt{K \ln \left[\left(1+ \frac{2^{s}18B_nKqA_\beta}{R_n}\right)\left(1+ \frac{2^{s}18B_nKA_\gamma}{R_n}\right)\left(1+ \frac{2^{s}18B_nKq \sqrt{q}}{a_{\Sigma}R_n}\right)\right]}\Bigg].
\end{align}
Notice from \eqref{eq_define_deviation8.Gaussian}, that $R_n \geq 2KB_n \max\left( A_\gamma, q A_\beta, \frac{ q\sqrt{q}}{a_{\Sigma}}\right)$. Moreover, it holds that  $1 \le 2^{s+3} $, and $\sum_{s=1}^S2^{-s} = 1-2^{-S} \le 1,\sum_{s=1}^S \left(\sqrt{e}/2\right)^s\le \sqrt{e}/\left(2-\sqrt{e}\right)$, and since for all $s\in \Ns, e^s \geq s$, and thus $2^{-s} \sqrt{s} \le \left(\sqrt{e}/2\right)^s$. Therefore, from \eqref{eq_define_deviation9.Gaussian}:
\begin{align}\label{Extend.CarolineMeynet.eq.5.20.Gaussian}
	&\sum_{s=1}^S 2^{-s}\sqrt{\ln \left[1+ M\left(2^{-s}R_n,F_m,\norm{\cdot}_n\right)\right]}\nn\\
	&\le \sum_{s=1}^S 2^{-s} \Bigg[\sqrt{\ln 2}+ \frac{2^{s}6 \sqrt{2} B_n q K m }{R_n} \sqrt{\ln \left(2p+1\right)}+\sqrt{K \ln \left[\left(2^{s+1}3^2\right)\left(2^{s+1}3^2\right)\left(2^{s+1}3^2\right)\right]}\Bigg]\nn\\
	&= \sum_{s=1}^S 2^{-s}\left[\sqrt{\ln 2}+\frac{2^{s}6 \sqrt{2} B_n q K m }{R_n}\sqrt{\ln(2p+1)}+ \sqrt{K} \sqrt{3\left(\left(s+1\right)\ln 2+2 \ln 3\right)} \right]\nn\\
	&\le\frac{6 \sqrt{2} B_n K q m }{R_n}\sqrt{\ln(2p+1)}S+\sqrt{K}\sqrt{3\ln2} \sum_{s=1}^S 2^{-s}\sqrt{s}+\sqrt{\ln 2}\left(1+\sqrt{3K}\right)+\sqrt{6\ln 3K}\nn\\
	&\le \frac{6 \sqrt{2} B_n K q m }{R_n}\sqrt{\ln(2p+1)}S+\sqrt{K}\sqrt{3\ln2} \sum_{s=1}^S \left(\frac{\sqrt{e}}{2}\right)^s+\sqrt{\ln 2}\left(1+\sqrt{3K}\right)+\sqrt{6\ln 3K}\nn\\
	&\le\frac{6 \sqrt{2} B_nq K m }{R_n}\sqrt{\ln(2p+1)}S+\sqrt{K\ln 2} \underbrace{\left(\frac{\sqrt{3e}}{2-\sqrt{e}}+1+\sqrt{3}+\sqrt{\frac{6\ln 3}{\ln 2}}\right)}_{=:C_1}.
\end{align}
Then, for all $S\in \Ns$, on the event $\cT$:
\begin{align} \label{Extend.CarolineMeynet.eq.5.21.Gaussian}
	&\E{}{\sup_{f_m \in \cF_m} \left|\frac{1}{n} \sum_{i=1}^n \epsilon_i f_m(Z_i)\right|}
	= \E{}{\sup_{f_m \in \cF_m} \left|\frac{1}{n} \sum_{i=1}^n \epsilon_i f_m(Z_i) {\Indi}_{\left|f_m(Z_i)\right| 
			\le R_n}\right|} \left(\text{using Lemma \ref{lem_control_fm_T}}\right)\nn\\
	&= \E{}{\sup_{f \in \cF} \left|\frac{1}{n} \sum_{i=1}^n \epsilon_i f(Z_i) \right|} \le R_n\left(\frac{6}{\sqrt{n}} \sum_{s=1}^S 2^{-s} \sqrt{\ln \left[1+M\left(2^{-s}R_n,\cF,\norm{.}_n\right)\right]}+2^{-S}\right)\left(\text{using Lemma \ref{Extend.CarolineMeynet.Lemma.5.3}}\right)\nn\\
	&\le R_n\left[\frac{6}{\sqrt{n}}\left(\frac{6 \sqrt{2} B_n K m q }{R_n}\sqrt{\ln(2p+1)}S+\sqrt{K\ln 2}C_1\right)+ 2^{-S}\right]\left(\text{using Lemma \ref{Extend.CarolineMeynet.eq.5.20.Gaussian}}\right).
\end{align}
%
We choose $S=\ln n/\ln 2$ so that the two terms depending on $S$ in \eqref{Extend.CarolineMeynet.eq.5.21.Gaussian} are of the same order. In particular, for this value of $S$, $2^{-S} \le 1/n$, and we deduce from \eqref{Extend.CarolineMeynet.eq.5.21.Gaussian} and \eqref{eq_define_deviation8.Gaussian} that on the event $\cT$,
\begin{align*}
	&\E{}{\sup_{f_m \in \cF_m} \left|\frac{1}{n} \sum_{i=1}^n \epsilon_i f_m(Z_i)\right|}\nn\\
	& \le \frac{36\sqrt{2}B_n K mq }{\sqrt{n}}\sqrt{\ln(2p+1)}\frac{\ln n}{\ln 2}+2KB_n \left(A_\gamma+  q A_\beta + \frac{ q\sqrt{q}}{a_{\Sigma}}\right)\left(6\sqrt{\ln 2} C_1\frac{\sqrt{K}}{\sqrt{n}}+ \frac{1}{n}\right)\nn\\
	& \le \frac{B_n K m q }{\sqrt{n}}\sqrt{\ln(2p+1)}\ln n\underbrace{\frac{36\sqrt{2}}{\ln 2}}_{\approx 73.45}+\frac{K\sqrt{K}}{\sqrt{n}}B_n \left(A_\gamma+  q A_\beta + \frac{ q\sqrt{q}}{a_{\Sigma}}\right)\underbrace{2\left(6\sqrt{\ln 2} C_1+ 1\right)}_{\approx 141.32}\nn\\
	& <  \frac{74 KB_n}{\sqrt{n}} \left[m q\sqrt{\ln(2p+1)}\ln n+2\sqrt{K} \left(A_\gamma+  q A_\beta + \frac{ q\sqrt{q}}{a_{\Sigma}}\right)\right].
\end{align*}

\subsection{Proof of Lemma \ref{lem_differentialEntropy_SGaME}}\label{sec_proof_lem_differentialEntropy_SGaME}

Since $\ln(z)$ is concave in $z$, Jensen's inequality
(see e.g. \citet{jensen1906fonctions,cover1999elements}) implies that $\ln\left(\E{Z}{Z}\right) \ge \E{Z}{\ln\left(Z\right)}$, where $Z$ is a random variable. 
Thus, for all $x\in\cX$, Jensen's inequality and Lemma \ref{eq_lem_product_2Gaussians} lead us the following upper bound
\begin{align}
	\int_{\R^q} \ln\left(s_0\left(y|x\right)\right) s_0\left(y|x\right) dy 
	&\le
	\sum_{k=1}^K {g}_{k}\left(x;\gamma_0\right) \ln\left[ \int_{\R^q} s_0\left(y|x\right) \cN\left(y;v_{0k}(x),\Sigma_{0k}\right) dy \right]   \nn\\
	%
	&\le
	\sum_{k=1}^K {g}_{k}\left(x;\gamma_0\right)  \ln\left[ \sum_{l=1}^K  {g}_{l}\left(x;\gamma_0\right) C_{s_0} \right]\nn\\
	&= \ln C_{s_0}  < \infty, \nn\\
	&\quad \text{where $C_{s_0}= \left(4\pi\right)^{-q/2} A_\Sigma^{q/2}, \left(\text{using Lemma \ref{eq_lem_product_2Gaussians}}\right) $.}
\end{align}
Therefore, we obtain 
\begin{align*}
	\max\left\{0,\sup_{x \in \cX} \int_{\R^q} \ln\left(s_0\left(y|x\right)\right) s_0\left(y|x\right) dy  \right\} \le \max\left\{0,\ln C_{s_0} \right\}=:H_{s_0} < \infty.
\end{align*}
%
%
Next, we state the following important Lemma \ref{eq_lem_product_2Gaussians}, which is used in the proof of Lemma \ref{lem_differentialEntropy_SGaME}.
\begin{lemma}\label{eq_lem_product_2Gaussians}
	There exists a positive constant $C_{s_0}:= \left(4\pi\right)^{-q/2} A_\Sigma^{q/2}$, $ 0 < C_{s_0} < \infty$, such that for all $k \in [K], l \in [L]$,
	\begin{align}\label{eq_upperbound_product2Gaussians}
		\int_{\R^q}  \cN\left(y;v_{0l}(x),\Sigma_{0l}\right) \cN\left(y;v_{0k}(x),\Sigma_{0k}\right) dy < C_{s_0}, \quad \forall x \in \cX.
	\end{align}
\end{lemma}
{\bf Proof of Lemma \ref{eq_lem_product_2Gaussians}.}
Firstly, for all $k \in [K], l \in [L]$, given $$c_{lk}(x) = C_{lk} \left[\Sigma_{0l}^{-1}v_{0l}(x)+\Sigma_{0k}^{-1} v_{0k}(x) \right], \quad C_{lk} = \left(\Sigma_{0l}^{-1}+\Sigma_{0k}^{-1}\right)^{-1},$$
Lemma \ref{lem_product_2Gaussians} leads to
\begin{align}\label{eq_applied_lem_product_2Gaussians}
	&\int_{\R^q}  \left[\cN\left(y;v_{0l}(x),\Sigma_{0l}\right) \cN\left(y;v_{0k}(x),\Sigma_{0k}\right) \right]	 dy 
	= Z^{-1}_{lk} \underbrace{\int_{\R^q} \cN\left(y;c_{lk}(x),C_{lk}\right)	 dy}_{=1}, \text{ where }\nn\\
	&= \left(2\pi\right)^{-q/2} \det\left(\Sigma_{0l}+\Sigma_{0k}\right)^{-1/2} \exp\left(-\frac{1}{2} \left(v_{0l}(x)-v_{0k}(x)\right)^\top \left(\Sigma_{0l}+\Sigma_{0k}\right)^{-1}\left(v_{0l}(x)-v_{0k}(x)\right)\right)
\end{align}	
Next, since the determinant is the product of the eigenvalues, counted with multiplicity, and Weyl's inequality, see \eg~Lemma \ref{lem_Weyl_inequality}, for all $k \in [K], l \in [L]$, we have 
\begin{align}
	\det\left(\Sigma_{0l}+\Sigma_{0k}\right) &\ge \left[m\left(\Sigma_{0l}+\Sigma_{0k}\right)\right]^q \nn\\
	&\ge \left[m\left(\Sigma_{0l}\right)+m\left(\Sigma_{0k}\right)\right]^q \left(\text{using \eqref{eq_Weyl_inequality} from Lemma \ref{lem_Weyl_inequality}}\right)\nn\\
	& = \left[M\left(\Sigma_{0l}^{-1}\right)^{-1}+M\left(\Sigma_{0k}^{-1}\right)^{-1}\right]^q\nn\\
	&\ge \left(2A_\Sigma^{-1}\right)^q \left(\text{using boundedness assumptions in \eqref{eq_defBoundedParameters}}\right).\nn
\end{align}
Therefore, for all $k \in [K], l \in [L]$, it holds that 
\begin{align}
	\det\left(\Sigma_{0l}+\Sigma_{0k}\right)^{-1/2}
	&\le 2^{-q/2}\left(A_\Sigma\right)^{q/2} \left(\text{using boundedness assumptions in \eqref{eq_defBoundedParameters}}\right). \label{eq_upperbound_detSumCovariances}
\end{align}
Since $\left(\Sigma_{0l}+\Sigma_{0k}\right)^{-1}$ is a positive definite matrix, it holds that
$$\left(v_{0l}(x)-v_{0k}(x)\right)^\top \left(\Sigma_{0l}+\Sigma_{0k}\right)^{-1}\left(v_{0l}(x)-v_{0k}(x)\right) \ge 0, \quad \forall x \in \cX, l\in[L],k\in[K].$$
Then, since the exponential function is increasing, $\forall x \in \cX, l\in[L],k\in[K],$ we have
\begin{align}\label{eq_upperBound_QuadraticForm}
	\exp\left(-\frac{1}{2} \left(v_{0l}(x)-v_{0k}(x)\right)^\top \left(\Sigma_{0l}+\Sigma_{0k}\right)^{-1}\left(v_{0l}(x)-v_{0k}(x)\right)\right) \le \exp(0) = 1.
\end{align}	
Finally, from \eqref{eq_applied_lem_product_2Gaussians}, \eqref{eq_upperbound_detSumCovariances} and \eqref{eq_upperBound_QuadraticForm}, we obtain
\begin{align*}
	\int_{\R^q}  \left[\cN\left(y;v_{0l}(x),\Sigma_{0l}\right) \cN\left(y;v_{0k}(x),\Sigma_{0k}\right) \right]	 dy 
	&\le \left(2\pi\right)^{-q/2} 2^{-q/2}A_\Sigma^{q/2} = \left(4\pi\right)^{-q/2} A_\Sigma^{q/2} =:C_{s_0} < \infty.
\end{align*}

\section{Further Technical results} \label{technicalResult}

We denote the vector space of all $q$-by-$q$ real matrices by $\R^{q\times q}$ ($q \in \Ns$):
\begin{align*}
	A \in \R^{q \times q} \Longleftrightarrow A =\left(a_{i,j}\right) = \begin{bmatrix}
		a_{1,1}&\cdots&a_{1,q}\\
		\vdots&&\vdots\\
		a_{q,1}&\cdots&a_{q,q}
	\end{bmatrix}, a_{i,j} \in \R.
\end{align*}
If a capital letter is used to denote a matrix (\eg~$A,B$), then the corresponding lower-case letter with subscript $i,j$ refers to the $\left(i,j\right)$th entry (\eg~$a_{i,j},b_{i,j}$). When required, we also designate the elements of a matrix with the notation $\left[A\right]_{i,j}$ or $A\left(i,j\right)$. Denote the $q$-by-$q$ identity and zero matrices by $\id_q$ and $\zero_q$, respectively.

\begin{lemma}[Derivative of quadratic form, \cf~\citet{magnus2019matrix}]\label{derivative.QuadraticPrecision}
	Assume that $X$ and $a$ are non-singular matrix in $\R^{q \times q}$ and vector in $\R^{q \times 1}$, respectively. Then
	\begin{align*}
		\frac{\partial a^\top X^{-1}a}{\partial X} = -X^{-1}aa^\top  X^{-1}.
	\end{align*}
\end{lemma}
\begin{lemma}[Jacobi's formula, \cf~Theorem 8.1 from~\citet{magnus2019matrix}]\label{Jacobi's formula}
	If $X$ is a differentiable map from the real numbers to $q$-by-$q$ matrices,
	\begin{align*}
		\frac{d}{dt} \det \left(X(t)\right)  = \tr\left(\adj \left(X(t)\right) \frac{d X(t)}{dt}\right),
		\frac{\partial \det \left(X\right) }{\partial X}  = \left(\adj \left(X\right)\right)^\top  = \det\left(X\right) \left(X^{-1}\right)^\top .
	\end{align*}
\end{lemma}

\begin{lemma}[Operator induced $p$-norm] \label{def.p-normMatrix}
	We recall an operator (induced) $p$-norms of a matrix $A \in \R^{q \times q}$ ($q \in \Ns, p \in \left\{1,2,\infty\right\}$),
	\begin{align} \label{eq.p-normMatrix}
		\norm{A}_p = \max_{x \neq 0} \frac{\norm{A x}_p}{\norm{x}_p} = \max_{x \neq 0} \left\|A \left(\frac{x}{\norm{x}_p}\right)\right\|_p =  \max_{\norm{x}_p = 1}\norm{A x}_p,
	\end{align}
	where for all $x \in \R^q$,
	\begin{align}
		\norm{x}_\infty \le \norm{x}_1 = \sum_{i=1}^q \left|x_i\right| \le q \norm{x}_\infty,\label{eq.vector.infinity1.equiv}\\
		\norm{x}_2 =\left(\sum_{i=1}^q \left|x_i\right|^2\right)^{\frac{1}{2}} = \left(x^\top x\right)^{\frac{1}{2}} \le \norm{x}_1 \le \sqrt{q}\norm{x}_2,\text{and} \label{eq.vector.12.equiv}\\
		\norm{x}_\infty =\max_{1 \le i \le q} \left|x_i\right| \le \norm{x}_2 \le \sqrt{q} \norm{x}_\infty \label{eq.vector.infinity2.equiv}.
	\end{align}
\end{lemma}

\begin{lemma}[Some matrix $p$-norm properties, \citet{golub2012matrix}]
	By definition, we always have the important property that for every $A \in \R^{q\times q}$ and $x \in \R^q$,
	\begin{align} 
		\label{norm.Equality}\norm{A x}_p \le \norm{A}_p\norm{x}_p,
	\end{align}
	and every induced $p$-norm is submultiplicative, \ie for every $A \in \R^{q\times q} $ and $\textbf{B} \in \R^{q\times q}$,
	\begin{align} 
		\label{submultiplicative.InducedNorm}
		\norm{A B}_p \le \norm{A}_p\norm{B}_p.
	\end{align}
	In particular, it holds that
	\begin{align}
		\norm{A}_1 &= \max_{1 \le j \le q} \sum_{i=1}^q \left|a_{ij}\right| \le \sum_{j=1}^q\sum_{i=1}^q \left|a_{ij}\right| := \left\| \vect(A)\right\|_1 \le q \norm{A}_1, \label{eq.Ineq.1-norm_Vect-Infty}\\
		\norm{\vect(A)}_\infty&:= \max_{1 \le i \le q,1 \le j \le q} \left|a_{ij}\right| \le \norm{A}_\infty = \max_{1 \le j \le q} \sum_{i=1}^q \left|a_{ij}\right| \le  q\norm{\vect(A)}_\infty  ,\label{eq.Ineq.Infty-norm_Vect-Infty}\\
		\norm{\vect(A)}_\infty &\le \norm{A}_2=\lambda_{\max}(A) \le q \norm{\vect(A)}_\infty, \label{eq.Ineq.2-norm_Vect-Infty}
	\end{align}
	where $\lambda_{\max}$ is the largest eigenvalue of a positive definite symmetric matrix $A$. The $p$-norms, when $p\in \left\{1,2,\infty\right\}$, satisfy
	\begin{align}
		\frac{1}{\sqrt{q}} \norm{A}_{\infty} &\le \norm{A}_{2} \le \sqrt{q}\norm{A}_{\infty}\label{eq.Ineq.Infty-2-norm},\\
		\frac{1}{\sqrt{q}} \norm{A}_{1} &\le \norm{A}_{2} \le \sqrt{q}\norm{A}_{1}.\label{eq.Ineq.1-2-norm}
	\end{align}
\end{lemma}

Given $\delta >0$, we need to define the $\delta$-packing number and $\delta$-covering number. 
\begin{definition}[$\delta$-packing number, \cf~Definition 5.4 from~\citet{wainwright2019high}] \label{def.delta-packing.Number}
	Let $\left(\cF,\norm{\cdot}\right)$ be a normed space and let $\cG \subset \cF$. With $\left(g_i\right)_{i=1,\ldots,m} \in \cG$, $\left\{g_1,\ldots,g_m\right\}$ is an $\delta$-packing of $\cG$ of size $m \in \Ns$, if $\norm{g_i - g_j} > \delta, \forall i \ne j, i,j \in \left\{1,\ldots,m\right\}$, or equivalently, $\bigcap_{i=1}^n B\left(g_i,\delta/2\right) = \emptyset$. Upon defining $\delta$-packing, we can measure the maximal number of disjoint closed balls with radius $\delta/2$ that can be ``packed'' into $\cG$. This number is called the \emph{$\delta$-packing number} and is defined as 
	\begin{align}\label{delta-packing.Number}
		M\left(\delta,\cG,\norm{\cdot}\right) := \max\left\{m \in \Ns: \exists\text{$\delta$-packing of $\cG$ of size $m$}\right\}.
	\end{align}
\end{definition}
\begin{definition}[$\delta$-covering number, \cf~Definition 5.1 from~\citet{wainwright2019high}]\label{def.delta-covering.Number}
	Let $\left(\cF,\norm{\cdot}\right)$ be a normed space and let $\cG \subset \cF$.  With $\left(g_i\right)_{i \in [n]} \in \cG$, $\left\{g_1,\ldots,g_n\right\}$ is an $\delta$-covering of $\cG$ of size $n$ if $\cG \subset \cup_{i=1}^n B\left(g_i,\delta\right)$, or equivalently, $\forall g \in \cG, \exists i$ such that $\norm{g-g_i} \le \delta$. Upon defining the $\delta$-covering, we can measure the minimal number of closed balls with radius $\delta$, which is necessary to cover $\cG$. This number is called the \emph{$\delta$-covering number} and is defined as 
	\begin{align}\label{delta-covering.Number}
		N\left(\delta,\cG,\norm{\cdot}\right) := \min\left\{n \in \Ns: \exists\text{$\delta$-covering of $\cG$ of size $n$}\right\}.
	\end{align}
	The covering entropy (metric entropy) is defined as follows $H_{\norm{.}}\left(\delta,\cG\right) = \ln 	\left(N\left(\delta,\cG,\norm{\cdot}\right)\right)$.
\end{definition}
The relation between the packing number and the covering number is described in the following lemma.
\begin{lemma}[Lemma 5.5 from~\citet{wainwright2019high}]\label{lem.delta-packingCovering.Lemma}
	Let $\left(\cF,\norm{\cdot}\right)$ be a normed space and let $\cG \subset \cF$. Then
	\begin{align*}
		M\left(2\delta,\cG,\norm{\cdot}\right) \le	N\left(\delta,\cG,\norm{\cdot}\right) \le M\left(\delta,\cG,\norm{\cdot}\right).
	\end{align*}
\end{lemma}
%

\begin{lemma}[Chernoff's inequality,~\eg Chapter 2 in~\citet{wainwright2019high}] \label{lem.Chernoff.Inequality} Assume that the random variable has a moment generating function in a neighborhood of zero, meaning that there is some constant $b>0$ such that the function $\varphi(\lambda) = \E{}{e^{\lambda(U-\mu)}}$ exists for all $\lambda \le |b|$. In such a case, we may apply Markov's inequality to the random variable $Y = e^{\lambda(U-\mu)}$, thereby obtaining the upper bound
	\begin{align*}
		\bP\left(U-\mu \geq a\right) = \bP\left(e^{\lambda\left(U-\mu\right)} \geq e^{\lambda t}\right) \le \frac{\E{}{e^{\lambda(U-\mu)}}}{e^{\lambda t}}.
	\end{align*}
	Optimizing our choice of $\lambda$ so as to obtain the tightest result yields the Chernoff bound
	\begin{align}
		\ln \left(\bP\left(U-\mu \geq a\right)\right) \le \sup_{\lambda \in [0,b]} \left\{\lambda t - \ln\left(\E{}{e^{\lambda(U-\mu)}}\right) \right\}. \label{eq.Chernoff-bound}
	\end{align}
	In particular, if $U \sim \cN(\mu,\sigma)$ is a Gaussian random variable with mean $\mu$ and variance $\sigma^2$. By a straightforward calculation, we find that $U$ has the moment generating function
	\begin{align*}
		\E{}{e^{\lambda U}} = e^{\mu\lambda + \frac{\sigma^2\lambda^2}{2}}, \text{valid for all }\lambda \in \R.
	\end{align*}
	Substituting this expression into the optimization problem defining the optimized Chernoff bound \eqref{eq.Chernoff-bound}, we obtain
	\begin{align*}
		\sup_{\lambda \geq 0} \left\{\lambda t - \ln\left(\E{}{e^{\lambda(U-\mu)}}\right) \right\} = \sup_{\lambda \geq 0} \left\{\lambda t - \frac{\sigma^2\lambda^2}{2} \right\} = -\frac{t^2}{2 \sigma^2},
	\end{align*}
	where we have taken derivatives in order to find the optimum of this quadratic function. So, \eqref{eq.Chernoff-bound} leads to
	\begin{align}
		\bP\left(X \geq \mu + t \right) \le e^{-\frac{t^2}{2 \sigma^2}}, \text{ for all } t\geq 0.\label{eq.Chernoff-bound.Gaussian}
	\end{align}
\end{lemma}
Recall that a multi-index $\alpha = \left(\alpha_1,\ldots,\alpha_p\right), \alpha_i \in \Ns, \forall i \in \left\{1,\ldots,p\right\}$ is an $p$-tuple of non-negative integers. Let
\begin{align*}
	\left|\alpha\right| = \sum_{i=1}^p \alpha_i, \h \alpha! = \prod_{i=1}^p \alpha_i!,
	x^\alpha = \prod_{i=1}^p x^{\alpha_i}_i, x \in \R^p, \partial^{\alpha} f= \partial_1^{\alpha_1}\partial_2^{\alpha_2}\cdots\partial_p^{\alpha_p} = \frac{\partial^{|\alpha|}f}{\partial x^{\alpha_1}_1\partial x^{\alpha_2}_2\cdots\partial x^{\alpha_p}_p}.
\end{align*}
The number $|\alpha|$ is called the \emph{order} or \emph{degree} of $\alpha$. Thus, the order of $\alpha$ is the same as the order of $x^\alpha$ as a monomial or the order of $\partial^{\alpha}$ as a partial derivative.
\begin{lemma}[Multivariate Taylor's Theorem from~\citet{duistermaat2004multidimensional}] \label{Taylor'sInequality} Suppose $f: \R^p \mapsto \R$ is in the class $C^{k+1}$, of continuously differentiable functions, on an open convex set $S$. If $a \in S$ and $a+h \in S$, then 
	\begin{align*}
		f(a+h) = \sum_{|\alpha| \le k} \frac{\partial^\alpha f(a)}{\alpha!}h^\alpha + R_{a,k}(h),
	\end{align*}
	where the remainder is given in Lagrange's form by
	\begin{align*}
		R_{a,k}(h) = \sum_{|\alpha|= k+1} \partial^\alpha f(a+ch)\frac{h^\alpha}{\alpha!}\text{ for some } c \in (0,1),
	\end{align*}
	or in integral form by
	\begin{align*}
		R_{a,k}(h) = (k+1)\sum_{|\alpha|= k+1} \frac{h^\alpha}{\alpha!} \int_0^1 (1-t)^k \partial^\alpha f(a+th)dt.
	\end{align*}
	In particular, we can estimate the remainder term if $\left|\partial^\alpha f(x)\right| \le M$ for $x \in S$ and $|\alpha|= k+1$,
	\begin{align*}
		\left|R_{a,k}(h)\right| \le \frac{M}{(k+1)!}\left\|h\right\|_1^{k+1}, \left\|h\right\|_1 = \sum_{i=1}^p |h_i|. 
	\end{align*}
\end{lemma}
Recall that the multivariate Gaussian (or Normal) distribution has a joint density given by
\begin{align}
	\cN\left(y;\mu;\Sigma\right) = \left(2 \pi\right)^{-q/2}\det\left(\Sigma\right)^{-1/2}\exp\left(-\frac{1}{2}\left(y - \mu\right)^\top\Sigma^{-1}\left(y-\mu\right)\right),
\end{align}
where $\mu$ is the mean vector (of length $q$) and $\Sigma$ is the symmetric, positive definite covariance matrix (of size $q\times q)$. Then, we have the following well-known Gaussian identity, see more in Lemma \ref{lem_product_2Gaussians}, which is proved in Equation (A.7) from \citet{williams2006gaussian}.
\begin{lemma}[Product of two Gaussians] \label{lem_product_2Gaussians}
	The product of two Gaussians gives another (un-normalized) Gaussian
	\begin{align}\label{eq_Product_2Gaussians}
		\cN\left(y;a,A\right)\cN\left(y;b,B\right) &= Z^{-1}\cN\left(y;c,C\right), \text{where, }\\
		c &= C \left(A^{-1}a+B^{-1}b\right) \text{ and } C = \left(A^{-1}+B^{-1}\right)^{-1},\nn\\
		Z^{-1} &= \left(2\pi\right)^{-q/2} \det\left(A+B\right)^{-1/2} \exp\left(-\frac{1}{2} \left(a-b\right)^\top \left(A+B\right)^{-1}\left(a-b\right)\right) \nn.
	\end{align}
\end{lemma}
We recall the following inequality of Hermann Weyl, see \eg~\citet[Theorem 4.3.1]{horn2012matrix}
\begin{lemma}[Weyl's inequality] \label{lem_Weyl_inequality}
	Let $A,B \in \R^{q\times q}$ be Hermitian and let the respective eigenvalues of $A,B$, and $A+B$ be $\left\{\lambda_i\left(A\right)\right\}_{i\in[q]}$, $\left\{\lambda_i\left(B\right)\right\}_{i\in[q]}$, and $\left\{\lambda_i\left(A+B\right)\right\}_{i\in[q]}$, each  algebraically nondecreasing order as:
	\begin{align*}
		m\left(A\right) = \lambda_1\left(A\right) \le \lambda_2\left(A\right) \le \ldots \le \lambda_q\left(A\right) = M\left(A\right).
	\end{align*}
	Then, for each $i \in [q]$,
	\begin{align*}
		\lambda_i\left(A+B\right) &\le \lambda_{i+j}\left(A\right) + \lambda_{q-j}\left(B\right),\quad j  \in \left\{0\right\} \cup [q-i],
		%
		\lambda_{i-j+1}\left(A\right) + \lambda_{j}\left(B\right) \le \lambda_i\left(A+B\right), \quad j  \in [i].
	\end{align*}
	In particular, we have
	\begin{align}\label{eq_Weyl_inequality}
		M\left(A+B\right) &\le M\left(A\right) +M\left(B\right), 
		%
		m\left(A+B\right) \ge m\left(A\right) +m\left(B\right).
	\end{align}
\end{lemma}


\end{document}